\documentclass{article}

\usepackage[numbers]{natbib}

\usepackage[final]{neurips_2022}

\usepackage[utf8]{inputenc} 
\usepackage[T1]{fontenc}    
\usepackage{hyperref}       
\usepackage{url}            
\usepackage{booktabs}       
\usepackage{amsfonts}       
\usepackage{nicefrac}       
\usepackage{microtype}      
\usepackage{xcolor}         
\usepackage{amsmath}
\usepackage{amssymb}
\usepackage{mathtools}
\usepackage{amsthm}
\usepackage{microtype}
\usepackage{graphicx}
\usepackage{subfigure}
\usepackage[font=small]{caption}
\usepackage{xcolor} 
\usepackage{bbm}
\usepackage[colorinlistoftodos,bordercolor=orange,backgroundcolor=orange!20,linecolor=orange,textsize=scriptsize]{todonotes}
\usepackage{enumitem}
\usepackage{pifont}
\usepackage{wrapfig}

\usepackage{algorithm}
\usepackage{algorithmic}

\allowdisplaybreaks
\newcommand*{\Scale}[2][4]{\scalebox{#1}{$#2$}}%

\def\R{\mathbb{R}}
\newcommand{\E}{{\mathbb E}}

\def\R{\mathbb R}
\def\E{\mathbb E}
\def\EE{\mathbb E}

\def\cD{\mathcal{D}}

\newcommand{\EndProof}{\begin{flushright}$\square$\end{flushright}}

\usepackage[capitalize,noabbrev]{cleveref}

\theoremstyle{plain}
\newtheorem{theorem}{Theorem}[section]

\newtheorem{lemma}[theorem]{Lemma}

\theoremstyle{definition}
\newtheorem{definition}[theorem]{Definition}
\newtheorem{assumption}[theorem]{Assumption}
\theoremstyle{remark}

\newcommand{\add}[1]{{\color{black}#1}}

\title{Decentralized Local Stochastic Extra-Gradient for Variational Inequalities}

\author{%
Aleksandr Beznosikov\\
  Innopolis University\thanks{Research Center for Artificial Intelligence, Innopolis University}, MIPT\thanks{Moscow Institute of Physics and Technology}, HSE University and Yandex\\ 
  \texttt{anbeznosikov@gmail.com}
  \And
  Pavel Dvurechensky \\
  WIAS\thanks{Weierstrass Institute for Applied Analysis and Stochastics}\\ 
  \texttt{pavel.dvurechensky@wias-berlin.de}
  \And
  Anastasia Koloskova \\
  EPFL\\ 
  \texttt{anastasia.koloskova@epfl.ch}
  \And
  Valentin Samokhin \\
  IITP RAS\thanks{Institute for Information Transmission Problems RAS}\\ 
  \texttt{samohin.vyu@phystech.edu}
  \And
  Sebastian U. Stich \\
  CISPA\thanks{CISPA Helmholtz Center for Information Security} \\ 
  \texttt{stich@cispa.de}
  \And
  Alexander Gasnikov \\
  MIPT, HSE University and IITP RAS\\ 
  \texttt{gasnikov@yandex.ru}
}

\begin{document}

\maketitle

\begin{abstract}
We consider distributed stochastic variational inequalities (VIs) on unbounded domains with the problem data that is heterogeneous (non-IID) and distributed across many devices. We make a very general assumption on the computational network that, in particular, covers the settings of fully decentralized calculations with time-varying networks and centralized topologies commonly used in Federated Learning. Moreover, multiple local updates on the workers can be made for reducing the communication frequency between the workers.
We extend the stochastic extragradient method to this very general setting and theoretically analyze its convergence rate in the strongly-monotone, monotone, and non-monotone (when a Minty solution exists) settings. The provided rates explicitly exhibit the dependence on network characteristics (e.g., mixing time), iteration counter, data heterogeneity, variance, number of devices, and other standard parameters. As a special case, our method and analysis apply to distributed stochastic saddle-point problems (SPP), e.g., to the training of Deep Generative Adversarial Networks (GANs) for which decentralized training has been reported to be extremely challenging. In experiments for the decentralized training of GANs we demonstrate the effectiveness of our proposed approach.
\end{abstract}

\section{Introduction}
In large-scale machine learning (ML) scenarios the training data is often split between many devices, such as geographically distributed datacenters or mobile devices~\cite{Kairouz2019:federated}. 
Decentralized training methods can learn an ML model with the same accuracy as if all the data would be aggregated on one single server~\cite{Lian2017:decentralizedSGD,Assran:2018sdggradpush}.
At the same time, training in a decentralized fashion has many advantages over traditional centralized approaches in such core aspects  as data ownership, privacy, fault tolerance, and scalability.
A particular instance of the decentralized learning setting is Federated Learning (FL), where the training is orchestrated by a single device or server that communicates with all the participating client devices~\cite{McMahan16:FedLearning,Kairouz2019:federated}.
In contrast, in fully decentralized learning (FD) scenarios the devices only communicate with their neighbors in the  communication network graph  with possibly  arbitrary topology~\cite{Lian2017:decentralizedSGD}. Thus, decentralized algorithms are important in scenarios where centralized communication is expensive, not desired, or impossible.

There have been tremendous advances recently in the development, design, and understanding of decentralized training schemes~\cite{nedic2009distributed,Wei2012:distributedadmm,%
shi2015extra,Lian2017:decentralizedSGD,%
scaman2017optimal,Uribe:2018uk,%
Tang2018:d2,Wang2018:cooperativeSGD,%
dvurechensky2018decentralize,rogozin2021accelerated,krawtschenko2020distributed%
}.
In particular, such aspects  as data-heterogeneity~\cite{Tang2018:d2,Pu2020:tracking,Lin2021:quasi}, communication efficiency (through local updates~\cite{Lan2018:decentralized,koloskova2020unified} or compression~\cite{Tang2019:squeeze,Koloskova:2019choco}), and personalization~\cite{Vanhaesebrouck2017:personalization,Bellet2018:p2p} have been studied.
However, all these advances were aimed at training with single-criterion loss functions leading to minimization problems, and they do not apply to more general problem classes. 
For example, the training of Generative Adversarial Networks (GANs)~\cite{goodfellow2014generative}
requires simultaneous competing optimization of the generator and the discriminator objectives, i.e.,\ solving a non-convex-non-concave saddle-point problem (SPP).
This problem structure makes GANs notoriously difficult to train even in the single-node setting~\cite{gidel2019variational,chavdarova2019,Chavdarova2020:tamingGANs}, not talking about training over decentralized datasets~\cite{liu2019decentralized,%
Mukherjee2020:decentralizedminmax,%
rogozin2021decentralized}.

Our goal in this paper is solving decentralized stochastic SPPs, and, more generally, decentralized stochastic Minty variational inequalities (MVIs)~\cite{minty62,juditsky2011solving}. 
In a decentralized stochastic MVI, the data is distributed over $M\geq 1$ devices/nodes and each device $m \in [M]$
has access to its local stochastic oracle $F_m(z,\xi_m)$
for the local operator $F_m(z):= \E_{\xi_m \sim \cD_m} F_m(z,\xi_m)$. The data $\xi_m$ in the device $m$ follows an unknown distribution $\cD_m$ that can be different for each device $m \in [M]$. The devices are connected via a communication network forming a graph such that two devices can exchange information if and only if the corresponding nodes are connected by an edge in this graph. The goal is, while respecting the communication constraints, to find cooperatively a point $z^* \in \mathbb{R}^n$ such that, for all $ z \in \mathbb{R}^n$,
\begin{align}
 {\small \textstyle \frac{1}{M}\sum_{m=1}^M \langle \E_{\xi_m \sim \cD_m} F_m(z,\xi_m),z^* - z \rangle \leq 0.}
 \label{eq:mvi}
\end{align}
A special instance of  decentralized stochastic MVIs is the decentralized stochastic SPP with local objectives $  f_m(x,y):= \E_{\xi_m \sim \cD_m} [f_m(x,y,\xi_m)]  $: 
\begin{align}
  \Scale[1]{ \min\limits_{x \in \mathbb{R}^{n_x}} \max\limits_{y \in \mathbb{R}^{n_y}} \left[ f(x,y) := \textstyle \frac{1}{M}\sum_{m=1}^M f_m(x,y) \right].} \label{eq:dec}
\end{align}
The relation to VI can be seen by considering the variable $z=\bigl[ {\tiny \begin{matrix}x \\ y \end{matrix}} \bigr]$ and the gradient field %
$F_m(z) = \bigl[ {\tiny \begin{matrix} \nabla_x f_m(x,y) \\ -\nabla_y f_m(x,y) \end{matrix}} \bigr]$.
In the special case when $f(x,y)$ is convex-concave, the corresponding operator $F(z)=\frac{1}{M}\sum_{m=1}^M \E_{\xi_m}F_m(z,\xi_m)$ is monotone. However, in the context of GANs training, where $x$ and $y$ are the parameters of the generator and the discriminator, respectively, the local losses $f_m(x,y)$ are possibly non-convex-non-concave in $x,y$ and one can not assume the monotonicity of $F$ in general, see also \cite{diakonikolas2021efficient}.

In this paper, we develop a novel algorithm for solving problems \eqref{eq:mvi} and \eqref{eq:dec}. 
Note that the gradient descent-ascent scheme for the problem~\eqref{eq:dec} may diverge even in the simple convex-concave setting with $M=1$ device~\cite{chavdarova2019}. Thus, unlike \cite{liu2019decentralized}, we use extragradient updates~\cite{%
korpelevich1976extragradient,juditsky2011solving,gidel2019variational} as a building block and combine them with a gossip-type communication protocol~\cite{%
Xiao2014:averaging,boyd2006randomized%
} on arbitrary, possibly time-varying, network topologies. 
One of the main challenges due to the communication constraints is a ``network error'' induced by the impossibility of all the devices to reach the exact consensus, i.e., to have exactly the same information about the current iterate of the algorithm. Thus, each device stores a local variable, and only approximate consensus among the devices can be achieved by gossip steps~\cite{Kong2021:cc}.
Unlike other decentralized algorithms~\cite{scaman2017optimal,liu2019decentralized}, our method avoids multiple gossip steps per iteration, which leads to better practical performance and the possibility to work on time-varying networks. Moreover, our method allows for multiple local updates between communication rounds to reduce the communication overhead. This also makes our approach suitable for communication- and privacy-restricted FL or fully decentralized settings~\cite{%
Zinkevich2010:local}.

\textbf{Our contributions.} 
1) Based on extragradient updates, we develop a novel algorithm for distributed stochastic MVIs (and, as a special case, for distributed stochastic SPPs) with heterogeneous data. Our scheme supports a very general communication protocol that covers centralized settings as in Federated Learning, fully decentralized settings, local steps in both the centralized/decentralized settings, and time-varying network topologies.
In particular, we are not aware of earlier works proposing or analyzing extragradient methods with local steps for the fully decentralized setting or decentralized algorithms for stochastic MVIs over time-varying networks.

2) Under the very general communication protocol and in the three settings of MVIs, i.e., with an operator that is strongly-monotone, monotone, or non-monotone under the Minty condition, we prove the convergence of our algorithm and give an explicit dependence of the rates on the problem parameters:\  characteristics of the network (e.g., mixing time), data heterogeneity, the variance of the data, number of devices, and other standard parameters. These theoretical results translate to the corresponding three settings of SPPs (strongly-convex-strongly-concave, convex-concave, non-convex-non-concave under Minty condition). 
All our theoretical results are valid in the important heterogeneous data regime and allow judging in a quantifiable way how different properties, e.g., data heterogeneity, the scale of the noise in the data, and network characteristics, influence the convergence rate of the algorithm. Even for decentralized settings, our results are novel for time-varying graphs and three different settings of monotonicity. See also Table \ref{tab:comparison} that gives more details on our contribution compared to the existing literature. 
The main challenge of our analysis is to deal with the very general assumption about the communication protocol and cope with the errors caused by the stochastic nature and heterogeneity of the data and limited information exchange between the nodes of the communication network.
As a byproduct of independent interest, we analyze the stochastic extragradient method with biased oracle on unbounded domains, which was not done so far in the literature.\looseness=-1

3) We verify our theoretical results in numerical experiments and demonstrate the practical effectiveness of the proposed scheme. In particular, we train the DCGAN \cite{radford2015unsupervised} architecture on   the CIFAR-10~\cite{cifar10}  dataset. \looseness=-1

\subsection{Related Work}
The research on MVIs dates back at least to 1962~\cite{minty62} with the classical book~\cite{Kinderlehrer} and the recent works 
\cite{liu2019decentralizedprox,Lin2020:nearoptimal,Bullins2020:higherorder,diakonikolas2021efficient}.
VIs arise in a broad variety of applications: image denoising~\cite{esser2010general,chambolle2011first},
game theory and optimal control~\cite{facchinei2007finite}, robust optimization~\cite{BenTal2009:book}, and non-smooth oprimization via smooth reformulations~\cite{nesterov2005smooth,nemirovski2004prox}.
In ML, MVIs and SPPs
 arise in GANs training ~\cite{Daskalakis2018:gan,chavdarova2019,Chavdarova2020:tamingGANs}, reinforcement learning~\cite{Omidshafiei2017:rl,Jin2020:mdp}, and adversarial training~\cite{Madry2017:adv}.\looseness=-1

\textbf{Extragradient.}
The extragradient method (EGM) was first proposed in~\cite{korpelevich1976extragradient}, generalized as the mirror-prox method for deterministic problems in~\cite{nemirovski2004prox} and for stochastic problems with bounded variance in~\cite{juditsky2011solving}. 
Yet, if the stochastic noise is not uniformly bounded, the EGM may diverge, see~\cite{chavdarova2019,Mishchenko2019:extragradient}.

\begin{table*}[ht!]
\resizebox{\linewidth}{!}{
\begin{tabular}{llccccccc}
\hline
Reference & base method & arbitrary network & time-varying & local updates & no multiple gossip steps  & SM & M & NM \\
\hline
Liu et al. 2019 \cite{liu2019decentralized} & Stoch. ES & \ding{52}& \ding{56}& \ding{56}& \ding{56} & \ding{56} & \ding{56}& \ding{52}$^\dagger$\\
Beznosikov et al. 2021\cite{beznosikov2021distributed} Alg. 2 & Stoch. ES & \ding{52}& \ding{56}& \ding{56}& \ding{56} & \ding{52} & \ding{52}& \ding{56}\\
Barazandeh et el. 2021 \cite{BARAZANDEH2021108245} & Stoch. ES & \ding{52}& \ding{56}& \ding{56}& \ding{56} & \ding{56} & \ding{56}& \ding{52}\\ 
Liu et al. 2019 \cite{liu2019decentralizedprox} & Deter. prox & \ding{52}& \ding{56}& \ding{56}& \ding{52} & \ding{56} & \ding{56}& \ding{52}\\
Mukherjee and Chakraborty 2020 \cite{Mukherjee2020:decentralizedminmax} &  Deter. ES & \ding{52}& \ding{56}& \ding{56}& \ding{52} & \ding{52} & \ding{56}& \ding{56}\\
Tsaknakis et al. 2020 \cite{9054056} &  Stoch. DA & \ding{52}& \ding{56}& \ding{56}& \ding{52} & \ding{56} & \ding{56}& \ding{52}$^\ddagger$\\
Rogozin et al. 2021\cite{rogozin2021decentralized} & Deter. ES & \ding{52}& \ding{56}& \ding{56}& \ding{52} & \ding{56} & \ding{52}& \ding{56}\\
Xian et al. 2021 \cite{NEURIPS2021_d994e372} & Stoch. DA & \ding{52}& \ding{56}& \ding{56}& \ding{52} & \ding{56} & \ding{56}& \ding{52}$^\ddagger$\\
Beznosikov et al. 2021 \cite{beznosikov2021distributed} Alg 3 & Stoch. ES & \ding{56}& \ding{56}& \ding{52}& -$^\mathsection$ & \ding{52} & \ding{56}& \ding{56}\\
Deng and Mahdavi 2021\cite{deng2021local}  & Stoch. DA & \ding{56}& \ding{56}& \ding{52}& - & \ding{52} & \ding{56}& \ding{52}$^\ddagger$\\
Hou et al. 2021 \cite{hou2021efficient} & Stoch. DA & \ding{56}& \ding{56}& \ding{52}& - & \ding{52} & \ding{56}& \ding{56}\\ \hline
Ours & Stoch. ES & \ding{52} & \ding{52} & \ding{52}& \ding{52}& \ding{52}& \ding{52}& \ding{52} \\
\hline
\end{tabular} 
}
\hspace*{2mm}\scriptsize{$^\dagger$ -- homogeneous case, $^\ddagger$ -- non-convex-concave SPP (other works use minty condition -- \eqref{as6}), $\mathsection$ -- this column does not apply to centralized algorithms.}
\caption{ Comparison of approaches for distributed strongly-monotone (SM), monotone (M), and non-monotone (NM) VIs or, respectively, strongly-convex-strongly-concave, convex-concave, non-convex-non-concave SPPs.
\\ Definitions of columns: 
\textbf{base method} --- the non-distributed algorithm that is taken as the basis for the distributed method, typically it is either the extragradient method (EGM) or the descent-ascent (DA); 
\textbf{arbitrary network} --- supporting fully decentralized vs. only centralized topology; 
\textbf{time-varying} --- decentralized method supporting time-varying network topology; 
\textbf{local updates} --- method supporting local steps between communications; 
\textbf{no multiple gossip steps} --- at one global iteration the method does not use many iterations of gossip averaging to reach a good consensus accuracy;
\textbf{SM, M, NM} --- monotonicity assumption, see  Assumption~\ref{a2}.}
\label{tab:comparison}
\end{table*}

\textbf{Decentralized algorithms for MVIs and SPPs} are the most closely related to our work. In Table~\ref{tab:comparison}, we summarize their features  and make a comparison with our algorithm, showing that, e.g., existing methods do not support arbitrary time-varying network typologies. 
The methods that use multiple rounds of gossip averaging (sparse communication) per iteration \cite{liu2019decentralized,beznosikov2021distributed,BARAZANDEH2021108245} can give near-optimal theoretical rates, but are often unstable in practice. Thus, it is preferable to have only one sparse communication per iteration~\cite{liu2019decentralizedprox,Mukherjee2020:decentralizedminmax,9054056, rogozin2021decentralized, NEURIPS2021_d994e372}.
The second column of the table refers to standard algorithms that are extended to distributed settings in the corresponding work. In particular, the algorithm of \cite{liu2019decentralizedprox} requires expensive proximal updates.
The closest work to ours is \cite{beznosikov2021distributed}, 
where a decentralized EGM without local steps is analyzed in the (strongly-)monotone setting. Unlike our more general algorithm with local steps, theirs require multiple gossip updates in each iteration which is not desired in practice. 
For the FL, i.e., centralized, setting, \cite{beznosikov2021distributed} studies the EGM with local steps in the strongly-monotone setting, and~\cite{deng2021local,hou2021efficient} study the descent-ascent method with local steps. Yet, all three works do not consider arbitrary time-varying graphs as in our work.

\section{Algorithm}

\add{In this section, we present and discuss the proposed algorithm (Algorithm~\ref{alg4}) that is based on two main ideas: (i) the extragradient step, as in the classical methods for VIs \cite{korpelevich1976extragradient,nemirovski2004prox}, and (ii) the gossip averaging \cite{boyd2006randomized,nedic2009distributed} widely used in decentralized optimization methods and in the literature on diffusion strategies in distributed learning \cite{SAYED2014323, sayed2013diffusion, 8272358, alghunaim2022unified, 9548800}. Unlike these papers that propose algorithms for optimization problems by exploiting gradient descent, our algorithm is based on the extragradient method and is designed to solve VIs and SPPs. Moreover, unlike the mentioned works, our method also allows for local steps in-between the communication rounds and for time-varying networks and has non-asymptotic theoretical convergence rate guarantees.}



Each step of Algorithm~\ref{alg4} can be divided into two phases. The local phase (lines 4--6) consists of a step of the stochastic extragradient method at each node using only local information. As in the non-distributed case, the nodes first make  an extrapolation step ``to look into the future'' and then an update based on the operator value at the ``future'' point. This is followed by the communication phase (gossip step) (line 7), during which the nodes share and average local iterates with their neighbors $\mathcal{N}^k_m$ in the communication network graph corresponding to the iteration $k$. The averaging process involves the weights $w^k_{m,i}$ which are the elements of the matrix $W^k$ 
called the mixing matrix:
\begin{definition}[Mixing matrix] \label{mix} 
We call a matrix $W \in [0;1]^{M \times M}$ a mixing matrix if it satisfies the following conditions: 1) $W$ is symmetric, 2) $W$ is doubly stochastic ($W \mathbf{1} = \mathbf{1}$, $\mathbf{1}^T W = \mathbf{1}^T$, where $\mathbf{1}$ denotes the vector of all ones), 3) $W$ is aligned with the network: $w_{ij} \neq 0$ if and only if $i=j$ or the edge $(i,j)$ is in the communication network graph.
\end{definition}
Reasonable choices of mixing matrices 
 are, for example,
(i) $W^{k} = I_M - \smash{\frac{L^k}{\lambda_{\max} (L^k)}}$, where $L^k$ is the Laplacian matrix of the network graph at the step $k$ and $I_M$ is the identity matrix, or
(ii) using some local rules in the graph, based on the degrees of the neighboring nodes~\cite{Xiao2014:averaging}.
Note that our setting has a great flexibility since in-between the 
iterations the topology of the communication graph is allowed to change, and the matrix $W^k$, that encodes the current structure of the network, changes accordingly. 
This is encoded in line 2, where the matrix $W^k$ is generated by some rule $\mathcal{W}^k$ which can have different nature. Examples include deterministic choice of a sequence of matrices $W^k$, sampling from a time-varying probability distribution on matrices. Even local steps without communication can be encoded with a diagonal matrix $W^k$.

\begin{algorithm}[h]
	\caption{Extra Step Time-Varying Gossip Method}
	\label{alg4}
	{\bf parameters:} stepsize $\gamma >0$, $\{\mathcal{W}^k\}_{k\geq 0}$ -- rules or distributions for mixing matrix in iteration $k$.\\
	{\bf initialize:} $z^0 \in\mathcal{Z}, \forall m :z^0_m=z^0$
	\begin{algorithmic}[1]
\FOR {$k=0,1, 2, \ldots$ } 
\STATE Sample matrix $W^k$ from $\mathcal{W}^k$ \label{line1}
\FOR {each node $m$}
\STATE Generate independently $\xi^k_m \sim \cD_k$, $\xi^{k+1/3}_m \sim \cD_k$  \label{line2}
 \STATE $z^{k+1/3}_m = z^k_m - \gamma F_m(z^k_m, \xi^k_m)$ \label{line3}
 \STATE $z^{k+2/3}_m = z^k_m - \gamma F_m(z^{k+1/3}_m, \xi^{k+1/3}_m)$   \label{line4}
 \STATE $z^{k+1}_m = \sum_{i \in \mathcal{N}^k_m} w^k_{m,i} z_i^{k+2/3}$ \label{line5}
 \ENDFOR
\ENDFOR
	\end{algorithmic}
\end{algorithm}
To ensure that it is possible to approach the consensus between the nodes, we need the following assumption  on the mixing properties of the matrix sequence $W^k$. \looseness=-1
\begin{assumption}[Expected Consensus Rate] \label{a5} 
We assume that there exist a constant $p \in (0, 1]$ and an integer $\tau \geq 1$ such that, after $K$ iterations, for all matrices $Z \in \mathbb{R}^{d\times M}$  and all integers $l \in \{0, \ldots, K/\tau\}$,
    \begin{align} 
    \label{as5} 
    \mathbb{E}_W[\| Z W_{l,\tau} - \bar Z \|_{F}^2] \leq (1 - p) \|Z - \bar Z \|^2_F, 
    \end{align}
    where $W_{l,\tau} = W^{l \tau} \cdot \ldots \cdot W^{(l+1)\tau - 1}$, we use the matrix notation $Z = [z_1, \ldots, z_M]$, $\bar Z = [\bar z, \ldots, \bar z]$ with $\bar z = \frac{1}{M} \sum_{m=1}^M z_m$, and the expectation $\mathbb{E}_W$ is taken over distributions of $W^{t}$ and indices $t \in \{l \tau,...,(l+1)\tau - 1\}$. 
\end{assumption}

This assumption ensures that, after $\tau$ gossip steps with such time-varying matrices, we improve the consensus between the nodes, i.e., how close each $z_m$ is to $\bar{z}$, by the factor of $\frac{1}{1-p}$. Importantly, in this case, some matrices $W^k$ can be, for example, the identity matrix (which corresponds to performing only local  steps  in iteration $k$). 

Assumption \ref{a5} has been recently quite popular in the literature on distributed optimization methods \cite{nedich2016geometrically, koloskova2020unified, kovalev2021lower}. Moreover, it is very general and covers many special cases of decentralized and centralized algorithms.
For example, if we fix $W^{k} = W$ for some fixed connected graph, we get a decentralized algorithm on this graph. If, at the same time, we set the matrix $W = \frac{1}{M} \mathbf{1} \mathbf{1}^T$, then it is easy to see that we get an analog of the centralized setting with the averaging over all nodes performed in each communication step. If we take  $W^{k} = W$ for some fixed connected graph at every $\tau$-th step  and in other steps use $W^k = I_M$, we have a decentralized (and, in particular, centralized) algorithm with local steps \cite{stich2018local, gorbunov2021local,  koloskova2020unified} and communications after each $\tau$ iterations. 
Generic Assumption \ref{a5} covers also many other settings of time-varying decentralized topologies, e.g., random topologies, cliques, $B$-connected graphs \cite{jadbabaie2003coordination,nedic2009distributedAveraging}. Below we show that, under an appropriate choice of the stepsize, our extragradient method provably converges under such a general assumption that covers centralized and decentralized settings, local steps in both centralized and decentralized settings, and changing topologies of the communication graph. Even for decentralized settings, this is novel for time-varying graphs and three different settings of monotonicity which we consider.

\section{Setting and Assumptions}

In this section, we introduce necessary assumptions that are used to analyze the proposed algorithm. 
\begin{assumption}[Lipschitzness] \label{a1} 
For all $m$, the operator $F_m(z)$ is Lipschitz with constant $L$, i.e., 
    \begin{align}
    \label{as1}
    \|F_m(z_1) - F_m(z_2)\| \leq L\|z_1-z_2\|, \quad \forall z_1, z_2. \tag{L}
    \end{align}
\end{assumption}
This is a standard assumption that is used in the analysis of all the methods displayed in Table \ref{tab:comparison}. 
%
\begin{assumption} We consider three scenarios for the operator $F$, namely, when $F$ is strongly-monotone, monotone and non-monotone, but with an additional assumption:    \label{a2} \\
\textbf{(SM) Strong monotonicity.} There exists $\mu > 0$ such that, for all $z_1, z_2$,
\begin{align}
    \label{as2}
    \langle F(z_1) - F(z_2), z_1 - z_2 \rangle \geq \mu\|z_1-z_2\|^2. \tag{SM}
\end{align}
\textbf{(M) Monotonicity.} For all $z_1, z_2 $, it holds that:
\begin{align}
    \label{as21}
    \langle F(z_1) - F(z_2), z_1 - z_2 \rangle \geq 0. \tag{M}
\end{align}
\textbf{(NM) Non-monotonicity (Minty).} There exists $z^*$ such that, for all $z$,
\begin{align}
    \label{as6}
    \langle F(z), z - z^* \rangle\geq 0. \tag{NM}
\end{align}
\end{assumption}
Assumptions \eqref{as2}, \eqref{as21} and \eqref{as1} are standard and classical assumptions in the literature on VIs.  Assumption \eqref{as6} is sometimes called the \textit{Minty or Variational Stability condition} and it has been widely used recently by the community as a structured variant of non-monotonicity \cite{dang2015convergence,iusem2017extragradient, mertikopoulos2018optimistic, liu2019decentralizedprox, kannan2019optimal, hsieh2020explore, diakonikolas2021efficient}, particularly, since it is appropriate in GANs training \cite{liu2019towards, liu2019decentralized, dou2021one, BARAZANDEH2021108245}.

The next assumption is standard for the stochastic setting.
\begin{assumption}[Bounded noise] \label{a3}  $F_m (z, \xi)$ is unbiased and has bounded variance, i.e.,  for all $z$,
    \begin{align} 
    \label{as3}
    \hspace{-1em}\mathbb{E}[F_m(z,\xi)] = F_m(z),\;\mathbb{E}[\|F_m(z,\xi)-F_m(z)\|^2]\leq \sigma^2.
    \end{align}
\end{assumption}

Our last assumption reflects the variability of the local operators compared to their mean and is usually called $D$-heterogeneity. This assumption is widely used in the analysis of local-steps (and not only) algorithms for minimization problems \cite{khaled2020tighter,woodworth2020minibatch,gorbunov2021local, pmlr-v32-shamir14, arjevani2015communication,agafonov2021accelerated, hendrikx2020statistically,daneshmand2021newton,dvurechensky2022hyperfast}. Moreover, \cite{deng2021local,hou2021efficient} use this assumption for the analysis of centralized local-steps methods for SPPs. The authors of \cite{liu2019decentralized} assume $D = 0$ for the decentralized training of GANs. Even in this case algorithms' analysis can be challenging.
\begin{assumption}[$D$-heterogeneity] \label{a4}  The values of the local operator have bounded variablility,  i.e., for all $z$,
    \begin{align} 
    \label{as4}
    \|F_m(z) - F(z) \|^2 \leq D^2.
    \end{align}
\end{assumption}
\section{Main Results}
In this section, we present the convergence rate results for the proposed method under different settings of Assumption \ref{a2}. To present the main result, we introduce notation $\bar z^k := \frac{1}{M} \sum_{m=1}^M z^{k}_m$, $\bar z^{k+1/3} := \frac{1}{M} \sum_{m=1}^M z^{k+1/3}_m$ for the averaged among the devices iterates and $\widehat{z}^{k}=\frac{1}{k+1}\sum_{i=0}^k\bar z^{i+1/3}$ for the averaged among the devices and iterates sequence, a.k.a.\ ergodic average. Finally, we denote $\Delta=\frac{  \tau}{ p}  \bigl( \frac{D^2 \tau}{p} + \sigma^2 \bigr)$ which plays the role of the consensus error, i.e., the error caused by the impossibility of reaching the exact consensus between the nodes. Note that the data heterogeneity appears in the convergence rates only through the quantity $\Delta$.
\begin{theorem}[Main theorem] \label{app:th1}
Let Assumptions \ref{a5}, \ref{a1}, \ref{a3}, \ref{a4} hold and the sequences $\bar z^k$, $\widehat{z}^{k}$ be generated by Algorithm \ref{alg4} that is run for $K > 0$ iterations. Then,

\noindent  $\bullet$ \textbf{Strongly-monotone case:} under Assumption \ref{a2}(SM), with $\gamma = \mathcal{\tilde O}\left( \min\left\{\frac{p}{\tau L}, \frac{1}{\mu K}\right\}\right)$ it holds that
\begin{align}
\label{t_sm}
 \textstyle{\E\left[\|\bar z^{K+1} - z^* \|^2\right] = \Scale[1]{ \mathcal{\tilde O} \left( \|z^{0} - z^* \|^2 \cdot \exp\left( - \frac{\mu K p }{240 L \tau}\right) + \frac{\sigma^2}{\mu^2 M K}  + \frac{L^2 \Delta}{\mu^4   K^2} \right);}}
  \end{align}
 $\bullet$ \textbf{Monotone case:} under Assumption \ref{a2}(M), for any convex compact $\mathcal{C}$ s.t. $z^0,z^* \in \mathcal{C}$ and $\max_{z,z' \in \mathcal{C}} \|z-z'\|\leq \Omega_{\mathcal{C}}$, with $\gamma =  \mathcal{O}\left(\min\left\{\frac{1}{L},\left(\frac{\Omega_{\mathcal{C}}^2M}{K\sigma^2}\right)^{\frac{1}{2}}, \left(\frac{\Omega_{\mathcal{C}}^2}{K^2L^2\Delta}\right)^{\frac{1}{4}}\right\}\right)$ it holds that 
\begin{align}
\label{t_m_1}
 \textstyle{ \sup_{z \in \mathcal{C}} \E\left[\left\langle  F(z), \widehat{z}^{K} - z \right\rangle\right] = \Scale[1]{
\mathcal{O} \left( \frac{L\Omega_{\mathcal{C}}^2}{K} + \frac{\sigma\Omega_{\mathcal{C}}}{\sqrt{MK}}+\frac{\sqrt{L\Omega_{\mathcal{C}}^3\sqrt{\Delta}}}{\sqrt{K}}+\sqrt{\frac{(\Delta+L^2\Omega_{\mathcal{C}}^2)\Omega_{\mathcal{C}} \sqrt{\Delta}}{KL}}\right)}.}
 \end{align}
Under the additional assumption that, for all $k$,  $\|\bar z^k\| \leq \Omega$, with $\gamma = \mathcal{O}\left( \min\left\{\frac{1}{L},\left(\frac{\Omega_{\mathcal{C}}^2M}{K\sigma^2}\right)^{\frac{1}{2}}, \left(\frac{\Omega_{\mathcal{C}}^2}{K^2L^2\Delta}\right)^{\frac{1}{4}},\left(\frac{\Omega_{\mathcal{C}}^2}{K((\Omega +\Omega_{\mathcal{C}})L \sqrt{\Delta} + \Delta)}\right)^{\frac{1}{2}}\right\}\right)$, we have that 
\begin{align}
\label{t_m_2}
\textstyle{ \sup_{z \in \mathcal{C}} \E\left[\left\langle  F(z), \widehat{z}^{K} - z \right\rangle\right]=\Scale[1]{
\mathcal{O} \left( \frac{L\Omega_{\mathcal{C}}^2}{K} + \frac{\sigma\Omega_{\mathcal{C}}}{\sqrt{MK}}+\frac{\sqrt{L\Omega_{\mathcal{C}}^3\sqrt{\Delta}}}{K^{3/4}}+\sqrt{\frac{((\Omega +\Omega_{\mathcal{C}})L \sqrt{\Delta} + \Delta)\Omega_{\mathcal{C}}^2 }{K}}\right);}}
 \end{align}
$\bullet$ \textbf{Non-monotone case:} under Assumption \ref{a2}(NM) and if $\|z^0\| \leq \Omega,\|z^*\| \leq \Omega$, with $\gamma = \mathcal{O}\left(\min\left\{\frac{1}{L}, \left(\frac{\Omega^2}{K^2L^2 \Delta}\right)^{\frac{1}{4}}\right\}\right)$:
\begin{align}
\label{t_nm_1}
 \textstyle{\E\left[ \frac{1}{K+1}\sum_{k=0}^K\|F(\bar z^k)\|^2\right] = \Scale[1]{ \mathcal{O} \left(\frac{L^2\Omega^2}{K}+\frac{\sigma^2}{M}+L\Omega\sqrt{\Delta}+\frac{\sqrt{L\Omega\Delta^{3/4}}}{\sqrt{K}}\right).}}
\end{align} 
Under the additional assumption that, for all $k$, $\|\bar z^k\| \leq \Omega$, with $\gamma = \mathcal{O}\left(\min\left\{\frac{1}{L}, \left(\frac{\Omega^2}{KL \Delta}\right)^{\frac{1}{3}}\right\}\right)$ ,  we have that 
\begin{align}
\label{t_nm_2}
\textstyle{\E\left[ \frac{1}{K+1}\sum_{k=0}^K\|F(\bar z^k)\|^2\right]= \Scale[1.1]{  \mathcal{O} \left(\frac{L^2 \Omega^2}{K} + \frac{\sigma^2}{M} + \frac{(L \Omega \Delta)^{2/3}}{K^{1/3}}   + L \Omega \sqrt{\Delta} \right).}}
\end{align}
\end{theorem}
The proof of the theorem is given in the supplementary material, where one can also find explicit dependence of the rates on the stepsize $\gamma$ before it is chosen optimally. 
We underline that the standard analysis \cite{juditsky2011solving} does not apply for the following reasons. Firstly, unlike \cite{juditsky2011solving}, in our problem \eqref{eq:mvi}, the feasible set is not bounded, which is especially important for the analysis in the monotone and non-monotone settings. Secondly, our algorithm has an additional communication step (line 7) between the computational nodes, which leads to the impossibility for all the nodes to have the same information about the global operator $F(z)$ and about the current iterate $z$. This, in order, leads to a biased oracle that, unlike existing works, has to be analyzed in the setting of an unbounded feasible set, which is quite challenging. 
To analyze our variant of the extragradient method, we successfully handle this challenge.  Our key steps are to bound the bias (see, e.g.,\ the last two terms in the r.h.s.\ of Lemma \ref{Lm:bar_g_lower} that are caused by the network errors), prove the boundedness in expectation of the sequence of the iterates for monotone (see Section \ref{S:app_monot_unbound} of the supplementary material) and non-monotone (see Section \ref{S:app_nonmonot_unbound} of the supplementary material) cases, which may be of independent interest and which we have not seen in the literature, even in the non-distributed setting with biased stochastic oracles. Proving the boundedness is challenging due to the noise caused by the stochasticity and heterogeneity of the data and network effects due to the imperfect exchange of information.
Surprisingly, in the end, we still manage to analyze our algorithm under the very general Assumption \ref{a5} and we are not aware of any results with similar generality of the settings: different network topologies (including time-varying), distributed architectures, different monotonicity assumptions.

The provided convergence rates have an explicit dependence on the problem parameters: the network that is characterized by the mixing time $\tau$ and the mixing factor $p$, the data heterogeneity $D$ (these three quantities appear in the convergence rates only through the quantity $\Delta$), the variance $\sigma^2$ of the noise in the data, the Lipschitz constant $L$, the strong monotonicity parameter $\mu$, the number of nodes/devices $M$. Thus, our rates allow judging how different properties, e.g., data heterogeneity, noise level, and network characteristics influence the convergence rates. This, in particular, opens up an opportunity for a meta-optimization process if we can design the network and change $M$, $\tau$, $p$ to achieve faster convergence. 

We now discuss the convergence results obtained in the theorem, and also compare them with already existing algorithms (see Table \ref{tab:comparison}) and their guarantees. Firstly, all the estimates have a similar several-term structure. The first term corresponds to the deterministic setting and is similar to existing methods for smooth VIs in the non-distributed setting. Only in the strongly-convex case, there is an additional factor $\tau/p$ that increases the condition number $L/\mu$ of the problem. The second (stochastic) term is also standard for the non-distributed setting and corresponds to the stochastic nature of the problem. Note that, for a very general distributed setting, we have managed to obtain the corresponding terms similar to the non-distributed setting.
Moreover, we can see the benefit of exploiting distributed computations: the leading stochastic term depends on $\sigma^2/M$ that decreases as the number $M$ of the nodes increases.
The other terms correspond to the consensus error $\Delta$ and are due to the imperfect communications between the nodes, i.e., that all the nodes can't have exactly the same information about the current iterate. Importantly, in all the cases, this error does not make the overall convergence worse since the dependence on $K$ is no worse in these terms than the dependence on $K$ in the stochastic term.
In the experimental section, we illustrate that the network error is not an artifact of the analysis but  is indeed present in practice.

Theorem \ref{app:th1} is formulated for a fixed budget of iterations $K$ and the corresponding stepsizes $\gamma$ that depend on $K$, which is pretty standard in the literature \cite{juditsky2011solving, stich2019unified, beznosikov2022stochastic}, where many algorithms fix the stepsize depending on the budget of the iterations. In Section \ref{sec:restart} of the supplementary material, we present a simple restarting procedure that allows to extend the results of Theorem \ref{app:th1} to any-time convergence without a-priori fixing $K$. The idea is to set $K_t=2^t$ for $t=0,1,\ldots$ and restart the algorithm after each $K_t$ iterations. We next make refined comments for each particular setting of monotonicity.


$\bullet$ \textbf{Strongly-monotone case:} 
In the centralized setting with local updates, our rate is slightly better than in \cite{beznosikov2021distributed}. 
Unlike our algorithm, centralized algorithms with local steps for SPPs in \cite{deng2021local,hou2021efficient} are based on the gradient descent-ascent method that may diverge in the stochastic setting even for bilinear problems.
Moreover, their analysis implies a very small stepsize $\gamma \sim \frac{\mu p}{L^2\tau}$ (cf. ours $\gamma \sim \frac{p}{L\tau}$), which greatly slows down the convergence of the algorithm. 

For the decentralized setting, \cite{beznosikov2021distributed} propose an optimal algorithm with the rate matching the lower bound which they also give. Our rate is worse probably because of the generality of the Assumption \ref{a5}. On the other hand, our algorithm is more practical since it avoids using multiple gossip steps at each iteration.
Also, our algorithm is more general, allowing us to work with time-varying topologies and local steps even in the decentralized setting.

$\bullet$ \textbf{Monotone case:} The quantity $\sup_{z \in \mathcal{C}} \E\left[\left\langle  F(z), \widehat{z}^{K} - z \right\rangle\right]$ in the convergence rates estimates reflects the stochastic nature of the problem and is a counterpart of the standard restricted gap (or merit) function \cite{nesterov2007dual}: $\text{Gap}_{\mathcal{C}} (u) := \sup_{z \in \mathcal{C}} \left[ \langle F(z),  u - z  \rangle \right]$. When $F$ is a monotone operator, if $\text{Gap}_{\mathcal{C}} (\hat{u}) = 0 $ and $\mathcal{C}$ contains a neighborhood of $\hat{u}$, then \cite{nesterov2007dual,antonakopoulos2019adaptive} $\hat{u}$ is a solution to \eqref{eq:mvi} and even more: it is a strong solution to the corresponding variational inequality, i.e., for all $z$, $\langle F(\hat{u}),\hat{u} - z \rangle \leq 0$.  
Thus, $\text{Gap}_{\mathcal{C}} (u)$ is an appropriate measure of suboptimality in this setting and \eqref{t_m_1} guarantees that after a sufficient number of iterations, we obtain an approximate solution in expectation. Importantly, for \eqref{t_m_1}, neither $z$ nor $\bar z^k$ are assumed to be bounded. As in the previous works on non-distributed algorithms for MVIs \cite{nesterov2007dual,antonakopoulos2019adaptive}, we use $\text{Gap}_{\mathcal{C}} (u)$ with an arbitrary compact set $\mathcal{C}$ that contains $z^0$ and $z^*$  (this can be a large set). Further, \eqref{t_m_2} is a refined version of the general result \eqref{t_m_1} under the additional assumption of the boundedness of the averaged iterates. If the boundedness does not hold, we still have \eqref{t_m_1}. 
Moreover, \eqref{t_m_1} and \eqref{t_m_2} hold for the same method, and to run the algorithm, there is no need to know in advance whether the generated sequence is bounded or not.

Only \cite{beznosikov2021distributed,rogozin2021decentralized} consider MVIs with monotone operator in distributed setting. Our algorithm is more general than theirs: our algorithm supports time-varying networks and local steps between communications. The algorithm in \cite{beznosikov2021distributed} uses multiple gossip steps between the updates of the iterates. On the one hand, this allows decreasing the consensus error $\Delta$. On the other hand, this leads to an additional factor in the number of communications compared to our estimates: the first term in their bound is $\sqrt{\chi}$ times larger than ours, where $\chi > 1$ is some condition number of the mixing matrix. Moreover, multiple gossip steps may be impractical if the communication is performed through unstable channels or is expensive for some reason. The paper \cite{rogozin2021decentralized} considers only deterministic setting.

$\bullet$ \textbf{Non-monotone case:} The same as in the previous case remark on the boundedness of $\bar z^k$, $z^*$ assumed to obtain \eqref{t_nm_2} applies in this case.
Further, in this setting, the convergence is guaranteed up to some accuracy that is governed by the stochastic nature of the problem (the $\sigma^2$-term) and by the distributed nature of the problem (the $\Delta$-terms). With this respect, the results are similar to the non-distributed stochastic extragradient method \cite{barazandeh2021solving} and the distributed method \cite{liu2019decentralized} analyzed in the homogeneous case $D=0$. To the best of our knowledge, convergence up to arbitrarily small accuracy can be guaranteed only for deterministic distributed methods \cite{liu2019decentralizedprox}, i.e., in a much simpler setting than ours. Moreover, the methods of \cite{liu2019decentralizedprox} are not the most robust since they require evaluating the proximal operator of a function and it is assumed that this can be done in a  closed form, which is computationally expensive and may not hold in practice.

Note that, based on our result, it is possible to achieve convergence up to arbitrarily small accuracy if one considers the homogeneous case with $D = 0$. Indeed, choosing the right batch size, for example, proportionally to $ K^{\alpha}$ with $\alpha > 0$, one can replace $\sigma^2$ by $\frac{\sigma^2}{K^{\alpha}}$ in \eqref{t_nm_1} and \eqref{t_nm_2} and get convergence guarantees. 

\section{Experiments}
In this section, we present two sets of experiments to validate the performance of Algorithm~\ref{alg4}. 
In Section~\ref{sec:exp1}, we verify the obtained  convergence guarantees on two examples: a strongly-monotone and a monotone bilinear problems, and in Section~\ref{sec:exp2}, we explore the non-monotone case with an application to GANs training. Extended details of the experimental setup can be found in the supplementary material.

\subsection{Verifying Theoretical Convergence Rate}
\label{sec:exp1}

\begin{figure}[h]
\begin{minipage}{\textwidth}
\hfill
\includegraphics[width =  0.3\linewidth ]{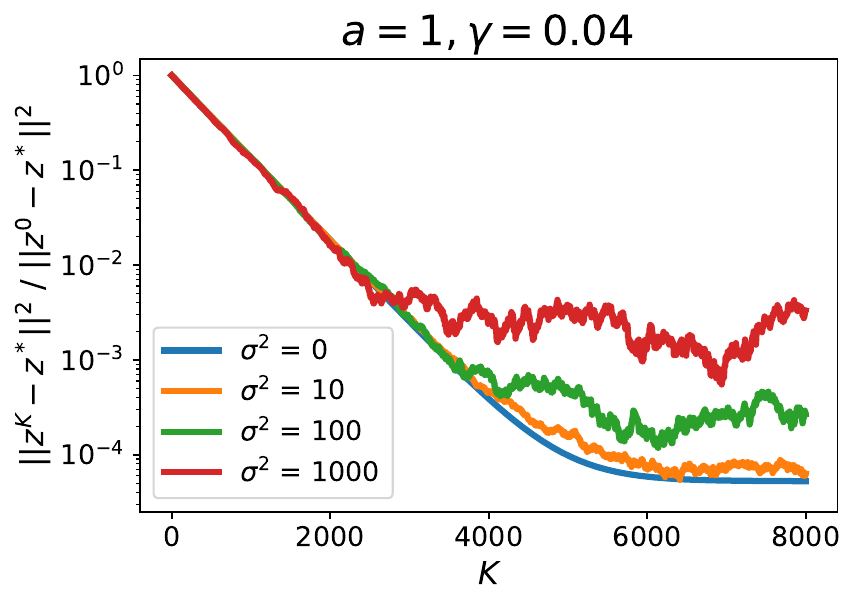}
\hfill
\includegraphics[width =  0.3\linewidth ]{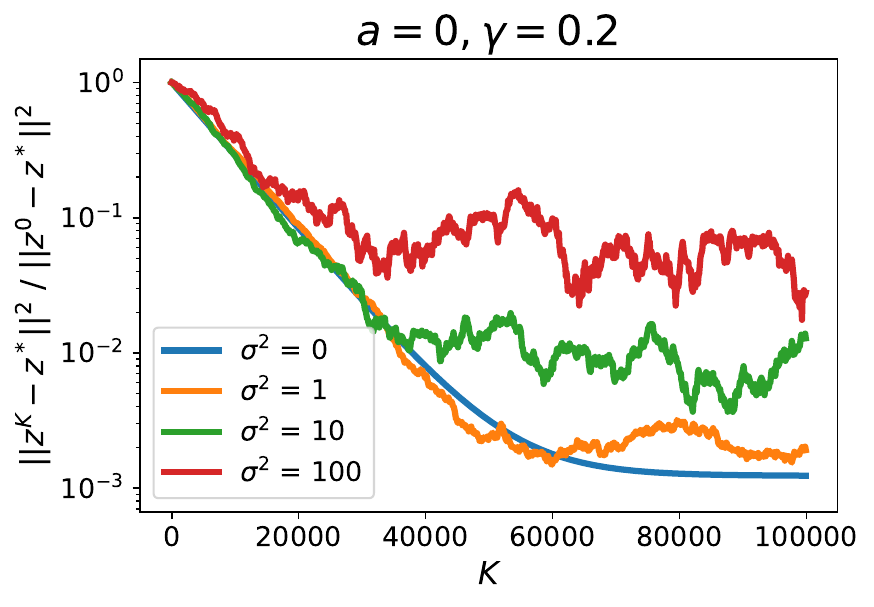}
\hfill\null\\
\vspace{-0.6cm}
\caption{Convergence of Algorithm~\ref{alg4} with constant stepsize in the presence of stochastic noise in strongly-monotone (left) and monotone (right) cases. We observe linear convergence up to an error floor depending on the noise variance and problem parameters  (cf.\ Theorem~\ref{app:th1}).
In Section~\ref{sec_stepsize} of the supplementary material we show convergence to arbitrary accuracy with decreasing stepsizes.}
\label{fig:curves}
\end{minipage}
\end{figure}


First, we focus on the verifying whether the actual behaviour of Algorithm~\ref{alg4}  is predicted by the theoretical convergence rate in Theorem~\ref{app:th1}.

\textbf{Setup.} We consider a distributed bilinear SPP \eqref{eq:dec} with the objective functions 
$
 f_m(x, y) = \frac{a}{2} \|x\|^2 + b x^\top y - \frac{a}{2}\|y\|^2 + c_m^\top x,
$
where $x, y, c_m \in \R^n$, $a,b \in \R$  and $m \in \{1, \dots, M \}$. This set of functions satisfy Assumptions~\ref{a1}, \ref{a2}, \ref{a4} with constants $\mu = a, L^2 = a^2 + b^2$, $D = \max_m \|c_m - \bar c\|$. In this section, we use a ring topology on $M = 20$ nodes with uniform averaging weights, and we set the dimension $n = 5$, $b = 1$, $D \approx 3$, and keep $\tau = 1$. The value of the parameter $p$ in this setting is approximately $0.288$ \citep[Table 1]{koloskova2019decentralized}. To satisfy Assumption~\ref{a3}, we generate stochastic gradients by adding to the exact gradients unbiased Gaussian noise with variance~$\sigma^2$. \looseness=-1

\textbf{Convergence Behaviour.} 
In Figure~\ref{fig:curves}, we show the convergence of Algorithm~\ref{alg4} with a fixed stepsize on the strongly-monotone ($a=1$) and monotone ($a=0$) instances. 
In the strongly-monotone case, we see a linear convergence up to some level defined by the heterogeneity parameter and the noise. The convergence for the non-strongly-monotone problem is slower, but we also see  a linear convergence up to some level (for bilinear problems this behavior is expected from the theoretical point of view \cite{korpelevich1976extragradient}). 
Note that the convergence to some limiting accuracy is expected since when a constant stepsize is used in stochastic optimization/stochastic variational inequalities with strong convexity/monotonicity,  algorithms are usually guaranteed to converge only to a vicinity of the solution, see, e.g., Theorem 2 in \cite{Mishchenko2019:extragradient}. This is also in accordance with Theorem \ref{app:th1} that, for a fixed stepsize, guarantees the convergence to some non-zero limit accuracy and says that, to achieve the zero error, one needs to choose a decreasing stepsize. We additionally validate in Section~\ref{sec_stepsize} of the supplementary material that with a decreasing stepsize, the algorithm can converge to the zero error.

\begin{figure}[h]
\begin{minipage}{\textwidth}
\hfill
\includegraphics[width =  0.3\linewidth ]{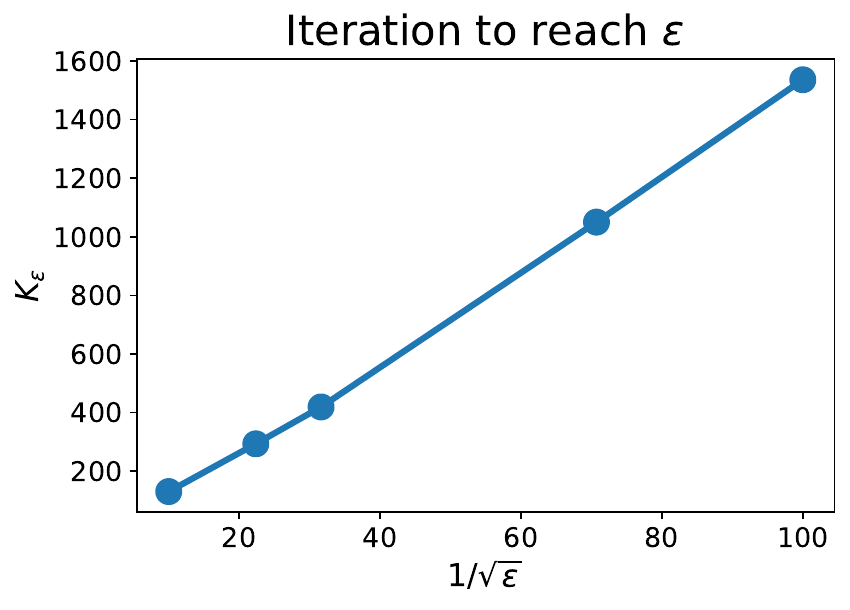}
\hfill
\includegraphics[width =  0.3\linewidth ]{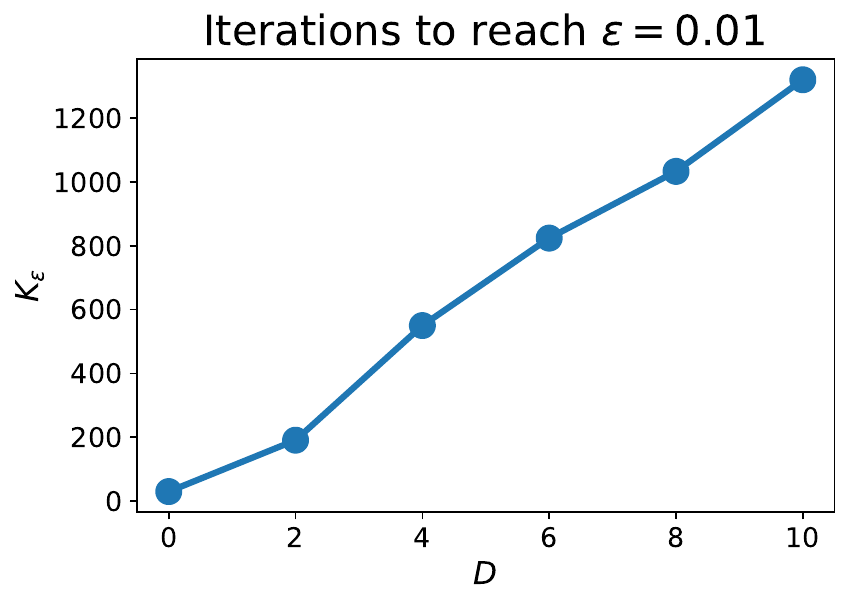}
\hfill
\includegraphics[width =  0.3\linewidth ]{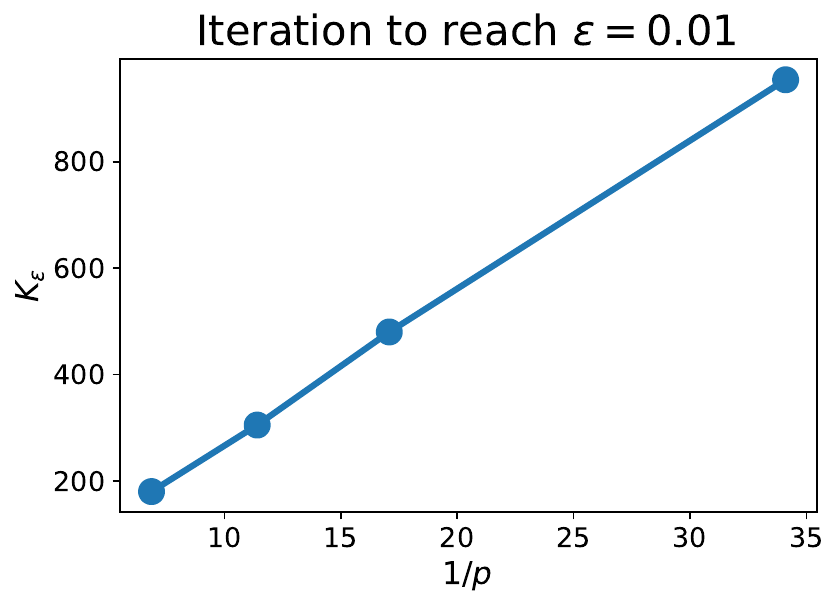}
\hfill\null\\
\vspace{-0.6cm}
\caption{Verifying the $ \mathcal{ O} \big(\frac{ D^2}{p^2  K^2}\big)$ convergence rate for the strongly-monotone noiseless ($ \sigma^2=0$) case.}
\label{fig:T2exp}
\end{minipage}
\end{figure}
\textbf{Dependence on the Heterogeneity parameter $D$.}
In the second set of  experiments, we aim to verify the dependence on the data heterogeneity parameter $D$.
Therefore, we consider the setting when $\sigma^2=0$.
From our theory, equation~\eqref{t_sm}, we predict that the  most significant term in the convergence rate when $\sigma^2=0$ scales as $\mathcal{ O} \bigl(\smash{\frac{ D^2}{p^2  K^2}}\bigr)$ (since the primary goal of this experiment is to study the dependence on $p$, $D$, $K$, we omit all the other fixed parameters for simplicity). We take $b = 1$, $a=1$ and conduct experiments 
with the number of iterations needed to achieve the error $\smash{\textstyle \frac{1}{M}\sum_{m=1}^M \|z_m^{k} - z^* \|^2} < \varepsilon $, for different $\varepsilon$. In all these experiments, the stepsize is tuned  individually.

First, we verify the power of $K$ in the bounds. For this experiment, we keep $D, p$ constant and vary the accuracy $\varepsilon$. As we  can see from the leftmost subplot in Figure~\ref{fig:T2exp}, the number of iterations scales as $K \propto \smash{\frac{1}{\sqrt{\varepsilon}}}$, confirming the predicted $\mathcal{ O}  \bigl(\frac{1}{K^2}\bigr)$ dependency of the error on $K$.
Next, we measure the number of iterations sufficient to reach the error $\varepsilon=0.01$ while varying $D$. The middle plot shows that the number of iterations scales proportionally to $D$ (showing $D \propto K$). Lastly, we depict the number of iterations to reach $\varepsilon=0.01$ while changing the graph parameter $p$ and again observe $\smash{\frac{1}{p}}\propto K$.
Summarizing, these experiments verify the $\mathcal{ O} \bigl(\smash{\frac{ D^2}{p^2  K^2}}\bigr)$ term in the convergence rate.

\subsection{Training GANs}
\label{sec:exp2}

Our algorithm allows combining in the distributed learning setting different communication graph topologies, as well as local steps. Thus, our goal in this section is to illustrate this empirically with the experiments on GANs training. In Section \ref{sec_exp_gans} of the supplementary material, we discuss to what extent our theoretical results hold for GANs training.

\textbf{Data and model.} We consider the CIFAR-10 \cite{cifar10} dataset  containing $60000$ images, equally distributed over $10$ classes. 
We increased the size of the dataset by 4 times using transformations and adding noise.
We simulate a distributed setup of 16 nodes on two GPUs and use Ray \cite{moritz2018ray}. To emulate the heterogeneous setting, we partition the dataset into $16$ subsets. For each subset, we select a major class that forms $20\%$ of the data, while the rest of the data split is filled uniformly by the other classes.
As a basic architecture we choose DCGAN \cite{radford2015unsupervised}, conditioned by class labels, similarly to \cite{cgan} (the network architecture can be found in Section \ref{sec_exp_gans}). We chose Adam \cite{kingma2014adam} as the optimizer. We make one local Adam step and then one gossip averaging step with time-varying matrices $W^k$---similar to Algorithm \ref{alg4}.

\textbf{Setting.} We compare the following three topologies (and the corresponding matrices $W^k$):\\
$\bullet$ \textbf{Full.} Full graph at the end of each epoch, otherwise local steps. This means that we make 120 communication rounds (by communication round we mean the exchange of information between a pair of devices) in an epoch.\\
$\bullet$ \textbf{Local.} Full graph at the end of each 5th epoch, otherwise local steps. This means that we make 24 communication rounds in an epoch (in average: 4 epochs without communications and 1 epoch with 120 rounds).\\
$\bullet$ \textbf{Clusters.} At the end of each epoch, clique clusters of size 4 are randomly formed (in total 4 cliques). This means that we make 24 communication rounds in an epoch.

Note that the communication budget of the first approach is 5 times larger.

We use the same learning rate equal to 0.002 for the generator and discriminator. The rest of the parameters and features of the architecture can be found in the supplementary material.

\begin{wrapfigure}[18]{r}{8.5cm}
\vspace{-0.7cm}
\begin{minipage}{0.6\textwidth}
\includegraphics[width =  0.49\linewidth ]{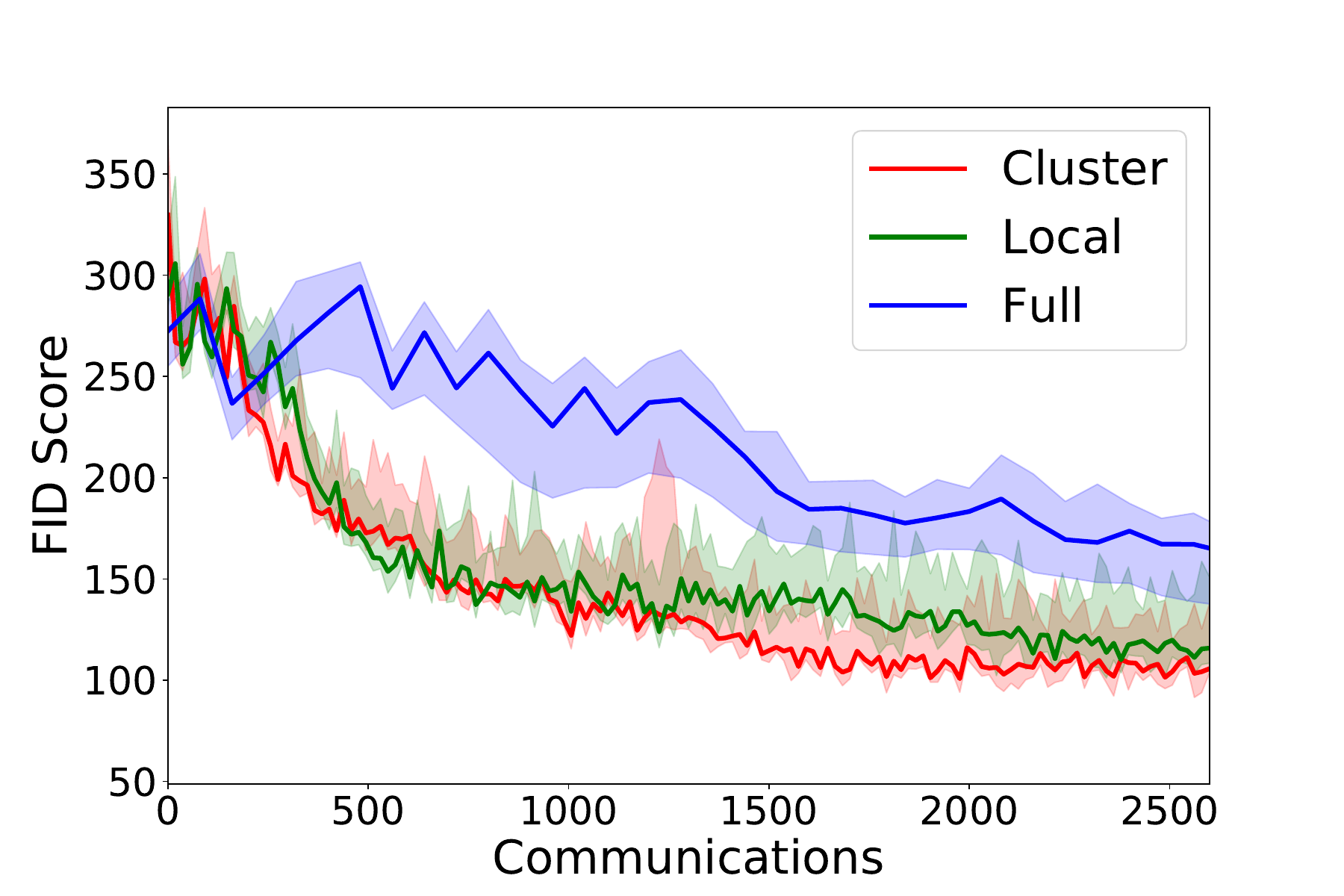}
\includegraphics[width =  0.49\linewidth ]{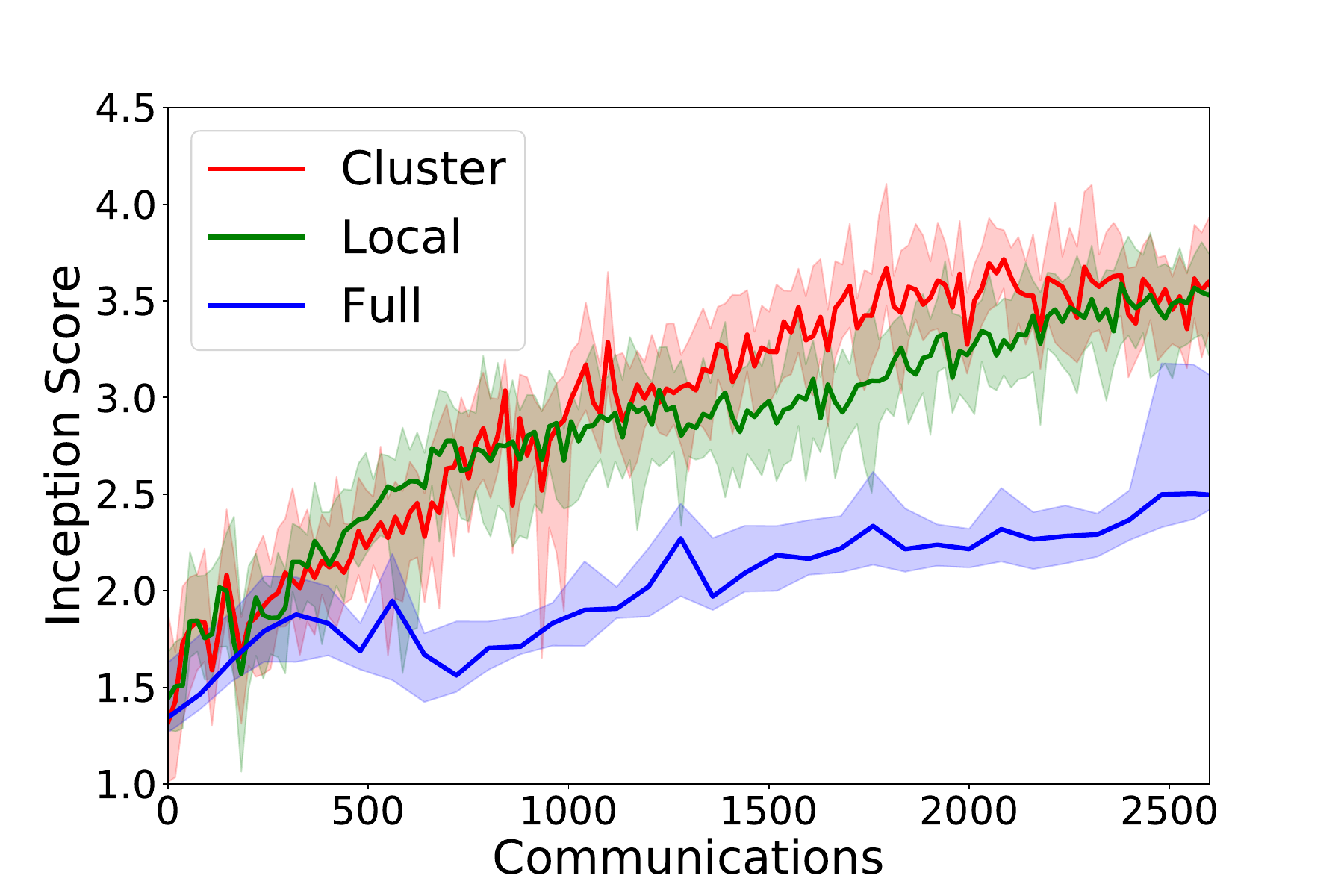}
\hfill\null\\
\includegraphics[width =  0.49\linewidth ]{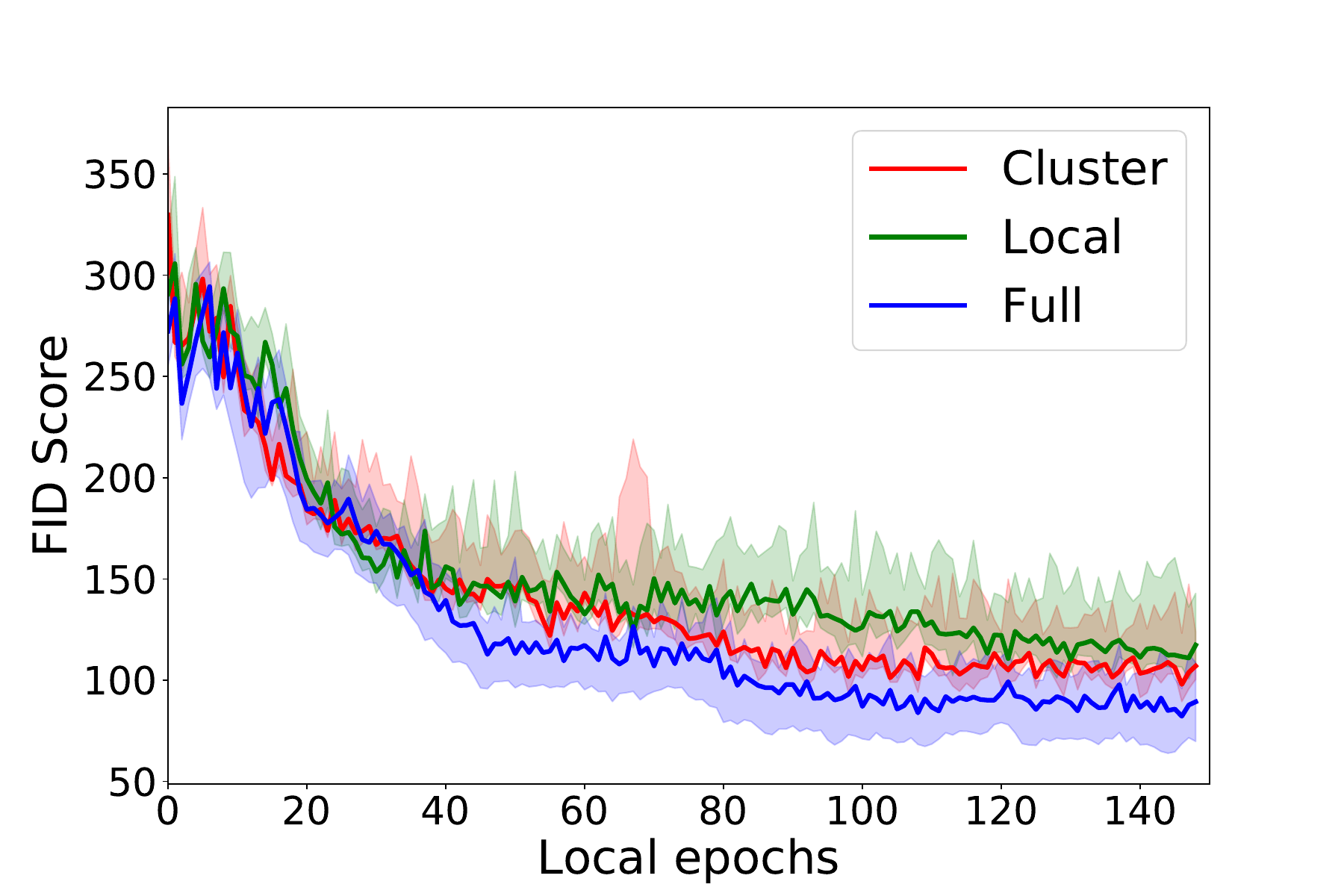}
\includegraphics[width =  0.49\linewidth ]{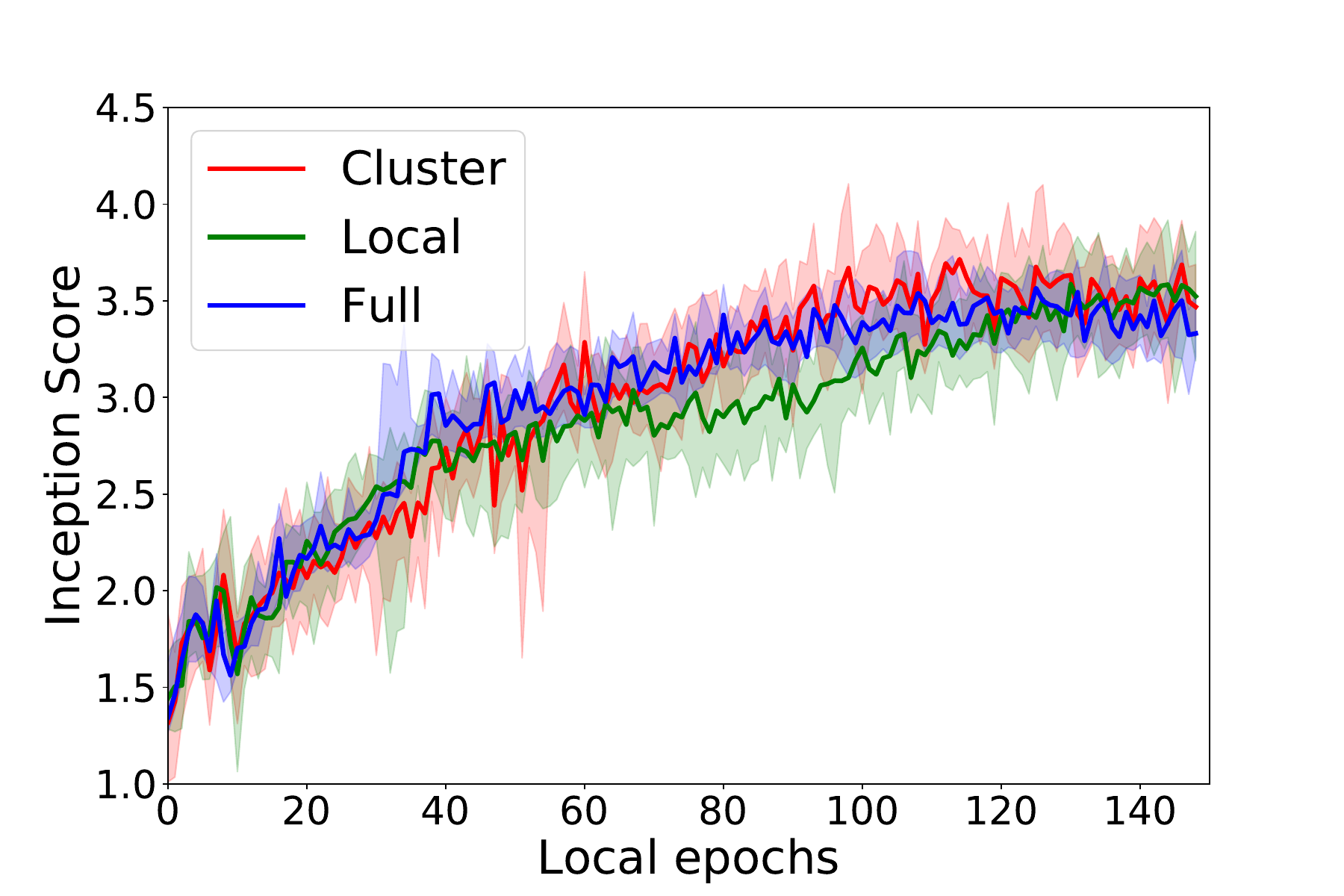}
\hfill\null\\
\vspace{-0.6cm}
\caption{Comparison of the three network topologies  in DCGAN distributed decentralized learning on CIFAR-10. FID Score and Inception Score vs the number of communications (two top), and  the same Scores vs the number of local epochs (two bottom). The experiment was repeated 5 times on different  random data splitting, the maximum and minimum deviations are depicted in the plots by the shade.}
\label{fig:gans1}
\end{minipage}
\end{wrapfigure}

\textbf{Results.} The results of the experiment are presented in Figure \ref{fig:gans1} and Figure \ref{fig:gans2} (Section \ref{sec:pic_gan}). In terms of the number of local epochs, all the methods converged quite close to each other and produced similar images. In terms of communications, Local and Cluster topologies lead to much better results, and the Cluster topology is slightly better than the Local.

\section{Conclusion}
We propose a novel efficient algorithm for solving decentralized stochastic MVIs and SPPs under a very general assumption on the network topology and communication constraints.
In particular, our method is the first decentralized extragradient method with local steps for time-varying network topologies. Moreover, for the proposed algorithm, we prove the convergence rate theorem in the SM, M and NM cases. In the numerical experiments, we verify that  the dependence of our rates on the data heterogeneity parameter $D$ is tight in the SM case, and cannot be further improved in general.
By training DCGAN on a decentralized topology, we demonstrate that our method is effective on practical DL tasks. As a future work it would be interesting to generalize such algorithms for infinite-dimensional problems.

\section*{Acknowledgments}

This research of A. Beznosikov has been supported by The Analytical Center for the Government of the Russian Federation (Agreement No. 70-2021-00143 dd. 01.11.2021, IGK 000000D730321P5Q0002).
The work by P. Dvurechensky in Section C.3 was funded by the Deutsche Forschungsgemeinschaft (DFG, German Research Foundation) under Germany's Excellence Strategy – The Berlin Mathematics Research Center MATH+ (EXC-2046/1, project ID: 390685689).


\bibliography{example_paper_1}
\bibliographystyle{plain}

\newpage
\appendix
\onecolumn

\part*{Supplementary Material}



\section{Experiments} 
In this section, we provide additional details about the experiments reported in the main text and additional experiments with a decreasing stepsize.
We implement all the methods using Python 3.8 using
PyTorch \cite{pytorch} and Ray \cite{ray} 
and run the experiments on a machine with 24 AMD EPYC 7552 @ 2.20GHz processors,
2 GPUs  NVIDIA A100-PCIE with 40536 Mb of memory each (Cuda 11.3).

\subsection{Additional Details of the Experiments with Training GANs} \label{sec_exp_gans}

\begin{figure}[h!]
\begin{minipage}{0.5\textwidth}
  \centering
\includegraphics[width =  0.95\textwidth ]{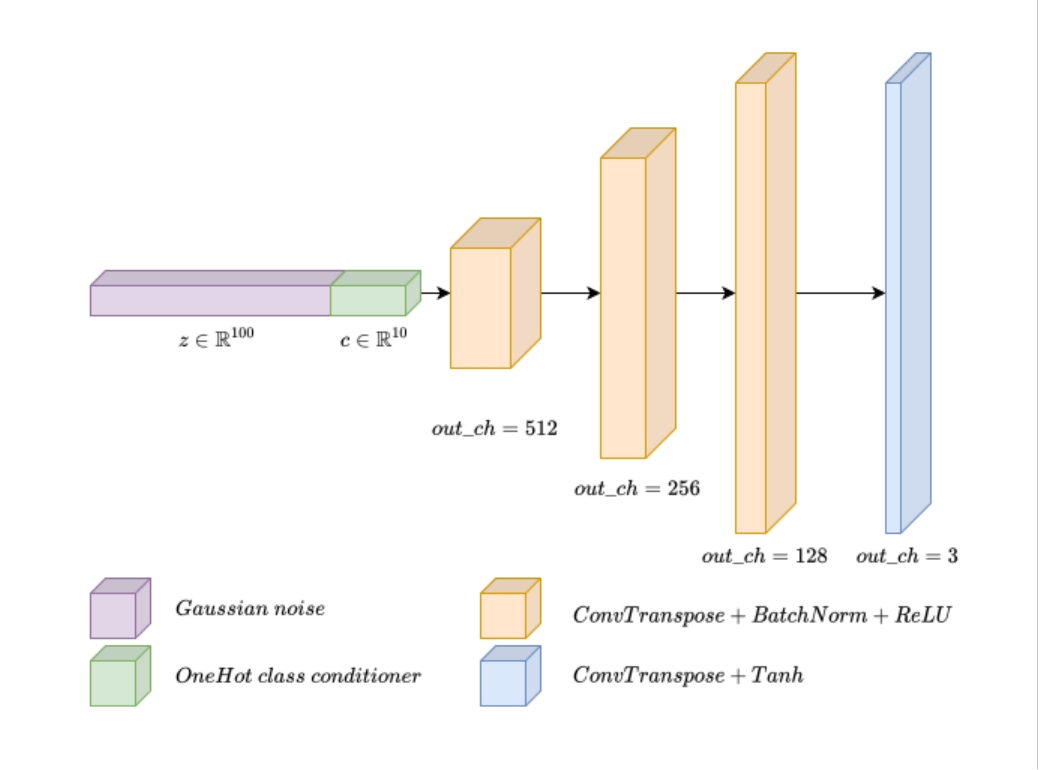}
\end{minipage}%
\begin{minipage}{0.5\textwidth}
  \centering
\includegraphics[width =  0.9\textwidth ]{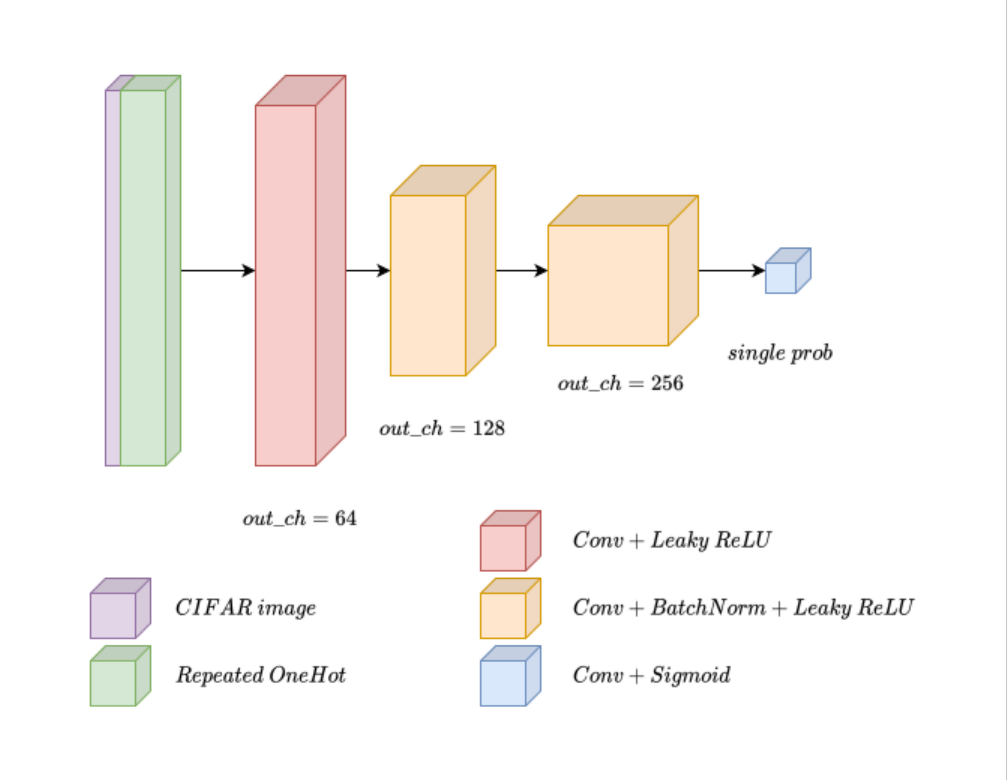}
\end{minipage}%
\\
\begin{minipage}{0.5\textwidth}
  \centering
a) Generator 
\end{minipage}%
\begin{minipage}{0.5\textwidth}
  \centering
b) Discriminator
\end{minipage}%
\caption{DCGAN architecture.}
\label{fig:gans3}
\end{figure}

As mentioned in the main text of the paper, we use DCGAN architecture  \cite{radford2015unsupervised}, conditioned by the class labels, similarly to \cite{cgan}. The illustration of the architecture is provided in Figure \ref{fig:gans3}. In Table 
\ref{tab:hyp}, we give the hyperparameters for all the experiments.

\begin{table}[h]
\centering
\begin{tabular}{ll}
\hline \text {Hyperparameters } & \\
\hline \text { Batch size } & =64 \\
\text { Weight clipping for the discriminator } & =0.01 \\
\text { Learning rate for generator and discriminator} & =0.002 \\
\text{ Initialization: } & \text{normal}\\
\text{ Other parameters: } &\text{default in PyTorch}\\
\hline
\end{tabular}
\vspace{0.2cm}
\caption{Hyperparameters for DCGAN training.}
\label{tab:hyp} 
\end{table}

Next, we comment on the relation between the main assumptions of the theoretical analysis and the example of training GANs.
First of all, the goal of the DCGAN training experiment is to study how the network topology influences the convergence of the algorithm. Even if the assumptions do not hold, we see that the algorithm performs quite well and is flexible w.r.t. the choice of the topology. 
Secondly, as we write in the main text, the made assumptions are quite standard and widely used in the literature.
In particular, Assumptions 3.1 and 3.3 are classical and are often used in the literature, see, e.g. \cite{juditsky2011solving}, including the literature on the neural networks training. Assumption 3.4 is also widely used \cite{khaled2020tighter,woodworth2020minibatch,gorbunov2021local, pmlr-v32-shamir14, arjevani2015communication, hendrikx2020statistically,deng2021local,hou2021efficient,liu2019decentralized}, and holds with a small constant $D$ when the data is uniformly split among the devices. Such splitting can be easily made when one uses a computational cluster with a large amount of data, e.g., images. In Section \ref{sec:exp2}, we deliberately consider a more difficult setup and make the distribution of images over the nodes not uniform, but heterogeneous. As we see, the results of the experiments are quite promising in this case.
Assumption 3.2 (NM) is also used in the literature on the algorithms for training GANs and their analysis \cite{liu2019decentralized,mertikopoulos2018optimistic,liu2019towards}. Moreover, this assumption is shown to hold in some nonconvex minimization problems, for example, when SGD is used for training neural networks \cite{li2017convergence,kleinberg2018alternative,zhou2019sgd}.


\subsection{Additional Experiments with Decreasing Stepsize}\label{sec_stepsize}
As it can be seen from the proofs in the next sections, our theoretical results in Theorem~\ref{app:th1} hold for the fixed stepsizes that optimize the error of the obtained approximate solution given the budget of $K$ iterations. Thus, given a target accuracy $\epsilon > 0$ and using the bounds in Theorem~\ref{app:th1}, we can choose the number of iterations $K=K(\epsilon)$ to guarantee the accuracy $\epsilon$. In turn, based on the value $K(\epsilon)$, we choose the fixed stepsize $\gamma=\gamma(K(\epsilon))$. This procedure of defining the stepsize based on the target accuracy and the corresponding budget of iterations is quite standard in the literature, see, e.g., \cite{juditsky2011solving, stich2019unified, beznosikov2022stochastic}.
In Section \ref{sec:restart}, we provide a generic technique that allows us not to fix the target accuracy in advance and construct a decreasing sequence of stepsizes. This is useful when the desired target accuracy $\epsilon$ is not known, or not determined. 
In this section, we numerically illustrate that  Algorithm~\ref{alg4} can reach arbitrarily small error when implemented with quite simple decreasing stepsizes. We note that, despite not analyzed theoretically, the used in the experiments decreasing stepsize leads to a good performance of the algorithm, which additionally illustrates the flexibility of our approach for practical purposes.

\begin{figure*}[th]
	\centering
	\begin{minipage}{\textwidth}
		\centering
		\begin{minipage}{0.5\textwidth}
    \centering
   \includegraphics[width=0.8\linewidth]{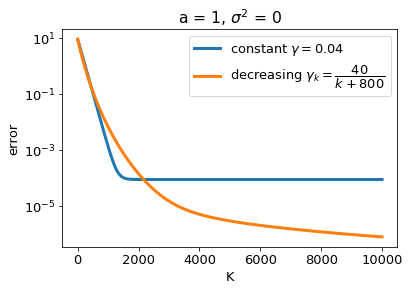}
    \end{minipage}%
	\begin{minipage}{0.5\textwidth}
    \centering
   \includegraphics[width=0.8\linewidth]{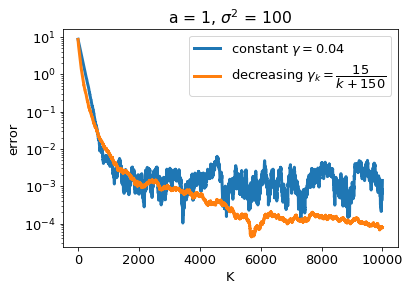}
    \end{minipage}%
		\caption{%
		Convergence of Algorithm~\ref{alg4} with the decreasing stepsizes in the noiseless (left) and stochastic (right) cases. 
		}
		\label{fig:curves_dec}
	\end{minipage}%
\end{figure*}

For the experiments, we consider the same setup as in  Section~\ref{sec:exp1} of the main text (see Figure~\ref{fig:curves}, left), i.e., strongly-monotone bilinear objective functions distributed over the network with the ring topology. We consider two cases: with and without stochastic noise, i.e., we fix either $\sigma = 0$, or  $\sigma = 100$. During the training, we decrease the stepsize as $\gamma_k = \frac{\alpha}{k + \beta}$, where $k$ is the current iteration number. We set $\alpha = 40, \beta = 800$ in the noiseless case and $\alpha = 15, \beta = 150$ when $\sigma^2 = 100$. In Figure~\ref{fig:curves_dec}, we can see that the error decreases to zero with a sublinear rate. This is in contrast to the limiting behavior which we observe in Figure~\ref{fig:curves}, when the algorithm is not able to optimize below a certain threshold. Sublinear convergence may be expleined by the second and third terms in the estimate \eqref{t_sm}.

\subsection{Images Generated by the Trained GAN}\label{sec:pic_gan}
\begin{figure*}[h!]
\begin{minipage}{0.32\textwidth}
  \centering
\includegraphics[width =  0.85\textwidth ]{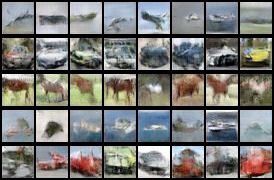}
\end{minipage}\hfill%
\begin{minipage}{0.32\textwidth}
  \centering
\includegraphics[width =  0.85\textwidth ]{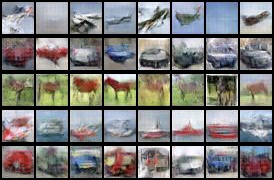}
\end{minipage}\hfill%
\begin{minipage}{0.32\textwidth}
  \centering
\includegraphics[width =  0.85\textwidth ]{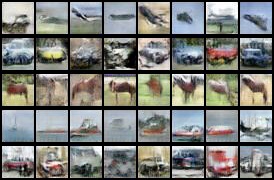}
\end{minipage}\hfill\null%
\\
\begin{minipage}{0.32\textwidth}
  \centering
  \scriptsize{(a) Cluster}
\end{minipage}\hfill%
\begin{minipage}{0.32\textwidth}
  \centering
  \scriptsize{(b) Local}
\end{minipage}\hfill%
\begin{minipage}{0.32\textwidth}
  \centering
  \scriptsize{(c) Full}
\end{minipage}\hfill\null\\
\vspace{-0.6cm}
\caption{Images generated by DCGAN trained distributedly using different communication graph topologies: (a)  Cluster, (b) Local, (c) Full.}
\label{fig:gans2}
\end{figure*}

\section{Useful Facts Used in the Proofs}
Before we start with the proofs, we give several simple facts that are used throughout the proofs of the main theorem.

\textbf{Upper bound for a squared sum.} For an arbitrary integer $n\ge 1$ and arbitrary set of vectors $a_1,\ldots,a_n$, we have
\begin{equation}
    \left\|\sum\limits_{i=1}^n a_i\right\|^2 \le n\sum\limits_{i=1}^n \|a_i\|^2.\label{eq:sqs}
\end{equation}

\textbf{Cauchy-Schwarz inequality.} For arbitrary vectors $a$ and $b$ and any constant $c >0$, we have
\begin{align}
    2 \langle a, b\rangle &\leq c \| a\|^2 + c^{-1}\| b\|^2\label{eq:cs}, \\
    \| a + b\|^2 &\leq (1 + c) \| a\|^2 + (1 + c^{-1}) \| b\|^2 \label{eq:cs1}.
\end{align}

\textbf{Cauchy-Schwarz inequality for random variables.} Let $\xi$ and $\eta$ be real-valued random variables such that $\EE[\xi^2] < \infty$ and $\EE[\eta^2] < \infty$. Then
\begin{equation}
    \EE[\xi\eta] \le \sqrt{\EE[\xi^2]\EE[\eta^2]}\label{eq:cauchy_schwarz_random}.
\end{equation}

\textbf{Frobenius norm of product.} For given matrices $A$ and $B$, it holds that
\begin{align}
    \label{eq:frob}
    \|AB\|_F \leq \| A\|_F \|B \|_2,
\end{align}
where $\|\cdot\|_F$ denotes the Frobenius norm of a matrix, and $\|\cdot \|_2$ is the spectral norm of a matrix, i.e., the maximal singular value.

\section{Missing Proofs for the Main Theorem}
In this section, we provide the proof of the main theorem.  For convenience, we  give its full statement, including the explicit expressions for the stepsizes.

\begin{theorem}[Theorem \ref{app:th1}] 
Let Assumptions \ref{a5}, \ref{a1}, \ref{a3}, \ref{a4} hold and the sequences $\bar z^k$, $\widehat{z}^{k}$ be generated by Algorithm \ref{alg4} that is run for $K > 0$ iterations. Then,

\noindent  $\bullet$ \textbf{Strongly-monotone case:} under Assumption \ref{a2}(SM), with $\gamma = \min\left\{\frac{p}{120\tau L} , \frac{2\ln\left( \max\{2, \mu^2 M r_0 K/(40\sigma^2) \} \right)}{\mu K}\right\}$ it holds that
\begin{align*}
{\E\left[\|\bar z^{K+1} - z^* \|^2\right] = { \mathcal{\tilde O} \left( \|z^{0} - z^* \|^2 \cdot \exp\left( - \frac{\mu K p }{240 L \tau}\right) + \frac{\sigma^2}{\mu^2 M K}  + \frac{L^2 \Delta}{\mu^4   K^2} \right);}}
  \end{align*}
 $\bullet$ \textbf{Monotone case:} under Assumption \ref{a2}(M), for any convex compact $\mathcal{C}$ s.t. $z^0,z^* \in \mathcal{C}$ and $\max_{z,z' \in \mathcal{C}} \|z-z'\|\leq \Omega_{\mathcal{C}}$, with $\gamma =  \min\left\{\frac{1}{3L},\left(\frac{2\Omega_{\mathcal{C}}^2M}{5(K+1)\sigma^2}\right)^{\frac{1}{2}}, \left(\frac{\Omega_{\mathcal{C}}^2}{6(K+1)^2L^2\Delta}\right)^{\frac{1}{4}}\right\}$ it holds that 
\begin{align*}
{ \sup_{z \in \mathcal{C}} \E\left[\left\langle  F(z), \widehat{z}^{K} - z \right\rangle\right] = {
\mathcal{O} \left( \frac{L\Omega_{\mathcal{C}}^2}{K} + \frac{\sigma\Omega_{\mathcal{C}}}{\sqrt{MK}}+\frac{\sqrt{L\Omega_{\mathcal{C}}^3\sqrt{\Delta}}}{\sqrt{K}}+\sqrt{\frac{(\Delta+L^2\Omega_{\mathcal{C}}^2)\Omega_{\mathcal{C}} \sqrt{\Delta}}{KL}}\right)}}.
 \end{align*}
Under the additional assumption that, for all $k$,  $\|\bar z^k\| \leq \Omega$, with $\gamma = \min\left\{\frac{1}{3L},\left(\frac{\Omega_{\mathcal{C}}^2M}{20(K+1)\sigma^2}\right)^{\frac{1}{2}}, \left(\frac{\Omega_{\mathcal{C}}^2}{60(K+1)^2L^2\Delta}\right)^{\frac{1}{4}},\left(\frac{\Omega_{\mathcal{C}}^2}{(K+1)((\Omega +\Omega_{\mathcal{C}})L \sqrt{\Delta} + \Delta)}\right)^{\frac{1}{2}}\right\}$, we have that 
\begin{align*}
{ \sup_{z \in \mathcal{C}} \E\left[\left\langle  F(z), \widehat{z}^{K} - z \right\rangle\right]={
\mathcal{O} \left( \frac{L\Omega_{\mathcal{C}}^2}{K} + \frac{\sigma\Omega_{\mathcal{C}}}{\sqrt{MK}}+\frac{\sqrt{L\Omega_{\mathcal{C}}^3\sqrt{\Delta}}}{K^{3/4}}+\sqrt{\frac{((\Omega +\Omega_{\mathcal{C}})L \sqrt{\Delta} + \Delta)\Omega_{\mathcal{C}}^2 }{K}}\right);}}
 \end{align*}
$\bullet$ \textbf{Non-monotone case:} under Assumption \ref{a2}(NM) and if $\|z^0\| \leq \Omega,\|z^*\| \leq \Omega$, with $\gamma = \min\left\{\frac{1}{5L}, \left(\frac{\|z^0 - z^* \|^2}{(K+1)^2L^2 \Delta}\right)^{1/4}\right\}$:
\begin{align*}
{\E\left[ \frac{1}{K+1}\sum_{k=0}^K\|F(\bar z^k)\|^2\right] = { \mathcal{O} \left(\frac{L^2\Omega^2}{K}+\frac{\sigma^2}{M}+L\Omega\sqrt{\Delta}+\frac{\sqrt{L\Omega\Delta^{3/4}}}{\sqrt{K}}\right).}}
\end{align*} 
Under an additional assumption that, for all $k$, $\|\bar z^k\| \leq \Omega$, with $\gamma =  \min\left\{\frac{1}{5L}, \left(\frac{\Omega^2}{(K+1)L \Delta}\right)^{1/3}\right\}$ ,  we have that 
\begin{align*}
{\E\left[ \frac{1}{K+1}\sum_{k=0}^K\|F(\bar z^k)\|^2\right]= {  \mathcal{O} \Big(\frac{L^2 \Omega^2}{K} + \frac{\sigma^2}{M} + \frac{(L \Omega \Delta)^{2/3}}{K^{1/3}}   + L \Omega \sqrt{\Delta} \Big).}}
\end{align*}
Here $\Delta=\frac{  \tau}{ p}  \bigl( \frac{D^2 \tau}{p} + \sigma^2 \bigr)$.
\end{theorem}

We start with some convenient notations, and then proceed with the proof on the case-by-case basis: strongly-monotone, monotone, non-monotone under the Minty condition.

\subsection{Notation}

We introduce some auxiliary notation as follows.

$\bullet$ Average across all the devices/nodes values of the iterates $z$ and the stochastic realizations of the operators $g$ at iteration $k$:
\begin{align}
\bar z^k &:= \frac{1}{M} \sum\limits_{m=1}^M z^{k}_m, \quad
\bar g^k := \frac{1}{M} \sum\limits_{m=1}^M g^k_m  =\frac{1}{M} \sum\limits_{m=1}^M F_m(z^k_m, \xi^k_m), \label{seq1_1}\\
\bar z^{k+1/3} &:= \frac{1}{M} \sum\limits_{m=1}^M z^{k+1/3}_m, \quad
\bar g^{k+1/3} := \frac{1}{M} \sum\limits_{m=1}^M g^{k+1/3}_m  =\frac{1}{M} \sum\limits_{m=1}^M F_m(z^{k+1/3}_m, \xi^{k+1/3}_m), \label{seq1_2} \\
\bar z^{k+1/3} &= \bar z^k - \gamma \bar g^k, \quad \bar z^{k+2/3} = \bar z^k - \gamma \bar g^{k+1/3}, \quad \bar z^{k+1} = \bar z^{k+2/3}. \label{seq1}
\end{align}
The last equality with $\bar z^{k+1} = \bar z^{k+2/3}$ follows from the fact that one step of the gossip procedure, i.e., step 7 of Algorithm~\ref{alg4} preserves the average over $m$ since the matrix $W$ is doubly stochastic (see Definition \ref{mix}).

$\bullet$ Matrix notation for the collection over all the nodes of the iterates $z$, of the averaged iterates $\bar z$, of the stochastic realizations of the operators $g$ and  of the averaged stochastic realizations of the operators $\bar g$ at iteration $k$:
\begin{align}
Z^k &:= [z^k_1, \ldots, z^k_M], \quad
\bar Z^k := [\bar z^k, \ldots, \bar z^k], \label{eq:matr_not_1}\\
G^k &:= [g^k_1, \ldots, g^k_M], \quad
\bar G^k := [\bar g^k, \ldots, \bar g^k] \label{eq:matr_not_2}, \\
\Phi^k &:= [F_1(z^k_1), \ldots, F_M(z^k_M)],  \quad \bar \Phi^k := \left[\frac{1}{M} \sum\limits_{m=1}^M F_m(z^k_m), \ldots,\frac{1}{M} \sum\limits_{m=1}^M F_m(z^k_m) \right] \label{eq:matr_not_3}.
\end{align}
This notation allows us to rewrite compactly the iterations of Algorithm 1 and the "averaged" dynamics given in \eqref{seq1}:
\begin{align}
\label{seq2}
Z^{k+1/3} &= Z^{k} - \gamma G^k, \quad \bar Z^{k+1/3} = \bar Z^{k} - \gamma \bar G^k, \notag\\
Z^{k+2/3} &= Z^{k} - \gamma G^{k+1/3}, \quad \bar Z^{k+2/3} = \bar Z^{k} - \gamma \bar G^{k+1/3}, \\
Z^{k+1} &= Z^{k+2/3} W^{k}, \quad \bar Z^{k+1} = \bar Z^{k+2/3}.\notag
\end{align}

$\bullet$ Average over devices consensus errors
\begin{align}
\label{seq3}
 \text{Err}(k) = \frac{1}{M}\sum\limits_{m=1}^M \|z_m^{k} - \bar z^{k} \|^2,  \;\;\; \text{Err}(k+1/3) = \frac{1}{M}\sum\limits_{m=1}^M \|z_m^{k+1/3} - \bar z^{k+1/3} \|^2,
 \end{align}
 which can be interpreted as the measure of the discrepancy in the values of the current iterate between the devices at iteration $k$.

\subsection{Proof of Theorem \ref{app:th1}, Strongly-Monotone Case.} 
\label{S:proof_SM}
\subsubsection{General Per-Iterate Estimate}
The goal of this subsection is to derive a general bound for the per-iteration progress of Algorithm \ref{alg4}. This bound will be further refined in the following subsections.
We begin the proof with the following lemma.
\begin{lemma} 
\label{l1}
Let $z,y \in \mathbb{R}^n$. Defining $z^+ = z - y$, for all $u \in \mathbb{R}^n$, we have
\begin{eqnarray*}
 \|z^+ - u\|^2 = \|z - u \|^2 - 2 \langle y, z^+ - u\rangle - \|z^+ - z \|^2.
 \end{eqnarray*}
\end{lemma}
\textbf{Proof:} Simple calculations give the following chain of equalities.
 \begin{eqnarray*}
 \|z^+ - u \|^2 &=& \|z^+ - z + z - u \|^2 \\
&=& \|z-u \|^2 + 2 \langle z^+ - z, z - u \rangle + \|z^+ -z \|^2 \\
&=& \|z-u \|^2 + 2 \langle z^+ - z, z^+ - u \rangle - \|z^+ - z \|^2 \\
&=& \|z-u \|^2 + 2 \langle z^+ - (z - y), z^+ - u \rangle - 2 \langle y, z^+ - u \rangle - \|z^+ - z \|^2 \\
&=& \|z - u\|^2 - 2 \langle y, z^+ - u \rangle - \|z^+ - z \|^2.
 \end{eqnarray*}
\EndProof

Applying the previous lemma with $z^+ = \bar z^{k+2/3}$,  $z=\bar z^{k}$, $u = z^*$ and $y = \gamma \bar g^{k+1/3}$, we have
\begin{align*}
 \|\bar z^{k+2/3} - z^* \|^2 = \|\bar z^k - z^* \|^2 - 2 \gamma \langle \bar g^{k+1/3}, \bar z^{k+2/3} - z^* \rangle - \|\bar z^{k+2/3} - \bar z^k \|^2.
 \end{align*}
From the same lemma, but with $z^+ = \bar z^{k+1/3}$,  $z=\bar z^{k}$, $u = z^{k+2/3}$, $y = \gamma \bar g^{k}$, we have
\begin{align*}
 \|\bar z^{k+1/3} - \bar z^{k+2/3} \|^2 = \|\bar z^k - \bar z^{k+2/3} \|^2 - 2 \gamma \langle \bar g^{k}, \bar z^{k+1/3} - \bar z^{k+2/3} \rangle - \|\bar z^{k+1/3} - \bar z^k \|^2.
 \end{align*}
Combining the two previous equalities, we obtain
 \begin{align*}
 \|\bar z^{k+2/3} - z^* \|^2 + \|\bar z^{k+1/3} - \bar z^{k+2/3} \|^2 =& \|\bar z^k - z^* \|^2 - \|\bar z^{k+1/3} - \bar z^k \|^2 \\
 &- 2 \gamma \langle \bar g^{k+1/3}, \bar z^{k+2/3} - z^* \rangle - 2 \gamma \langle \bar g^{k}, \bar z^{k+1/3} - \bar z^{k+2/3} \rangle.
 \end{align*}
A small rearrangement of the terms gives
\begin{align*}
 \|\bar z^{k+2/3} - z^* \|^2 &+ \|\bar z^{k+1/3} - \bar z^{k+2/3} \|^2 \\
 &= \|\bar z^k - z^* \|^2  - \|\bar z^{k+1/3} - \bar z^k \|^2 \\ 
 &\hspace{0.4cm}- 2 \gamma \langle \bar g^{k+1/3}, \bar z^{k+1/3} - z^* \rangle + 2 \gamma \langle \bar g^{k+1/3} - \bar g^{k}, \bar z^{k+1/3} - \bar z^{k+2/3} \rangle \\
 &\overset{\eqref{eq:cs}}{\leq} \|\bar z^k - z^* \|^2  - \|\bar z^{k+1/3} - \bar z^k \|^2 \\ 
 &\hspace{0.4cm}- 2 \gamma \langle \bar g^{k+1/3}, \bar z^{k+1/3} - z^* \rangle  + \gamma^2 \|\bar g^{k+1/3} - \bar g^{k}\|^2 + \| \bar z^{k+1/3} - \bar z^{k+2/3}\|^2.
 \end{align*}
Taking the full expectation, we get
\begin{align*}
 \E\left[\|\bar z^{k+2/3} - z^* \|^2\right]  &\leq \E\left[\|\bar z^k - z^* \|^2\right] - \E\left[\|\bar z^{k+1/3} - \bar z^k \|^2\right] \notag\\
 &\hspace{0.4cm}- 2 \gamma \E\left[\langle \bar g^{k+1/3}, \bar z^{k+1/3} - z^* \rangle\right] + \gamma^2 \E\left[\|\bar g^{k+1/3} - \bar g^{k}\|^2\right].  
 \end{align*}
Since $\bar z^{k+1} = \bar z^{k+2/3}$, we obtain the following  general bound for the per-iteration progress of Algorithm \ref{alg4}:
 \begin{align}
 \E\left[\|\bar z^{k+1} - z^* \|^2\right]  &\leq \E\left[\|\bar z^k - z^* \|^2\right] - \E\left[\|\bar z^{k+1/3} - \bar z^k \|^2\right] \notag\\
 &\hspace{0.4cm}- 2 \gamma \E\left[\langle \bar g^{k+1/3}, \bar z^{k+1/3} - z^* \rangle\right] + \gamma^2 \E\left[\|\bar g^{k+1/3} - \bar g^{k}\|^2\right]. \label{temp1}
 \end{align}

Our next goal is to consider in more details the last two terms, i.e., $- 2 \gamma \E\left[\langle \bar g^{k+1/3}, \bar z^{k+1/3} - z^* \rangle\right]$ and $\gamma^2 \E\left[\|\bar g^{k+1/3} - \bar g^{k}\|^2\right]$, and estimate them from above. 
 
\subsubsection{Two Auxiliary Estimates for the General Bound} 
In the following two auxiliary lemmas, we prove the estimates mentioned in the end of the previous subsection.

\begin{lemma} \label{l2}
Under Assumptions \ref{a1}, \ref{a2}(SM), \ref{a3}  it holds that
\begin{align}
 \label{temp2}
 - 2 \gamma \E\left[\langle \bar g^{k+1/3}, \bar z^{k+1/3} - z^* \rangle\right] \leq - \gamma \mu \E\left[\| \bar z^{k+1/3} - z^* \|^2\right]  + \frac{\gamma L^2}{\mu} \E\left[\text{ \rm Err}(k+1/3) \right].
  \end{align}
 \end{lemma}

\textbf{Proof:} First, we use the independence of all random vectors $\xi^{i} = (\xi^{i}_1, \ldots , \xi^{i}_m)$ and take the conditional expectation $\E_{\xi^{k+1/3}}$ w.r.t. the vector $\xi^{k+1/3}$, conditioned on the other randomness. This gives:
 \begin{align}
  - 2 \gamma \E \left[\langle \bar g^{k+1/3}, \bar z^{k+1/3} - z^* \rangle\right]
  &\stackrel{\eqref{seq1_2}}{=} - 2 \gamma \E\left[\left\langle  \frac{1}{M} \sum\limits_{m=1}^M \E_{\xi^{k+1/3}}[F_m(z_m^{k+1/3}, \xi_m^{k+1/3})], \bar z^{k+1/3} - z^* \right\rangle \right] \notag\\
  &\overset{\eqref{as3}}{=} -2 \gamma \E\left[\left\langle  \frac{1}{M} \sum\limits_{m=1}^M F_m(z_m^{k+1/3}), \bar z^{k+1/3} - z^* \right\rangle\right] \notag\\
  &= - 2 \gamma \E\left[\left\langle  \frac{1}{M} \sum\limits_{m=1}^M F_m(\bar z^{k+1/3}), \bar z^{k+1/3} - z^* \right\rangle\right] \notag\\
  &\hspace{0.4cm}+ 2 \gamma \E\left[\left\langle  \frac{1}{M} \sum\limits_{m=1}^M [F_m(\bar z^{k+1/3}) - F_m(z_m^{k+1/3})], \bar z^{k+1/3} - z^* \right\rangle\right] \notag\\
  &= - 2 \gamma \E\left[\left\langle  F(\bar z^{k+1/3}), \bar z^{k+1/3} - z^* \right\rangle\right] \notag\\
  &\hspace{0.4cm}+ 2 \gamma \E\left[\left\langle  \frac{1}{M} \sum\limits_{m=1}^M [F_m(\bar z^{k+1/3}) - F_m(z_m^{k+1/3})], \bar z^{k+1/3} - z^* \right\rangle\right]. \label{eq:l2_proof_1}
  \end{align}
Next, we prove that $\langle F(z^*),\bar z^{k+1/3}-z^*\rangle \geq 0$  by contradiction. To that end, assume that $\langle F(z^*),\bar z^{k+1/3}-z^*\rangle < 0$. By Assumption \ref{a1}, $F$ is $L$-Lipschitz continuous, and, hence, for a small enough $\alpha > 0$ it holds that $\langle F(\hat z),\bar z^{k+1/3}-z^*\rangle < 0$, where $\hat z = z^* + \alpha (\bar z^{k+1/3} - z^*)$. We substitute $\alpha \bar z^{k+1/3} = \hat z - (1 - \alpha) z^*$ to the inequality $\langle F(\hat z),\alpha \bar z^{k+1/3}-\alpha z^*\rangle < 0$ and get $\langle F(\hat z),\hat z- z^*\rangle < 0$. But, this contradicts the definition of the solution $z^*$ in \eqref{eq:mvi} which implies that $\langle F(\hat z), \hat z - z^*\rangle \geq 0$. Thus, we have that
$\langle F(z^*),\bar z^{k+1/3}-z^*\rangle \geq 0$. 
We combine this inequality with \eqref{eq:l2_proof_1} and obtain
 \begin{align*}
  - 2 \gamma \E \left[\langle \bar g^{k+1/3}, \bar z^{k+1/3} - z^* \rangle\right] 
  &\stackrel{\eqref{eq:l2_proof_1}}{\leq} - 2 \gamma \E\left[\left\langle  F(\bar z^{k+1/3}) - F(z^*), \bar z^{k+1/3} - z^* \right\rangle\right] \\
  &\hspace{0.4cm}+ 2 \gamma \E\left[\left\langle  \frac{1}{M} \sum\limits_{m=1}^M [F_m(\bar z^{k+1/3}) - F_m(z_m^{k+1/3})], \bar z^{k+1/3} - z^* \right\rangle\right] \\
  &\overset{\eqref{as2}}{\leq} - 2 \gamma \mu \E\left[\| \bar z^{k+1/3} - z^* \|^2\right] \\
  &\hspace{0.4cm}+ 2 \gamma \E\left[\left\langle  \frac{1}{M} \sum\limits_{m=1}^M [F_m(\bar z^{k+1/3}) - F_m(z_m^{k+1/3})], \bar z^{k+1/3} - z^* \right\rangle\right].
  \end{align*}
Applying \eqref{eq:cs} with $c = \mu >0$, we further get
\begin{align*}
  - 2 \gamma \E &\left[\langle \bar g^{k+1/3}, \bar z^{k+1/3} - z^* \rangle\right] \leq - 2 \gamma \mu \E\left[\| \bar z^{k+1/3} - z^* \|^2\right] \\
  &\hspace{0.4cm}+ \gamma \mu \E\left[\left\| \bar z^{k+1/3} - z^* \right\|^2\right] + \frac{\gamma}{\mu} \E\left[\left\| \frac{1}{M} \sum\limits_{m=1}^M [F_m(\bar z^{k+1/3}) - F_m(z_m^{k+1/3})] \right\|^2\right] \\
  &= -\gamma \mu \E\left[\| \bar z^{k+1/3} - z^* \|^2\right]   + \frac{\gamma}{\mu M^2} \E\left[\left\| \sum\limits_{m=1}^M [F_m(\bar z^{k+1/3}) - F_m(z_m^{k+1/3})] \right\|^2\right] \\
  &\overset{\eqref{eq:sqs}}{\leq} -\gamma \mu \E\left[\| \bar z^{k+1/3} - z^* \|^2\right]   + \frac{\gamma}{\mu M} \E\left[\sum\limits_{m=1}^M\left\| F_m(\bar z^{k+1/3}) - F_m(z_m^{k+1/3})\right\|^2\right] \\
  &\overset{\eqref{as1}}{\leq} -\gamma \mu \E\left[\| \bar z^{k+1/3} - z^* \|^2\right]  + \frac{\gamma L^2}{\mu M} \E\left[ \sum\limits_{m=1}^M\left\|\bar z^{k+1/3} - z_m^{k+1/3}\right\|^2 \right].
  \end{align*}
Applying \eqref{seq3} to the last term, we finish the proof.
\EndProof

\begin{lemma} Under Assumptions \ref{a1}, \ref{a3}  it holds that 
\label{l3}
\begin{align}
 \label{temp3}
 \E\left[\|\bar g^{k+1/3} - \bar g^{k}\|^2\right] 
 &\leq 5L^2 \E\left[ \|\bar z^{k+1/3} - \bar z^{k}\|^2\right] + \frac{10 \sigma^2}{M} \notag\\
 &\hspace{0.4cm}+ 5L^2\E\left[\text{\rm Err}(k+1/3) \right] +5L^2\E\left[\text{\rm Err}(k) \right].
  \end{align}
\end{lemma}
 
\textbf{Proof:} Consider the following chain of inequalities:
 \begin{align*}
 \E&\left[\|\bar g^{k+1/3} - \bar g^{k}\|^2\right] \stackrel{\eqref{seq1_1},\eqref{seq1_2}}{=} \E\left[\left\|\frac{1}{M}\sum\limits_{m=1}^M F_m(z^{k+1/3}_m, \xi^{k+1/3}_m) - \frac{1}{M}\sum\limits_{m=1}^M F_m(z^k_m, \xi^{k}_m)\right\|^2\right] \\
 &\leq 5 \E\left[\left\|\frac{1}{M}\sum\limits_{m=1}^M [F_m(z^{k+1/3}_m, \xi^{k+1/3}_m) - F_m(z^{k+1/3}_m)  ]\right\|^2\right]  \\
 &\hspace{0.4cm} +5\E\left[\left\|\frac{1}{M}\sum\limits_{m=1}^M [ F_m(z^{k+1/3}_m) - F_m(\bar z^{k+1/3})  ]\right\|^2\right] + 5\E\left[\left\|\frac{1}{M}\sum\limits_{m=1}^M [ F_m(\bar z^{k+1/3}) - F_m(\bar z^{k})  ]\right\|^2\right]  \\
 &\hspace{0.4cm} + 5\E\left[\left\|\frac{1}{M}\sum\limits_{m=1}^M [ F_m(z^{k}_m) - F_m(\bar z^{k})  ]\right\|^2\right] + 5  \E\left[\left\|\frac{1}{M}\sum\limits_{m=1}^M [F_m(z^{k}_m, \xi^k_m) - F_m(z^{k}_m)  ]\right\|^2\right] \\
 &\overset{\eqref{eq:sqs}}{\leq} 5 \E\left[\left\|\frac{1}{M}\sum\limits_{m=1}^M [F_m(z^{k+1/3}_m, \xi^{k+1/3}_m) - F_m(z^{k+1/3}_m)  ]\right\|^2\right] \\
 &\hspace{0.4cm}+ 5  \E\left[\left\|\frac{1}{M}\sum\limits_{m=1}^M [F_m(z^{k}_m, \xi^k_m) - F_m(z^{k}_m)  ]\right\|^2\right]  \\
 &\hspace{0.4cm}+\frac{5}{M}\sum\limits_{m=1}^M \E\left[\left\| F_m(z^{k+1/3}_m) - F_m(\bar z^{k+1/3})\right\|^2\right]  + \frac{5}{M} \sum\limits_{m=1}^M\E\left[\left\| F_m(z^{k}_m) - F_m(\bar z^{k})  \right\|^2\right] \\
 &\hspace{0.4cm}+ 5\E\left[\left\|F(\bar z^{k+1/3}) - F(\bar z^{k})\right\|^2\right] \\
 &\overset{\eqref{as1},\eqref{seq3}}{\leq} 5\E\left[\left\|\frac{1}{M}\sum\limits_{m=1}^M [F_m(z^{k+1/3}_m, \xi^{k+1/3}_m) - F_m(z^{k+1/3}_m)  ]\right\|^2\right]  \\
 &\hspace{0.4cm}+ 5  \E\left[\left\|\frac{1}{M}\sum\limits_{m=1}^M [F_m(z^{k}_m, \xi^k_m) - F_m(z^{k}_m)  ]\right\|^2\right] \\
 &\hspace{0.4cm} +5L^2\E\left[\text{Err}(k+1/3) \right] +5L^2\E\left[\text{Err}(k) \right] + 5 L^2\E\left[ \|\bar z^{k+1/3} - \bar z^{k}\|^2 \right] \\
 &= 5\E\left[\E_{\xi^{k+1/3}}\left[\left\|\frac{1}{M}\sum\limits_{m=1}^M [F_m(z^{k+1/3}_m, \xi^{k+1/3}_m) - F_m(z^{k+1/3}_m)  ]\right\|^2\right]\right]  \\
 &\hspace{0.4cm}+ 5  \E\left[\E_{\xi^{k}}\left[\left\|\frac{1}{M}\sum\limits_{m=1}^M [F_m(z^{k}_m, \xi^k_m) - F_m(z^{k}_m)  ]\right\|^2\right]\right] \\
 &\hspace{0.4cm} +5L^2\E\left[\text{Err}(k+1/3) \right] +5L^2\E\left[\text{Err}(k) \right] + 5 L^2\E\left[ \|\bar z^{k+1/3} - \bar z^{k}\|^2 \right].
 \end{align*}
 Using the independence of the realizations of $\xi^{k},\xi^{k+1/3}$ in each node and \eqref{as3}, we get:
  \begin{align*}
 \E\left[\|\bar g^{k+1/3} - \bar g^{k}\|^2\right]  \leq \frac{10 \sigma^2}{M} +5L^2\E\left[\text{Err}(k+1/3) \right] +5L^2\E\left[\text{Err}(k) \right] + 5 L^2\E\left[ \|\bar z^{k+1/3} - \bar z^{k}\|^2 \right].
  \end{align*}
\EndProof

\subsubsection{Refined General Bound}
We now return to the general bound for the per-iteration progress and combine the general bound \eqref{temp1} with the two estimates \eqref{temp2} and \eqref{temp3} obtained in the previous subsection. In this way, we obtain
\begin{align*}
 \E&\left[\|\bar z^{k+1} - z^* \|^2\right]  \leq \E\left[\|\bar z^k - z^* \|^2\right] - \E\left[\|\bar z^{k+1/3} - \bar z^k \|^2\right] \notag\\
 &\hspace{0.4cm}- \gamma \mu \E\left[\| \bar z^{k+1/3} - z^* \|^2\right]  + \frac{\gamma L^2}{\mu} \E\left[\text{Err}(k+1/3) \right]\\
 &\hspace{0.4cm}+ \gamma^2 \left( 5L^2 \E\left[ \|\bar z^{k+1/3} - \bar z^{k}\|^2\right] + \frac{10 \sigma^2}{M}  + 5L^2\E\left[\text{Err}(k+1/3) \right] +5L^2\E\left[\text{Err}(k) \right] \right).
  \end{align*}
Using \eqref{eq:cs1} with $c = 1$, $a = \bar z^{k+1/3} - z^*$ and $b = \bar z^{k+1/3} - \bar z^k$, we get 
\[
-\| \bar z^{k+1/3} - z^* \|^2 \leq -\tfrac{1}{2} \| \bar z^{k} - z^* \|^2 + \| \bar z^{k+1/3} - z^k \|^2,
\]
which in combination with the previous inequality gives
\begin{align}
\label{temp16}
 \E&\left[\|\bar z^{k+1} - z^* \|^2\right]  \leq \left(1 - \frac{\gamma \mu}{2} \right)\E\left[\|\bar z^k - z^* \|^2\right] - (1 - 5\gamma^2 L^2 - \gamma \mu)\E\left[\|\bar z^{k+1/3} - \bar z^k \|^2\right] \notag\\
 &\hspace{0.4cm}+\left(\frac{\gamma L^2}{\mu} + 5\gamma^2 L^2\right) \E\left[\text{Err}(k+1/3) \right] +  5\gamma^2 L^2\E\left[\text{Err}(k) \right] + \frac{10 \gamma^2 \sigma^2}{M}.  
  \end{align}
Since, by the Theorem assumptions, we have $\gamma \leq \frac{1}{3L}$, the refined general bound for the per-iteration progress becomes
\begin{align}
\label{temp41}
 \E&\left[\|\bar z^{k+1} - z^* \|^2\right]  \leq \left(1 - \frac{\gamma \mu}{2} \right)\E\left[\|\bar z^k - z^* \|^2\right] \notag\\
 &\hspace{0.4cm}+\left(\frac{\gamma L^2}{\mu} + 5\gamma^2 L^2\right) \E\left[\text{Err}(k+1/3) \right] +  5\gamma^2 L^2\E\left[\text{Err}(k) \right] + \frac{10 \gamma^2 \sigma^2}{M}.  
  \end{align}

Our next goal is to bound the consensus error terms $\E\left[\text{Err}(k) \right]$ and $\E\left[\text{Err}(k+1/3) \right]$. 

\subsubsection{Bounds for the Consensus Errors }

The bounds for the consensus error terms are proved in the following two technical lemmas that give recursions for these error terms.

\begin{lemma} \label{l4} Under Assumptions \ref{a1}, \ref{a3}, \ref{a4}, \ref{a5}, for $h = \lfloor k / \tau \rfloor - 1$,  it holds that 
 \begin{align}
 \label{temp4}
\E\left[\text{\rm Err}(k)\right] &\leq 
\left( 1 - \frac{3p}{4}\right) \E [\text{\rm Err}({h\tau}) ] + \frac{144 \gamma^2 L^2 \tau}{p} \sum\limits_{j=h\tau}^{k-1} \E\left[\text{\rm Err}(j+1/3)\right] \notag\\
&\hspace{0.4cm}+ \left( \frac{72 D^2 \tau}{p} +  8 \sigma^2\right) \sum\limits_{j=h\tau}^{k-1} \gamma^2, \\
 \label{temp5}
\E\left[\text{\rm Err}(k+1/3)\right] &\leq 
\left( 1 - \frac{3p}{4}\right) \E [\text{\rm Err}({h\tau}) ] + \frac{216 \gamma^2 L^2 \tau}{p} \sum\limits_{j=h\tau}^{k-1} \E\left[\text{\rm Err}(j+1/3)\right] + \frac{216 \gamma^2 L^2 \tau}{p} \E\left[\text{\rm Err}(k)\right] 
 \notag\\
&\hspace{0.4cm} + \left( \frac{108 D^2 \tau}{p} +  12 \sigma^2\right) \sum\limits_{j=h\tau}^{k-1} \gamma^2 
+ \left( \frac{108 D^2 \tau}{p} +  12 \sigma^2\right) \gamma^2.
\end{align}
  \end{lemma}
  
\textbf{Proof:}
Using the matrix notation introduced in \eqref{eq:matr_not_1}, \eqref{eq:matr_not_2}, \eqref{eq:matr_not_3}, \eqref{seq2}, we rewrite the error $\text{Err}(k)$ as follows:
\begin{align*}
M\cdot \E&\left[\text{Err}(k)\right] =\EE \|Z^{k} - \bar Z^{k} \|^2_F = \EE \| Z^k - \bar Z^{h\tau} - \bar Z^k +  \bar Z^{h\tau}\|^2_F \\ 
&= \E \Bigg[\Bigg\|Z^{h\tau} \prod\limits_{i = h\tau}^{k-1} W^i - \bar Z^{h\tau} - \gamma \sum\limits_{j=h\tau}^{k-1} G^{j+1/3} \prod\limits_{i = j}^{k-1} W^i \\
&\hspace{0.4cm}- \left(  \bar Z^{h\tau} \prod\limits_{i = h\tau}^{k-1} W^i - \bar Z^{h\tau} - \gamma \sum\limits_{j=h\tau}^{k-1} \bar G^{j+1/3} \prod\limits_{i = j}^{k-1} W^i \right)\Bigg\|^2_F\Bigg] \\
&= \E \Bigg[\E_{\xi^{(k-1)+1/3}} \Bigg[\Bigg\|Z^{h\tau} \prod\limits_{i = h\tau}^{k-1} W^i - \bar Z^{h\tau } - \left(  \bar Z^{h\tau} \prod\limits_{i = h\tau}^{k-1} W^i - \bar Z^{h\tau} \right) \\
&\hspace{0.4cm}- \gamma \sum\limits_{j=h\tau}^{k-1} (\Phi^{j+1/3} - \bar \Phi^{j+1/3} ) \prod\limits_{i = j}^{k-1} W^i \\
&\hspace{0.4cm}- \gamma \sum\limits_{j=h\tau}^{k-1} (G^{j+1/3} - \Phi^{j+1/3} - \bar G^{j+1/3} +\bar \Phi^{j+1/3} ) \prod\limits_{i = j}^{k-1} W^i \Bigg\|^2_F \Bigg]  \Bigg].
\end{align*}
Since only $G^{(k - 1) + 1/3}$ and $\Phi^{(k-1)+1/3}$ depend on $\xi^{(k -1) + 1/3}$, and $\E_{\xi^{(k-1)+1/3}} G^{(k -1) + 1/3} = \Phi^{(k-1)+1/3}$, $\E_{\xi^{(k-1)+1/3}} \bar G^{(k -1) + 1/3} = \bar \Phi^{(k-1)+1/3}$ (stochastic oracle is unbiased, see \eqref{as3}), we have
\begin{align*}
M\cdot \E\left[\text{Err}(k)\right] 
&= \E \Bigg[\Bigg\|Z^{h\tau} \prod\limits_{i = h\tau}^{k-1} W^i - \bar Z^{h\tau } - \left(  \bar Z^{h\tau} \prod\limits_{i = h\tau}^{k-1} W^i - \bar Z^{h\tau} \right) \\
&\hspace{0.4cm}- \gamma {\color{blue}\sum\limits_{j=h\tau}^{k-1} (\Phi^{j+1/3} - \bar \Phi^{j+1/3} ) \prod\limits_{i = j}^{k-1} W^i} \\
&\hspace{0.4cm}- \gamma \sum\limits_{j=h\tau}^{k-2} (G^{j+1/3} - \Phi^{j+1/3} - \bar G^{j+1/3} +\bar \Phi^{j+1/3} ) \prod\limits_{i = j}^{k-1} W^i \Bigg\|^2_F \Bigg] \\
&\hspace{0.4cm}+\gamma^2\E \Bigg[\left\|(G^{(k -1) + 1/3} - \Phi^{(k -1) + 1/3} - \bar G^{(k -1) + 1/3} +\bar \Phi^{(k -1) + 1/3} ) W^{k-1} \right\|^2_F \Bigg] .
\end{align*}
Next, we  consider the blue term in the previous display and apply \eqref{eq:cs1} with $c = \beta_1$, where the constant $\beta_1$ is defined below, $a=\Phi^{(k -1) + 1/3}- \bar \Phi^{(k -1) + 1/3}$ and $b$ collecting all the other terms in the first squared norm. This, combined with \eqref{eq:frob} and the fact that $\| W^{k-1}\|_2 \leq 1$, gives us
\begin{align*}
M\cdot \E\left[\text{Err}(k)\right] &\leq {\color{blue}(1 + \beta_1)}\E \Bigg[\Bigg\|Z^{h\tau} \prod\limits_{i = h\tau}^{k-1} W^i - \bar Z^{h\tau } - \left(  \bar Z^{h\tau} \prod\limits_{i = h\tau}^{k-1} W^i - \bar Z^{h\tau} \right) \\
&\hspace{0.4cm}- \gamma {\color{blue}\sum\limits_{j=h\tau}^{k-2} (\Phi^{j+1/3} - \bar \Phi^{j+1/3} ) \prod\limits_{i = j}^{k-1} W^i} \\
&\hspace{0.4cm}- \gamma \sum\limits_{j=h\tau}^{k-2} (G^{j+1/3} - \Phi^{j+1/3} - \bar G^{j+1/3} +\bar \Phi^{j+1/3} ) \prod\limits_{i = j}^{k-1} W^i \Bigg\|^2_F \Bigg] \\
&\hspace{0.4cm}+{\color{blue}(1 + \beta_1^{-1})\gamma^2\E \Bigg[\left\|\Phi^{(k -1) + 1/3} -\bar \Phi^{(k -1) + 1/3}\right\|^2_F \Bigg] }\\
&\hspace{0.4cm}+\gamma^2\E \Bigg[\left\|G^{(k -1) + 1/3} - \Phi^{(k -1) + 1/3} - \bar G^{(k -1) + 1/3} +\bar \Phi^{(k -1) + 1/3} \right\|^2_F \Bigg] .
\end{align*}
In the same way, again using the unbiasedness, we separate the terms with index $(k-2) + 1/3$ using also \eqref{eq:cs1} with $c = \beta_2$, where the constant $\beta_2$ is defined below 
\begin{align*}
M\cdot \E\left[\text{Err}(k)\right] 
 &\leq (1 + \beta_1)\E \Bigg[\Bigg\|Z^{h\tau} \prod\limits_{i = h\tau}^{k-1} W^i - \bar Z^{h\tau } - \left(  \bar Z^{h\tau} \prod\limits_{i = h\tau}^{k-1} W^i - \bar Z^{h\tau} \right) \\
&\hspace{0.4cm}- \gamma \sum\limits_{j=h\tau}^{k-2} (\Phi^{j+1/3} - \bar \Phi^{j+1/3} ) \prod\limits_{i = j}^{k-1} W^i \\
&\hspace{0.4cm}- \gamma \sum\limits_{j=h\tau}^{k-3} (G^{j+1/3} - \Phi^{j+1/3} - \bar G^{j+1/3} +\bar \Phi^{j+1/3} ) \prod\limits_{i = j}^{k-1} W^i \Bigg\|^2_F \Bigg] \\
&\hspace{0.4cm}+(1 + \beta_1^{-1})\gamma^2\E \Bigg[\left\|\Phi^{(k -1) + 1/3} -\bar \Phi^{(k -1) + 1/3}\right\|^2_F \Bigg]\\
&\hspace{0.4cm}+(1 + \beta_1)\gamma^2\E \Bigg[\left\|G^{(k-2)+1/3} - \Phi^{(k-2)+1/3} - \bar G^{(k-2)+1/3} +\bar \Phi^{(k-2)+1/3}  \right\|^2_F \Bigg]  \\
&\hspace{0.4cm}+\gamma^2\E \Bigg[\left\|G^{(k -1) + 1/3} - \Phi^{(k -1) + 1/3} - \bar G^{(k -1) + 1/3} +\bar \Phi^{(k -1) + 1/3} \right\|^2_F \Bigg] \\
&\leq (1 + \beta_1)(1 + \beta_2)\E \Bigg[\Bigg\|Z^{h\tau} \prod\limits_{i = h\tau}^{k-1} W^i - \bar Z^{h\tau } - \left(  \bar Z^{h\tau} \prod\limits_{i = h\tau}^{k-1} W^i - \bar Z^{h\tau} \right) \\
&\hspace{0.4cm}- \gamma \sum\limits_{j=h\tau}^{k-3} (\Phi^{j+1/3} - \bar \Phi^{j+1/3} ) \prod\limits_{i = j}^{k-1} W^i \\
&\hspace{0.4cm}- \gamma \sum\limits_{j=h\tau}^{k-3} (G^{j+1/3} - \Phi^{j+1/3} - \bar G^{j+1/3} +\bar \Phi^{j+1/3} ) \prod\limits_{i = j}^{k-1} W^i \Bigg\|^2_F \Bigg] \\
&\hspace{0.4cm}+(1 + \beta_1^{-1})\gamma^2\E \Bigg[\left\|\Phi^{(k-1)+1/3} -\bar \Phi^{(k-1)+1/3}  \right\|^2_F \Bigg] \\
&\hspace{0.4cm}+(1 + \beta_1)(1 + \beta_2^{-1})\gamma^2\E \Bigg[\left\|\Phi^{(k-2)+1/3} -\bar \Phi^{(k-2)+1/3}  \right\|^2_F \Bigg] \\
&\hspace{0.4cm}+(1 + \beta_1)\gamma^2\E \Bigg[\left\|G^{(k-2)+1/3} - \Phi^{(k-2)+1/3} - \bar G^{(k-2)+1/3} +\bar \Phi^{(k-2)+1/3}  \right\|^2_F \Bigg]  \\
&\hspace{0.4cm}+\gamma^2\E \Bigg[\left\|G^{(k-1)+1/3} - \Phi^{(k-1)+1/3} - \bar G^{(k-1)+1/3} +\bar \Phi^{(k-1)+1/3}  \right\|^2_F \Bigg] .
\end{align*}

Proceeding in a similar way for all the terms, we obtain 
\begin{align}
M\cdot \E\left[\text{Err}(k)\right] &\leq (1 + \beta_1)\ldots (1 + \beta_{k-1 - h\tau}) \E \Bigg[\Bigg\|Z^{h\tau} \prod\limits_{i = h\tau}^{k-1} W^i - \bar Z^{h\tau } - \left(  \bar Z^{h\tau} \prod\limits_{i = h\tau}^{k-1} W^i - \bar Z^{h\tau} \right) \Bigg\|^2_F \Bigg] \notag\\
&\hspace{0.4cm}+\gamma^2  \sum\limits_{j=h\tau}^{k-1} (1 + \beta_1)\ldots (1 + \beta_{k-j-1}) (1 + \beta^{-1}_{k-j}) \E \Bigg[\left\|\Phi^{j+1/3} -\bar \Phi^{j+1/3}  \right\|^2_F \Bigg] \notag \\
&\hspace{0.4cm}+\gamma^2 \sum\limits_{j=h\tau}^{k-1} (1 + \beta_1)\ldots (1 + \beta_{k-1-j}) \E \Bigg[\left\|G^{j+1/3} - \Phi^{j+1/3} - \bar G^{j+1/3} +\bar \Phi^{j+1/3}  \right\|^2_F \Bigg] . \label{eq:err_estim_proof_1}
\end{align}
Setting $\beta_i = \frac{1}{\alpha - i}$, where $\alpha \geq 4 \tau$, gives, for all $i = 0, \ldots, (k-1 - h\tau)$, that 
$$ (1 + \beta_1)(1 + \beta_2)\ldots (1 + \beta_{i}) = \frac{\alpha}{\alpha - i}.$$ 
By the definition of $h$, we have $k-1 - h\tau \leq 2\tau$. Hence, for all $i = 0, \ldots, (k-1 - h\tau)$,
$$(1 + \beta_1)(1 + \beta_2)\ldots (1 + \beta_{i}) \leq (1 + \beta_1)(1 + \beta_2)\ldots (1 + \beta_{k-1 - h\tau}) \leq \frac{\alpha}{\alpha - 2 \tau} \leq 2.$$ 
Moreover, $1 + \beta_i^{-1} \leq \alpha$. Substituting these estimates into \eqref{eq:err_estim_proof_1}, we obtain
\begin{align*}
M\cdot \E\left[\text{Err}(k)\right] &\leq \frac{\alpha}{\alpha - 2 \tau} \E \Bigg[\Bigg\|Z^{h\tau} \prod\limits_{i = h\tau}^{k-1} W^i - \bar Z^{h\tau } - \left(  \bar Z^{h\tau} \prod\limits_{i = h\tau}^{k-1} W^i - \bar Z^{h\tau} \right) \Bigg\|^2_F \Bigg] \\
&\hspace{0.4cm}+2\gamma^2 \alpha \sum\limits_{j=h\tau}^{k-1} \E \Bigg[\left\|\Phi^{j+1/3} -\bar \Phi^{j+1/3}  \right\|^2_F \Bigg] \\
&\hspace{0.4cm}+2\gamma^2 \sum\limits_{j=h\tau}^{k-1} \E \Bigg[\left\|G^{j+1/3} - \Phi^{j+1/3} - \bar G^{j+1/3} +\bar \Phi^{j+1/3}  \right\|^2_F \Bigg] .
\end{align*}
Choosing $\alpha = 4\tau \left(1 + \frac{2}{p}\right)$, we get
\begin{align*}
M\cdot \E\left[\text{Err}(k)\right] &\leq \left( 1 + \frac{1}{1 + \frac{4}{p}} \right) \E \Bigg[\Bigg\|Z^{h\tau} \prod\limits_{i = h\tau}^{k-1} W^i - \bar Z^{h\tau } - \left(  \bar Z^{h\tau} \prod\limits_{i = h\tau}^{k-1} W^i - \bar Z^{h\tau} \right) \Bigg\|^2_F \Bigg] \\
&\hspace{0.4cm}+\frac{24\gamma^2 \tau}{p} \sum\limits_{j=h\tau}^{k-1} \E \Bigg[\left\|\Phi^{j+1/3} -\bar \Phi^{j+1/3}  \right\|^2_F \Bigg] \\
&\hspace{0.4cm}+2\gamma^2 \sum\limits_{j=h\tau}^{k-1} \E \Bigg[\left\|G^{j+1/3} - \Phi^{j+1/3} - \bar G^{j+1/3} +\bar \Phi^{j+1/3}  \right\|^2_F \Bigg] .
\end{align*}
Noticing that, for a matrix $A \in \R^{n \times M}$ with columns $A_i$, $\|A - \bar A \|^2_F = \sum_{i=1}^M\|A_i - \bar A_i\|^2  \leq \sum_{i=1}^M\|A_i\|^2 =  \|A \|^2_F$, we further obtain
\begin{align}
\label{temp11}
M\cdot \E\left[\text{Err}(k)\right] &\leq \left( 1 + \frac{1}{1 + \frac{4}{p}} \right) \E \Bigg[\Bigg\|Z^{h\tau} \prod\limits_{i = h\tau}^{k-1} W^i - \bar Z^{h\tau }\Bigg\|^2_F \Bigg] \notag\\
&\hspace{0.4cm}+\frac{24\gamma^2 \tau}{p} \sum\limits_{j=h\tau}^{k-1} \E \Bigg[\left\|\Phi^{j+1/3} -\bar \Phi^{j+1/3}  \right\|^2_F \Bigg] \notag\\
&\hspace{0.4cm}+2\gamma^2 \sum\limits_{j=h\tau}^{k-1} \E \Bigg[\left\|G^{j+1/3} - \Phi^{j+1/3} - \bar G^{j+1/3} +\bar \Phi^{j+1/3}  \right\|^2_F \Bigg] \notag\\
&\overset{\eqref{as5}}{\leq}
\left( 1 - p\right) \left( 1 + \frac{1}{1 + \frac{4}{p}} \right) \E \Bigg[\Bigg\|Z^{h\tau} - \bar Z^{h\tau }\Bigg\|^2_F \Bigg] \notag\\
&\hspace{0.4cm}+\frac{24\gamma^2 \tau}{p} \sum\limits_{j=h\tau}^{k-1} \E \Bigg[\left\|\Phi^{j+1/3} -\bar \Phi^{j+1/3}  \right\|^2_F \Bigg] \notag\\
&\hspace{0.4cm}+2\gamma^2 \sum\limits_{j=h\tau}^{k-1} \E \Bigg[\left\|G^{j+1/3} - \Phi^{j+1/3} - \bar G^{j+1/3} +\bar \Phi^{j+1/3}  \right\|^2_F \Bigg] . 
\end{align}
It is easy to see that $\left( 1 - p\right) \left( 1 + \frac{1}{1 + \frac{4}{p}} \right) \leq \left( 1 - p\right) \left( 1 +  \frac{p}{4} \right) \leq \left( 1 - \frac{3p}{4}\right)$.
It remains to estimate the last two terms in the r.h.s. of \eqref{temp11}. For the last but one term, we have
\begin{align*}
 \E &\left[\left\| \Phi^{j+1/3} - \bar \Phi^{j+1/3} \right\|^2_F\right] =  \sum\limits_{m=1}^M  \left[\E\left\|F_m(z_m^{j+1/3})-\frac{1}{M} \sum\limits_{i=1}^M F_i(z_i^{j+1/3}) \right\|^2 \right] \\
 &\overset{\eqref{eq:sqs}}{\leq} 3\sum\limits_{m=1}^M \Bigg[ \E\left\| F_m(z_m^{j+1/3}) - F_m(\bar z^{j+1/3})\right\|^2 + \E\left\|F_m(\bar z^{j+1/3}) -  \frac{1}{M} \sum\limits_{i=1}^M F_i(\bar z^{j+1/3})\right\|^2 \\
 &\hspace{0.4cm}+ \E\left\|\frac{1}{M} \sum\limits_{i=1}^M F_i(\bar z^{j+1/3})- \frac{1}{M} \sum\limits_{i=1}^M F_i(z_i^{j+1/3})\right\|^2 \Bigg] \\
 &\overset{\eqref{as4}}{\leq} 3 \sum\limits_{m=1}^M \Bigg[ D^2  + \E\left\|\frac{1}{M} \sum\limits_{i=1}^M F_i(\bar z^{j+1/3})- \frac{1}{M} \sum\limits_{i=1}^M F_i(z_i^{j+1/3})\right\|^2 \\
 &\hspace{0.4cm}+ \E\left\| F_m(z_m^{j+1/3}) - F_m(\bar z^{j+1/3})\right\|^2 \Bigg] \\
 &\overset{\eqref{as1}}{\leq} 6ML^2 \E\left[\text{Err}(j+1/3)\right] + 3 MD^2.
  \end{align*}
For the last term, we have
\begin{align*}
 \E &\left[\left\| G^{j+1/3} - \Phi^{j+1/3} - \bar G^{j+1/3} +\bar \Phi^{j+1/3}  \right\|^2_F\right] \\
 &=  \sum\limits_{m=1}^M  \left[\E\left\|F_m(z_m^{j+1/3}, \xi_m^{j+1/3}) - F_m(z_m^{j+1/3})-\frac{1}{M} \sum\limits_{i=1}^M \left( F_i(z_i^{j+1/3}, \xi_i^{j+1/3}) - F_i(z_i^{j+1/3}) \right) \right\|^2 \right] \\
 &\overset{\eqref{eq:sqs}}{\leq} 2\sum\limits_{m=1}^M \Bigg[ \E\left\| F_m(z_m^{j+1/3}, \xi_m^{j+1/3}) - F_m(z_m^{j+1/3}) \right\|^2 \\
 &\hspace{0.4cm}+ \E\left\|\frac{1}{M} \sum\limits_{i=1}^M \left( F_i(z_i^{j+1/3}, \xi_i^{j+1/3}) - F_i(z_i^{j+1/3}) \right) \right\|^2 \Bigg]\\
 &\overset{\eqref{as3}}{\leq} 4 M\sigma^2.
  \end{align*}
Substituting the last two bounds into \eqref{temp11}, we obtain \eqref{temp4}:
\begin{align*}
\E\left[\text{Err}(k)\right] &\leq 
\left( 1 - \frac{3p}{4}\right) \E [\text{Err}({h\tau}) ] + \frac{144 \gamma^2 L^2 \tau}{p} \sum\limits_{j=h\tau}^{k-1} \E\left[\text{Err}(j+1/3)\right] 
+ \left( \frac{72 D^2 \tau}{p} +  8 \sigma^2\right) \sum\limits_{j=h\tau}^{k-1} \gamma^2 .
\end{align*}
The estimate for $\E\left[\text{Err}(k+1/3)\right]$ is obtained in a similar way. Indeed, it is sufficient to note that $M\E\left[\text{Err}(k+1/3)\right] =\EE \|Z^{k} - \gamma G^{k} - \bar Z^{k} + \gamma \bar G^{k} \|^2_F$. Then, in the proof, we take $\alpha =  4\tau \left(1 + \frac{2}{p}\right) - 1$ and use also $\beta_0 = \frac{1}{\alpha}$ for the term associated with $G^{k} - \bar G^{k}$. In this way, we obtain
$(1 + \beta_0)(1 + \beta_1)(1 + \beta_2)\ldots (1 + \beta_{i}) \leq (1 + \beta_0)(1 + \beta_1)(1 + \beta_2)\ldots (1 + \beta_{k-1 - h\tau}) \leq \frac{\alpha + 1}{\alpha - 2 \tau} \leq 3$, $(1 + \beta_i^{-1}) \leq \alpha + 1$. This gives us the final bound \eqref{temp5}:
\begin{align*}
\E\left[\text{Err}(k+1/3)\right] &\leq 
\left( 1 - \frac{3p}{4}\right) \E [\text{Err}({h\tau}) ] + \frac{216 \gamma^2 L^2 \tau}{p} \sum\limits_{j=h\tau}^{k-1} \E\left[\text{Err}(j+1/3)\right] + \frac{216 \gamma^2 L^2 \tau}{p} \E\left[\text{Err}(k)\right] 
 \\
&\hspace{0.4cm} + \left( \frac{108 D^2 \tau}{p} +  12 \sigma^2\right) \sum\limits_{j=h\tau}^{k-1} \gamma^2 
+ \left( \frac{108 D^2 \tau}{p} +  12 \sigma^2\right) \gamma^2.
\end{align*}
\EndProof
We now notice that the r.h.s. of \eqref{temp4} and \eqref{temp5} involve the terms $\sum\limits_{j=h\tau}^{k-1} \E\left[\text{\rm Err}(j+1/3)\right]$. Thus, in order to resolve the recurrences in \eqref{temp4} and \eqref{temp5}, we need also the bounds for $\E\left[\text{\rm Err}(j+1/3)\right]$ for all $h\tau \leq j \leq k-1$, where $h = \lfloor k / \tau \rfloor - 1$. If $ (h+1) \tau \leq j \leq k-1$, then we can use the same bounds \eqref{temp4} and \eqref{temp5} changing $k$ to $j$ since for such values of $j$ we have $\lfloor j / \tau \rfloor - 1=h$. Thus, it remains to consider such $j$ that $ h\tau \leq j < (h+1) \tau $. This is done in the second technical lemma of this subsection. 
\begin{lemma} \label{l5} Under Assumptions \ref{a1}, \ref{a3}, \ref{a4}, \ref{a5}, for $(h+1) \tau > j \geq h\tau$ with $h = \lfloor k / \tau \rfloor - 1$, it holds that 
 \begin{align}
 \label{temp9}
\E\left[\text{\rm Err}(j)\right] &\leq 
\left( 1 + \frac{p}{4}\right) \E [\text{\rm Err}({h\tau}) ] + \frac{144 \gamma^2 L^2 \tau}{p} \sum\limits_{i=h\tau}^{j-1} \E\left[\text{\rm Err}(i+1/3)\right] \notag\\
&\hspace{0.4cm}+ \left( \frac{72 D^2 \tau}{p} +  8 \sigma^2\right) \sum\limits_{i=h\tau}^{j-1} \gamma^2, \\
 \label{temp10}
\E\left[\text{\rm Err}(j+1/3)\right] &\leq 
\left( 1 + \frac{p}{4}\right) \E [\text{\rm Err}({h\tau}) ] + \frac{216 \gamma^2 L^2 \tau}{p} \sum\limits_{i=h\tau}^{j-1} \E\left[\text{\rm Err}(i+1/3)\right] + \frac{216 \gamma^2 L^2 \tau}{p} \E\left[\text{\rm Err}(j)\right] 
 \notag\\
&\hspace{0.4cm} + \left( \frac{108 D^2 \tau}{p} +  12 \sigma^2\right) \sum\limits_{i=h\tau}^{j-1} \gamma^2 
+ \left( \frac{108 D^2 \tau}{p} +  12 \sigma^2\right) \gamma^2.
\end{align}
 \end{lemma}
  
\textbf{Proof:} 
The proof follows the same lines as the proof of Lemma \ref{l4} until \eqref{temp11}, which needs to be modified since in the situation of the current Lemma, we can not use \eqref{as5} for small $j$'s. The modification is as follows:
\begin{align*}
M\cdot \E\left[\text{Err}(j)\right] 
&{\leq}
\left( 1 + \frac{1}{1 + \frac{4}{p}} \right) \E \Bigg[\Bigg\|\left(Z^{h\tau} - \bar Z^{h\tau }\right) \prod\limits_{i = h\tau}^{j-1} W^i \Bigg\|^2_F \Bigg] \notag\\
&\hspace{0.4cm}+\frac{24\gamma^2 \tau}{p} \sum\limits_{i=h\tau}^{j-1} \E \Bigg[\left\|\Phi^{i+1/3} -\bar \Phi^{i+1/3}  \right\|^2_F \Bigg] \notag\\
&\hspace{0.4cm}+2\gamma^2 \sum\limits_{i=h\tau}^{j-1} \E \Bigg[\left\|G^{i+1/3} - \Phi^{i+1/3} - \bar G^{i+1/3} +\bar \Phi^{i+1/3}  \right\|^2_F \Bigg]\\
&\overset{\eqref{eq:frob}}{\leq}
\left( 1 + \frac{1}{1 + \frac{4}{p}} \right) \E \Bigg[\Bigg\|Z^{h\tau} - \bar Z^{h\tau } \Bigg\|^2_F \Bigg] \notag\\
&\hspace{0.4cm}+\frac{24\gamma^2 \tau}{p} \sum\limits_{i=h\tau}^{j-1} \E \Bigg[\left\|\Phi^{i+1/3} -\bar \Phi^{i+1/3}  \right\|^2_F \Bigg] \notag\\
&\hspace{0.4cm}+2\gamma^2 \sum\limits_{i=h\tau}^{j-1} \E \Bigg[\left\|G^{i+1/3} - \Phi^{i+1/3} - \bar G^{i+1/3} +\bar \Phi^{i+1/3}  \right\|^2_F \Bigg],
\end{align*}
where we also used that $\| W^{i}\|_2 \leq 1$. The rest of the proof is similar to the proof of Lemma \ref{l4}.
\EndProof

\subsubsection{Combining the Building Blocks for the Final Bound}
We are finally ready to combine the building blocks and obtain the convergence rate result for Algorithm \ref{alg4} in the strongly-monotone case.
We combine the refined general bound for the per-iteration progress \eqref{temp41} with the bounds for the consensus error terms \eqref{temp4}, \eqref{temp5}, \eqref{temp9} and \eqref{temp10}.  
We note that, in general, $\E\left[\text{Err}(k+1/3)\right]$ may be smaller than $\E\left[\text{Err}(k)\right]$, but since the r.h.s. of \eqref{temp5} upper bounds the r.h.s. of \eqref{temp4}, we assume, for simplicity, that $\E\left[\text{Err}(k+1/3)\right] \geq \E\left[\text{Err}(k)\right]$. We additionally use that, by the Theorem assumptions, $\gamma \leq \frac{p}{120\tau L}$ and $\gamma \leq \frac{1}{3 L}$, and write the resulting recurrences as follows. 

$\bullet$ Using that $\gamma \leq \frac{1}{3 L}$ and the assumption that $\E\left[\text{Err}(k+1/3)\right] \geq \E\left[\text{Err}(k)\right]$, the recurrence \eqref{temp41} transforms into
\begin{align*}
\E\left[\|\bar z^{k+1} - z^* \|^2\right]  &\leq \left(1 - \frac{\gamma \mu}{2} \right)\E\left[\|\bar z^k - z^* \|^2\right] \notag\\
 &\hspace{0.4cm}+\left(\frac{\gamma L^2}{\mu} + 5\gamma^2 L^2\right) \E\left[\text{Err}(k+1/3) \right] +  5\gamma^2 L^2\E\left[\text{Err}(k) \right] + \frac{10 \gamma^2 \sigma^2}{M} \\
&\leq \left(1 - \frac{\gamma \mu}{2} \right)\E\left[\|\bar z^k - z^* \|^2\right] +\left(\frac{\gamma L^2}{\mu} + 10\gamma^2 L^2\right) \E\left[\text{Err}(k+1/3) \right] + \frac{10 \gamma^2 \sigma^2}{M} \\
 &\leq \left(1 - \frac{\gamma \mu}{2} \right)\E\left[\|\bar z^k - z^* \|^2\right] +\frac{2\gamma L^2}{\mu} \E\left[\text{Err}(k+1/3) \right] + \frac{10 \gamma^2 \sigma^2}{M}. 
  \end{align*}
$\bullet$ Using the assumption that $\E\left[\text{Err}(k+1/3)\right] \geq \E\left[\text{Err}(k)\right]$, the recurrence \eqref{temp5} takes the form
\begin{align*}
\E\left[\text{\rm Err}(k+1/3)\right] &\leq 
\left( 1 - \frac{3p}{4}\right) \E [\text{\rm Err}({h\tau}) ] + \frac{216 \gamma^2 L^2 \tau}{p} \sum\limits_{j=h\tau}^{k-1} \E\left[\text{\rm Err}(j+1/3)\right] + \frac{216 \gamma^2 L^2 \tau}{p} \E\left[\text{\rm Err}(k)\right] 
 \notag\\
&\hspace{0.4cm} + \left( \frac{108 D^2 \tau}{p} +  12 \sigma^2\right) \sum\limits_{j=h\tau}^{k-1} \gamma^2 
+ \left( \frac{108 D^2 \tau}{p} +  12 \sigma^2\right) \gamma^2 \\
&\leq 
\left( 1 - \frac{3p}{4}\right) \E [\text{\rm Err}({h\tau + 1/3}) ] + \frac{216 \gamma^2 L^2 \tau}{p} \sum\limits_{j=h\tau}^{k-1} \E\left[\text{\rm Err}(j+1/3)\right] + \frac{216 \gamma^2 L^2 \tau}{p} \E\left[\text{\rm Err}(k + 1/3)\right] 
 \notag\\
&\hspace{0.4cm} + \left( \frac{108 D^2 \tau}{p} +  12 \sigma^2\right) \sum\limits_{j=h\tau}^{k-1} \gamma^2 
+ \left( \frac{108 D^2 \tau}{p} +  12 \sigma^2\right) \gamma^2.
\end{align*}  
Rearranging, we obtain
\begin{align*}
\left(1 - \frac{216 \gamma^2 L^2 \tau}{p}\right)\E\left[\text{Err}(k+1/3)\right] &\leq 
\left( 1 - \frac{3p}{4}\right) \E [\text{Err}({h\tau} + 1/3) ] + \frac{216 \gamma^2 L^2 \tau}{p}  \sum\limits_{j=h\tau}^{k-1} \E\left[\text{Err}(j+1/3)\right] \\
&\hspace{0.4cm} + \left( \frac{216 D^2 \tau}{p} + 24 \sigma^2\right) \sum\limits_{j=h\tau}^{k-1} \gamma^2.
\end{align*}
Using that $\gamma \leq \frac{p}{120\tau L}$, we get
\begin{align*}
\left(1 - \frac{p}{64}\right)\E\left[\text{Err}(k+1/3)\right] &\leq 
\left( 1 - \frac{3p}{4}\right) \E [\text{Err}({h\tau} + 1/3) ] + \frac{p }{66 \tau}  \sum\limits_{j=h\tau}^{k-1} \E\left[\text{Err}(j+1/3)\right] \\
&\hspace{0.4cm} + \left( \frac{216 D^2 \tau}{p} + 24 \sigma^2\right) \sum\limits_{j=h\tau}^{k-1} \gamma^2.
\end{align*}
Finally, using that $0 < p \leq 1$, and, hence, $\left(1 - \frac{3p}{4}\right)\left(1 - \frac{p}{64}\right)^{-1} \leq 1 - \frac{p}{2}$ and $\left(1 - \frac{p}{64}\right)^{-1} \leq \frac{64}{63} \leq \frac{66}{64}$, we obtain
\begin{align*}
\E\left[\text{Err}(k+1/3)\right] &\leq 
\left( 1 - \frac{p}{2}\right) \E [\text{Err}({h\tau} + 1/3) ] + \frac{p }{64 \tau}  \sum\limits_{j=h\tau}^{k-1} \E\left[\text{Err}(j+1/3)\right] \\
&\hspace{0.4cm} + \left( \frac{225 D^2 \tau}{p} + 25 \sigma^2\right) \sum\limits_{j=h\tau}^{k-1} \gamma^2.
\end{align*}
Note that the above inequality holds for $h=\lfloor k / \tau \rfloor - 1$ since we started with  \eqref{temp5}. In the same way, for $ k-1\geq j \geq (h+1)\tau$, since in this case $h=\lfloor j / \tau \rfloor - 1$, we have the inequality
\begin{align*}
\E\left[\text{Err}(j+1/3)\right] &\leq 
\left( 1 - \frac{p}{2}\right) \E [\text{Err}({h\tau} + 1/3) ] + \frac{p }{64 \tau}  \sum\limits_{i=h\tau}^{j-1} \E\left[\text{Err}(i+1/3)\right] \\
&\hspace{0.4cm} + \left( \frac{225 D^2 \tau}{p} + 25 \sigma^2\right) \sum\limits_{i=h\tau}^{j-1} \gamma^2.
\end{align*}
$\bullet$ 
We repeat the same derivations for \eqref{temp10}. First,
\begin{align*}
\E\left[\text{\rm Err}(j+1/3)\right] &\leq 
\left( 1 + \frac{p}{4}\right) \E [\text{\rm Err}({h\tau}) ] + \frac{216 \gamma^2 L^2 \tau}{p} \sum\limits_{i=h\tau}^{j-1} \E\left[\text{\rm Err}(i+1/3)\right] + \frac{216 \gamma^2 L^2 \tau}{p} \E\left[\text{\rm Err}(j)\right] 
 \notag\\
&\hspace{0.4cm} + \left( \frac{108 D^2 \tau}{p} +  12 \sigma^2\right) \sum\limits_{i=h\tau}^{j-1} \gamma^2 
+ \left( \frac{108 D^2 \tau}{p} +  12 \sigma^2\right) \gamma^2 \\
&\leq 
\left( 1 + \frac{p}{4}\right) \E [\text{\rm Err}({h\tau} + 1/3) ] + \frac{216 \gamma^2 L^2 \tau}{p} \sum\limits_{i=h\tau}^{j-1} \E\left[\text{\rm Err}(i+1/3)\right] + \frac{216 \gamma^2 L^2 \tau}{p} \E\left[\text{\rm Err}(j + 1/3)\right] 
 \notag\\
&\hspace{0.4cm} + \left( \frac{108 D^2 \tau}{p} +  12 \sigma^2\right) \sum\limits_{i=h\tau}^{j-1} \gamma^2 
+ \left( \frac{108 D^2 \tau}{p} +  12 \sigma^2\right) \gamma^2.
\end{align*}
Further, rearranging, we obtain
\begin{align*}
\left(1 - \frac{216 \gamma^2 L^2 \tau}{p}\right) \E\left[\text{\rm Err}(j+1/3)\right] &\leq 
\left( 1 + \frac{p}{4}\right) \E [\text{\rm Err}({h\tau} + 1/3) ] + \frac{216 \gamma^2 L^2 \tau}{p} \sum\limits_{i=h\tau}^{j-1} \E\left[\text{\rm Err}(i+1/3)\right] 
 \notag\\
&\hspace{0.4cm} + \left( \frac{216 D^2 \tau}{p} +  24 \sigma^2\right) \sum\limits_{i=h\tau}^{j-1} \gamma^2.
\end{align*}
Since $\gamma \leq \frac{p}{120\tau L}$, we get
\begin{align*}
\left(1 - \frac{p}{64}\right) \E\left[\text{\rm Err}(j+1/3)\right] &\leq 
\left( 1 + \frac{p}{4}\right) \E [\text{\rm Err}({h\tau} + 1/3) ] + \frac{p}{66 \tau} \sum\limits_{i=h\tau}^{j-1} \E\left[\text{\rm Err}(i+1/3)\right] 
 \notag\\
&\hspace{0.4cm} + \left( \frac{216 D^2 \tau}{p} +  24 \sigma^2\right) \sum\limits_{i=h\tau}^{j-1} \gamma^2.
\end{align*}
Since $0 < p \leq 1$, and, hence, $\left(1 + \frac{p}{4}\right)\left(1 - \frac{p}{64}\right)^{-1} \leq 1 + \frac{p}{2}$ and $\left(1 - \frac{p}{64}\right)^{-1} \leq \frac{64}{63} \leq \frac{66}{64}$, we obtain
\begin{align*}
\E\left[\text{\rm Err}(j+1/3)\right] &\leq 
\left( 1 + \frac{p}{2}\right) \E [\text{\rm Err}({h\tau} + 1/3) ] + \frac{p}{64 \tau} \sum\limits_{i=h\tau}^{j-1} \E\left[\text{\rm Err}(i+1/3)\right] 
 \notag\\
&\hspace{0.4cm} + \left( \frac{225 D^2 \tau}{p} + 25 \sigma^2\right) \sum\limits_{i=h\tau}^{j-1} \gamma^2.
\end{align*}
Summarizing the above three bullet points, we have the following recurrences
\begin{align*}
\E\left[\|\bar z^{k+1} - z^* \|^2\right]
 &\leq \left(1 - \frac{\gamma \mu}{2} \right)\E\left[\|\bar z^k - z^* \|^2\right] +\frac{2\gamma L^2}{\mu} \E\left[\text{Err}(k+1/3) \right] + \frac{10 \gamma^2 \sigma^2}{M}, 
  \end{align*}
\begin{align*}
\E\left[\text{Err}(k+1/3)\right] &\leq 
\left( 1 - \frac{p}{2}\right) \E [\text{Err}({h\tau} + 1/3) ] + \frac{p }{64 \tau}  \sum\limits_{j=h\tau}^{k-1} \E\left[\text{Err}(j+1/3)\right] \\
&\hspace{0.4cm} + \left( \frac{225 D^2 \tau}{p} + 25 \sigma^2\right) \sum\limits_{j=h\tau}^{k-1} \gamma^2, \quad h = \lfloor k / \tau \rfloor - 1,
\end{align*}
\begin{align*}
\E\left[\text{Err}(j+1/3)\right] &\leq 
\left( 1 - \frac{p}{2}\right) \E [\text{Err}({h\tau} + 1/3) ] + \frac{p }{64 \tau}  \sum\limits_{i=h\tau}^{j-1} \E\left[\text{Err}(i+1/3)\right] \\
&\hspace{0.4cm} + \left( \frac{225 D^2 \tau}{p} + 25 \sigma^2\right) \sum\limits_{i=h\tau}^{j-1} \gamma^2,\quad k-1 \geq j \geq (h+1)\tau,
\end{align*}
\begin{align*}
\E\left[\text{\rm Err}(j+1/3)\right] &\leq 
\left( 1 + \frac{p}{2}\right) \E [\text{\rm Err}({h\tau} + 1/3) ] + \frac{p}{64\tau} \sum\limits_{i=h\tau}^{j-1} \E\left[\text{\rm Err}(i+1/3)\right] 
 \notag\\
&\hspace{0.4cm} + \left( \frac{225 D^2 \tau}{p} + 25 \sigma^2\right) \sum\limits_{i=h\tau}^{j-1} \gamma^2,\quad h\tau \leq j < (h+1)\tau.
\end{align*}

To simplify the further derivations, we introduce shortcut notations: $r_k = \E\left[\|\bar z^{k} - z^* \|^2\right]$, $e_k = \E\left[\text{Err}(k+1/3) \right]$, $a = \frac{\mu}{2}$, $A = \frac{225 D^2 \tau}{p} + 25 \sigma^2$, $B = \frac{2L^2}{\mu}$, $C = \frac{10 \sigma^2}{M}$, and $h = \lfloor k / \tau \rfloor - 1$. Then, the previous recurrences can be written as
\begin{align}
\label{temp7}
 r_{k+1} &\leq \left(1 - \gamma a \right) r_k+ \gamma B e_k + \gamma^2 C,  
  \end{align}
\begin{align}
\label{temp13}
e_k &\leq 
\left( 1 - \frac{ p}{2}\right) e_{h\tau} + \frac{p}{64\tau} \sum\limits_{j=h\tau}^{k-1} e_{j} + A \sum\limits_{j=h\tau}^{k-1} \gamma^2, \quad h = \lfloor k / \tau \rfloor - 1,
\end{align}
\begin{align}
\label{temp8}
e_j &\leq 
\left( 1 - \frac{ p}{2}\right) e_{h\tau} + \frac{p}{64\tau} \sum\limits_{i=h\tau}^{j-1} e_{i} + A \sum\limits_{i=h\tau}^{j-1} \gamma^2, \quad k-1 \geq j \geq (h+1)\tau,
\end{align}
\begin{align}
\label{temp12}
e_j &\leq 
\left( 1 + \frac{p}{2}\right) e_{h\tau} + \frac{p}{64\tau} \sum\limits_{i=h\tau}^{j-1} e_{i} + A \sum\limits_{i=h\tau}^{j-1} \gamma^2, \quad h\tau \leq j < (h+1)\tau.
\end{align}
The last two recurrent inequalities can be resolved w.r.t. $e_k$ using the following lemma.
\begin{lemma} \label{l6}
If a non-negative sequence $\{e_k\}_{k\geq 0}$ satisfies \eqref{temp13}, \eqref{temp8}, and \eqref{temp12}
with some constants $0<p \leq 1$, $\tau \geq 1$, $A \geq 0$, then it holds that
\begin{align*}
e_k &\leq \frac{16\gamma^2 A \tau}{ p}.
\end{align*}
\end{lemma}
\textbf{Proof:} We depart from \eqref{temp13} and iteratively substitute all $e_j$ for $j \geq (h+1) \tau$ starting from $k-1$ and finishing with $(h+1) \tau$:
\begin{align*}
e_k &\leq 
\left( 1 - \frac{ p}{2}\right) \cdot \left( 1 + \frac{p}{64\tau}\right) e_{h\tau}  + \frac{p}{64\tau} \left( 1 + \frac{p}{64\tau}\right) \sum\limits_{j=h\tau}^{k-2} e_{j} + A \sum\limits_{j=h\tau}^{k-1} \gamma^2 + \frac{p}{64\tau} \cdot  A \sum\limits_{j=h\tau}^{k-2} \gamma^2 \\
&\leq\left( 1 - \frac{ p}{2}\right) \cdot \left( 1 + \frac{p}{64\tau}\right)^{\tau} e_{h\tau}  + \frac{p}{64\tau} \left( 1 + \frac{p}{64\tau}\right)^{\tau} \sum\limits_{j=h\tau}^{(h+1)\tau - 1} e_{j} \\
&\hspace{0.4cm} + A \left( 1 + \frac{p}{64\tau}\right)^{k - (h+1)\tau} \sum\limits_{j=h\tau}^{(h+1)\tau-1} \gamma^2 + A  \sum\limits_{j=(h+1)\tau}^{k-1}  \left( 1 + \frac{p}{64\tau}\right)^{k - 1 - j} \gamma^2 .
\end{align*}
Next, using \eqref{temp12}, we substitute all $e_j$ for $j$ such that $h\tau \leq j < (h+1)\tau$:
\begin{align*}
e_k &\leq\left( 1 - \frac{ p}{2} + \frac{p}{64\tau} \left( 1 + \frac{p}{2}\right) \right) \cdot \left( 1 + \frac{p}{64\tau}\right)^{\tau} e_{h\tau} +\frac{p}{64\tau} \left( 1 + \frac{p}{64\tau}\right)^{\tau + 1} \sum\limits_{j=h\tau}^{(h+1)\tau - 2} e_{j} \\
&\hspace{0.4cm} + A \left( 1 + \frac{p}{64\tau}\right)^{k - (h+1)\tau + 1} \sum\limits_{j=h\tau}^{(h+1)\tau-2} \gamma^2 + A  \sum\limits_{j=(h+1)\tau - 1}^{k-1}  \left( 1 + \frac{p}{64\tau}\right)^{k - 1 - j} \gamma^2 .
\end{align*}
Since $\frac{p}{64\tau} \left( 1 + \frac{p}{2}\right) \leq \frac{p}{16\tau} \left( 1 - \frac{p}{2}\right)$, we get
\begin{align*}
e_k &\leq\left( 1 - \frac{ p}{2}  \right) \left( 1 + \frac{p}{16 \tau}\right) \left( 1 + \frac{p}{64\tau}\right)^{\tau} e_{h\tau} +\frac{p}{64\tau} \left( 1 + \frac{p}{64\tau}\right)^{\tau + 1} \sum\limits_{j=h\tau}^{(h+1)\tau - 2} e_{j} \\
&\hspace{0.4cm} + A \left( 1 + \frac{p}{64\tau}\right)^{k - (h+1)\tau + 1} \sum\limits_{j=h\tau}^{(h+1)\tau-2} \gamma^2 + A  \sum\limits_{j=(h+1)\tau - 1}^{k-1}  \left( 1 + \frac{p}{64\tau}\right)^{k - 1 - j} \gamma^2 .
\end{align*}
Proceeding in the same way for the rest of $e_j$, we have
\begin{align*}
e_k &\leq\left( 1 - \frac{ p}{2}  \right) \left( 1 + \frac{p}{16\tau}\right)^{2\tau} e_{h\tau} + A  \sum\limits_{j=h\tau}^{k-1}  \left( 1 + \frac{p}{64\tau}\right)^{k - 1 - j} \gamma^2 .
\end{align*}
Noting that $\left( 1 + \frac{p}{64\tau}\right)^{k - 1 - j} \leq \left( 1 + \frac{p}{16\tau}\right)^{2\tau} \leq \exp(p/8) \leq 1 + \frac{p}{4}$ for $p \leq 1$, we further derive

\begin{align*}
e_k &\leq\left( 1 - \frac{ p}{4}  \right) e_{h\tau} + 2A  \sum\limits_{j=h\tau}^{k-1} \gamma^2 \leq \left( 1 - \frac{ p}{4}  \right) e_{h\tau} + 2A\gamma^2 (k - h\tau) \leq \left( 1 - \frac{ p}{4}  \right) e_{h\tau} + 4 A \gamma^2 \tau,
\end{align*}
using also that $h = \lfloor k / \tau \rfloor - 1$. 
It remains to apply the recursion for $e_{h\tau}$ (with $e_0 = 0$), and obtain
\begin{align*}
e_k &\leq\left( 1 - \frac{ p}{4}  \right)^2 e_{(h-1)\tau} + \left( 1 - \frac{ p}{4}  \right) \cdot 4 A \gamma^2 \tau  +  4 A \gamma^2 \tau 
\leq \ldots \\
&\leq 4A \gamma^2 \tau \sum\limits_{j=0}^{h} \left( 1 - \frac{ p}{4}  \right)^j \leq 4A \gamma^2 \tau \sum\limits_{j=0}^{\infty} \left( 1 - \frac{ p}{4}  \right)^j \\
&\leq \frac{16\gamma^2 A\tau}{p}.
\end{align*}

This finishes the proof of the bound for $e_k$.
\EndProof

We proceed with the proof of the main estimate for the strongly-monotone case by substituting the above estimate for $e_k$ into \eqref{temp7}:
\begin{align*}
 r_{k+1} &\leq \left(1 - \gamma a \right) r_k+ \frac{16\gamma^3 A B\tau}{p}  + \gamma^2 C.  
\end{align*}
Running the recursion from $0$ to $K$ gives
\begin{align}
\label{eq:SM_proof_fin_bound}
 r_{K+1} &\leq \left(1 - \gamma a \right)^K r_0+ \frac{16\gamma^2 A B\tau}{a p}  + \frac{\gamma C}{a} \leq \exp\left(- \gamma a K\right) r_0+ \frac{16\gamma^2 A B\tau}{a p}  + \frac{\gamma C}{a}.  
\end{align}
Throughout the proof of this bound we used that $\gamma \leq \frac{p}{120\tau L}$ and $\gamma \leq \frac{1}{3 L}$.
Let us denote $\frac{1}{d} = \frac{p}{120 L\tau}$. 
By the definitions of $\tau$ and $p$, we have that $\frac{1}{d}=\frac{p}{120 L\tau} < \frac{1}{3 L}$, and the choice 
\[
\gamma = \min\left\{\frac{p}{120\tau L} , \frac{2\ln\left( \max\{2, \mu^2 M r_0 K/(40\sigma^2) \} \right)}{\mu K}\right\} = \min\left\{\frac{1}{d} , \frac{\ln\left( \max\{2, a^2 r_0 K/C \} \right)}{aK}\right\}
\]
 in the Theorem assumptions is valid to obtain \eqref{eq:SM_proof_fin_bound}. 
Further, this choice leads to the following convergence rate estimates.

$\bullet$ If $\frac{1}{d} \geq \frac{\ln\left( \max\{2, a^2 r_0 K/C \} \right)}{aK}$ then $\gamma = \frac{\ln\left( \max\{2, a^2 r_0 K/C \} \right)}{aK}$ gives
    \begin{align*}
    \mathcal{O} \left( \exp\left(- a^2 K \cdot \frac{\ln\left( \max\{2, a^2 r_0 K/C \} \right)}{aK} \right) r_0+ \frac{A B\tau}{a p} \cdot \frac{\ln^2\left( \max\{2, a^2 r_0 K/C \} \right)}{a^2 K^2} + \frac{C}{a} \cdot \frac{\ln\left( \max\{2, a^2 r_0 K/C \} \right)}{aK} \right) \\
    = \mathcal{\tilde O} \left( \exp\left(- \ln\left( \max\{2, a^2 r_0 K/C \} \right) \right) r_0+ \frac{A B\tau}{a^3 p K^2}  + \frac{C}{a^2K}\right) \leq \mathcal{\tilde O} \left( \frac{C}{a^2 K} + \frac{A B\tau}{a^3 p K^2}  + \frac{C}{a^2 K}\right).  
    \end{align*}
$\bullet$ If $\frac{1}{d} \leq \frac{\ln\left( \max\{2, a^2 r_0 K/C \} \right)}{aK}$ then $\gamma = \frac{1}{d}$ gives
    \begin{align*}
    \mathcal{ O} \left( \exp\left(- a K \cdot \frac{1}{d}\right) r_0+ \frac{A B\tau}{a p} \cdot \frac{1}{d^2}  + \frac{C}{a} \cdot \frac{1}{d}\right) \hspace{7cm} \\
    \leq \mathcal{ O} \left( \exp\left(- \frac{a K}{d}\right) r_0+ \frac{A B\tau}{a p} \cdot \frac{\ln^2\left( \max\{2, a^2 r_0 K/C \} \right)}{a^2 K^2} + \frac{C}{a} \cdot \frac{\ln\left( \max\{2, a^2 r_0 K/C \} \right)}{aK}\right) \\
    = \mathcal{\tilde O} \left( \exp\left(- \frac{a K}{d}\right) r_0+ \frac{A B\tau}{a^3 p K^2} + \frac{C}{a^2 K}\right).  
    \end{align*}
Thus, our choice of $\gamma$ leads to the desired estimates for the convergence rate
\begin{align*}
 r_{k+1} = \mathcal{\tilde O} \left( \exp\left(- \frac{a K}{d}\right) r_0+ \frac{A B\tau}{a^3 p K^2}  + \frac{C}{a^2 K}\right).  
\end{align*}
Finally, we substitute constants $A, B, C, a, d$ and obtain
\begin{align*}
 r_{k+1} &= \mathcal{\tilde O} \left( \exp\left(- \frac{p \mu K}{240\tau L}\right) r_0+ \frac{\tau^2 D^2 L^2}{p^2 \mu^4  K^2} + \frac{\tau \sigma^2 L^2}{p \mu^4  K^2}  + \frac{\sigma^2}{\mu^2 M K}\right).  
\end{align*}
Using the definition $\Delta=\frac{  \tau}{ p}  \bigl( \frac{D^2 \tau}{p} + \sigma^2 \bigr)$, from the last estimate we obtain the bound in the theorem statement.

This completes the proof in the strongly-monotone case. \EndProof

\subsection{Proof of Theorem \ref{app:th1}, Monotone Case}
The proof partially relies on the general per-iterate estimate in the previous section, which we slightly modify and refine for the purposes of this section.
We first note that in the proof of inequality \eqref{temp1} we can take an arbitrary $z$ instead of $z^*$. Rearranging the terms, we obtain for an arbitrary $z$:
 \begin{eqnarray}
2 \gamma \E\left[\langle \bar g^{k+1/3}, \bar z^{k+1/3} - z \rangle\right] &\leq&  
 \E\left[\|\bar z^k - z \|^2\right] -  \E\left[\|\bar z^{k+1} - z \|^2\right]   \notag\\
 && - \E\left[\|\bar z^{k+1/3} - \bar z^k \|^2\right] + \gamma^2 \E\left[\|\bar g^{k+1/3} - \bar g^{k}\|^2\right].   \label{eq:monot_1lsgd}
 \end{eqnarray}
 Next, we need two bounds: a lower bound for the l.h.s. that relates it with the true operator $F$, and an upper bound for the last term in the r.h.s. that is given by Lemma \ref{l3}.
 
 The lower bound is given by the following result.

 \begin{lemma} 
 \label{Lm:bar_g_lower}
 Let the operator $F$ satisfy Assumptions  \ref{a1}, \ref{a3}. Then, for any fixed $z$, we have
\begin{eqnarray} \label{eq:monot_linear_lower_bound}
\E\left[\langle \bar g^{k+1/3}, \bar z^{k+1/3} - z \rangle\right] &\geq& \E\left[\left\langle  F(\bar z^{k+1/3}), \bar z^{k+1/3} - z \right\rangle\right]\\
&& \hspace{-12em}- \frac{\gamma L^2}{2}\E\left[\| \bar z^{k+1/3} - \bar z^{k} \|^2\right] - \frac{1}{2 \gamma} \E \text{\rm Err}(k+1/3) - L\sqrt{\E \text{\rm Err}(k+1/3)} \sqrt{\E\| \bar z^{k} -  z \|^2}
\end{eqnarray}
\end{lemma}
\textbf{Proof:}
We take into account the independence of all random vectors $\xi^{i} = (\xi^{i}_1, \ldots , \xi^{i}_m)$ and take only the conditional expectation $\E_{\xi^{k+1/3}}$ w.r.t. the vector $\xi^{k+1/3}$ conditioned on all the other randomness:
 \begin{eqnarray*}
  \E\left[\langle \bar g^{k+1/3}, \bar z^{k+1/3} - z \rangle\right] &\stackrel{\eqref{seq1_2}}{=}& \E\left[\left\langle  \frac{1}{M} \sum\limits_{m=1}^M \E_{\xi^{k+1/3}}[F_m(z_m^{k+1/3}, \xi_m^{k+1/3})], \bar z^{k+1/3} - z \right\rangle \right] \\
  &\overset{\eqref{as3}}{=}& \E\left[\left\langle  \frac{1}{M} \sum\limits_{m=1}^M F_m(z_m^{k+1/3}), \bar z^{k+1/3} - z \right\rangle\right] \\
  &=&  \E\left[\left\langle  \frac{1}{M} \sum\limits_{m=1}^M F_m(\bar z^{k+1/3}), \bar z^{k+1/3} - z\right\rangle\right] \\
  &&- \E\left[\left\langle  \frac{1}{M} \sum\limits_{m=1}^M [F_m(\bar z^{k+1/3}) - F_m(z_m^{k+1/3})], \bar z^{k+1/3} - z \right\rangle\right] \\
  &=& \E\left[\left\langle  F(\bar z^{k+1/3}), \bar z^{k+1/3} - z \right\rangle\right] \\
  &&- \E\left[\left\langle  \frac{1}{M} \sum\limits_{m=1}^M [F_m(\bar z^{k+1/3}) - F_m(z_m^{k+1/3})], \bar z^{k+1/3} - z \right\rangle\right]. 
  \end{eqnarray*}

Next, we estimate from below the last term in the r.h.s. Using \eqref{eq:cs} with $c = \gamma L^2$, we obtain
\begin{eqnarray*}
&&\hspace{-2em}- \E\left\langle  \frac{1}{M} \sum\limits_{m=1}^M [F_m(\bar z^{k+1/3}) - F_m(z_m^{k+1/3})], \bar z^{k+1/3} - \bar z^{k} + \bar z^{k}- z \right\rangle \\
  &\geq& - \frac{\gamma L^2}{2}\E \| \bar z^{k+1/3} - \bar z^{k} \|^2 -  \frac{1}{2\gamma L^2} \E \left\| \frac{1}{M} \sum\limits_{m=1}^M [F_m(\bar z^{k+1/3}) - F_m(z_m^{k+1/3})] \right\|^2  \\
&& - \E \left[\left\| \frac{1}{M} \sum\limits_{m=1}^M [F_m(\bar z^{k+1/3}) - F_m(z_m^{k+1/3})] \right\| \| \bar z^{k} -  z \| \right]\\
&\overset{\eqref{as1}}{\geq} & - \frac{\gamma L^2}{2}\E\left[\| \bar z^{k+1/3} - \bar z^{k} \|^2\right]  - \frac{L^2}{2 M \gamma L^2} \E\left[ \sum\limits_{m=1}^M\left\|\bar z^{k+1/3} - z_m^{k+1/3}\right\|^2 \right]  \\
&& -  \frac{L}{M}\E \left[\sum\limits_{m=1}^M\left\|  \bar z^{k+1/3} - z_m^{k+1/3} \right\| \| \bar z^{k} -  z \|\right] \\
  &\overset{\eqref{seq3}}{\geq}& - \frac{\gamma L^2}{2}\E\left[\| \bar z^{k+1/3} - \bar z^{k} \|^2\right] - \frac{1}{2 \gamma} \E\text{Err}(k+1/3) - L\sqrt{\E \text{Err}(k+1/3)} \sqrt{\E\| \bar z^{k} -  z \|^2} ,
  \end{eqnarray*} 
where in the last inequality we used also that
\begin{eqnarray*}
 \E  \left[\frac{1}{M}\sum\limits_{m=1}^M\left\|  \bar z^{k+1/3} - z_m^{k+1/3} \right\| \| \bar z^{k} -  z \|\right] 
 & \leq \sqrt{\E  \left(\frac{1}{M}\sum\limits_{m=1}^M\left\|  \bar z^{k+1/3} - z_m^{k+1/3} \right\|\right)^2} \sqrt{\E\| \bar z^{k} -  z \|^2} 
 \\
 & \leq \sqrt{\frac{1}{M}\E  \sum\limits_{m=1}^M\left\|  \bar z^{k+1/3} - z_m^{k+1/3} \right\|^2} \sqrt{\E\| \bar z^{k} -  z \|^2} 
 \\
 & \overset{\eqref{seq3}}{=} \sqrt{\E \text{Err}(k+1/3)} \sqrt{\E\| \bar z^{k} -  z \|^2}.
\end{eqnarray*} 
Combining the above, we obtain the statement of the Lemma.
\EndProof
 
Combining inequality \eqref{eq:monot_1lsgd} with Lemma \ref{l3} and Lemma \ref{Lm:bar_g_lower}, rearranging the terms, and using the monotonicity of the operator $F$, i.e. \eqref{as2} with $\mu=0$, we obtain, for any $z$
\begin{eqnarray}
2 \gamma \E\left[\left\langle  F(z), \bar z^{k+1/3} - z \right\rangle\right] \notag
&\leq&   
2 \gamma \E\left[\left\langle  F(\bar z^{k+1/3}), \bar z^{k+1/3} - z \right\rangle\right] \notag\\
&\leq&   \E\left[\|\bar z^k - z \|^2\right] -  \E\left[\|\bar z^{k+1} - z \|^2\right]   \notag\\
 && - \E\left[\|\bar z^{k+1/3} - \bar z^k \|^2\right] +  5\gamma^2L^2 \E\left[ \|\bar z^{k+1/3} - \bar z^{k}\|^2\right]  \notag\\
&&+ \frac{10 \sigma^2\gamma^2}{M}   + 5L^2\gamma^2\E\left[\text{Err}(k+1/3) \right] +5L^2\gamma^2\E\left[\text{Err}(k) \right]\notag\\
&& + L^2\gamma^2\E\left[\| \bar z^{k+1/3} - \bar z^{k} \|^2\right] +  \text{Err}(k+1/3) \notag\\
&& + 
2\gamma L \sqrt{\E \text{Err}(k+1/3)} \sqrt{\E\| \bar z^{k} -  z \|^2}\notag\\
&\leq&   \E\left[\|\bar z^k - z \|^2\right] -  \E\left[\|\bar z^{k+1} - z \|^2\right] + \frac{10 \sigma^2\gamma^2}{M} \notag\\
&&+ 5L^2\gamma^2\E\left[\text{Err}(k+1/3) \right] +5L^2\gamma^2\E\left[\text{Err}(k) \right]\notag\\
&& +  \text{Err}(k+1/3) + 2\gamma L \sqrt{\E \text{Err}(k+1/3)} \sqrt{\E\| \bar z^{k} -  z \|^2} \notag\\
&\leq& \E\left[\|\bar z^k - z \|^2\right] -  \E\left[\|\bar z^{k+1} - z \|^2\right] + \frac{10 \sigma^2\gamma^2}{M} + (1+5\gamma^2L^2) \E\left[\text{Err}(k+1/3) \right] \notag\\
&&  + 5\gamma^2L^2\E\left[\text{Err}(k) \right] + 2\gamma L \sqrt{\E \text{Err}(k+1/3)} \sqrt{\E\| \bar z^{k} -  z \|^2},
 \end{eqnarray}
where in the last but one inequality we used that, by the Theorem assumptions, $\gamma \leq \frac{1}{L \sqrt{6}}$.
The above is a refined general bound for the per-iteration progress in the monotone setting. 
Further, by Lemma \ref{l6}, we have, for any $z$,
 \begin{eqnarray}
2 \gamma \E\left[\left\langle  F(z), \bar z^{k+1/3} - z \right\rangle\right] \notag
&\leq&   
\E\left[\|\bar z^k - z \|^2\right] -  \E\left[\|\bar z^{k+1} - z \|^2\right] \notag\\
&&+ \frac{10 \sigma^2\gamma^2}{M} + (1+5\gamma^2L^2) \cdot \frac{8 \gamma^2\tau}{ p} \cdot \left( \frac{225 D^2 \tau}{p} + 25 \sigma^2 \right) \notag\\
&&  + 5\gamma^2L^2 \cdot  \frac{8\gamma^2 \tau}{ p} \cdot \left( \frac{225 D^2 \tau}{p} + 25 \sigma^2 \right) \notag\\
&& + 2 \gamma L \sqrt{ \frac{8 \gamma^2 \tau}{ p} \cdot \left( \frac{225 D^2 \tau}{p} + 25 \sigma^2 \right)} \sqrt{\E\| \bar z^{k} -  z \|^2}\notag\\
&\leq& \E\left[\|\bar z^k - z \|^2\right] -  \E\left[\|\bar z^{k+1} - z \|^2\right] + \frac{10 \sigma^2\gamma^2}{M}\notag\\
&&+ (1+10\gamma^2L^2) \cdot  \frac{8\gamma^2\tau}{ p} \cdot \left( \frac{225 D^2 \tau}{p} + 25 \sigma^2 \right)\notag\\
&& + \gamma L \sqrt{ \frac{32 \gamma^2 \tau}{ p} \cdot \left( \frac{225 D^2 \tau}{p} + 25 \sigma^2 \right)} \sqrt{\E\| \bar z^{k} -  z \|^2}\notag\\
&\leq&\E\left[\|\bar z^k - z \|^2\right] -  \E\left[\|\bar z^{k+1} - z \|^2\right] + \xi \notag\\
&& + \sqrt{\eta}\sqrt{\E\| \bar z^{k} -  z \|^2},
\label{eq:monot_proof_1}
 \end{eqnarray}
where we denote $\Delta:=32 \cdot  \frac{ \tau}{ p} \cdot \left( \frac{225 D^2 \tau}{p} + 25 \sigma^2 \right)$, $\xi:=(1+10\gamma^2L^2) \gamma^2 \Delta + \frac{10 \sigma^2\gamma^2}{M}$, $\eta = \gamma^4 L^2 \Delta$.
Our next goal is to analyze the recurrence \eqref{eq:monot_proof_1} in two situations: bounded iterates and unbounded iterates.

\subsubsection{Unbounded Iterates}
\label{S:app_monot_unbound}
First, we consider the general case when the iterates $\bar z^{k}$ are not assumed to be bounded. We carefully analyze this sequence and prove that this sequence can not go too far from any solution to the variational inequality. This allows us to obtain the final convergence rate bound. Let us denote $r_k(z):=\sqrt{\E\| \bar z^{k} -  z \|^2}$ and let $z^*$ be a solution to the variational inequality. Then, we have
 \begin{eqnarray*}
r_k(z) &\leq& \sqrt{2\E\| \bar z^{k} -  z^* \|^2 + 2 \| z -  z^* \|^2} \leq \sqrt{2\E\| \bar z^{k} -  z^* \|^2} + \sqrt{2 \| z -  z^* \|^2}\\
&=& \sqrt{2}r_k(z^*) + \sqrt{2}\| z -  z^* \|,\\
(r_k(z))^2 &\leq& 2 \E\| \bar z^{k} -  z^* \|^2  +  \| z -  z^* \|^2 = 2(r_k(z^*))^2 + 2\| z -  z^* \|^2.
 \end{eqnarray*}
Thus, from \eqref{eq:monot_proof_1}, we have, for any $z$ and any $k\geq 0$,
 \begin{eqnarray}
2 \gamma \E\left[\left\langle  F(z), \bar z^{k+1/3} - z \right\rangle\right]
&\leq& r_k(z)^2 -  r_{k+1}(z)^2 + \xi + \sqrt{\eta}r_k(z).
\label{eq:monot_proof_2}
 \end{eqnarray}
Summing these inequalities from $k=0$ to $K$, we obtain, for any $z$,
\begin{eqnarray}
2 \gamma (K+1) \E\left[\left\langle  F(z), \widehat{z}^{K} - z \right\rangle\right]
&\leq& r_0(z)^2 + (K+1)\xi + \sqrt{\eta}\sum_{k=0}^Kr_k(z) \notag\\
&\leq& 2 r_0(z^*)^2 + 2 \| z -  z^* \|^2 + (K+1)\xi \notag\\
&&+\sqrt{\eta}\left(\sqrt{2} (K+1)\| z -  z^* \|+ \sqrt{2}\sum_{k=0}^Kr_k(z^*) \right), 
\label{eq:monot_proof_3}
 \end{eqnarray}
where $\widehat{z}^{K}=\frac{1}{K+1}\sum\limits_{k=0}^K\bar z^{k+1/3}$.

Our next goal is to bound from above the expression
\[
r_0(z^*)^2 + (K+1)\xi + \sqrt{2\eta}\sum_{k=0}^Kr_k(z^*).
\]
 Taking $z=z^*$ in \eqref{eq:monot_proof_2} and using the fact that $z^*$ is a solution to the variational inequality, we obtain, for any $k\geq 0$
\begin{eqnarray*}
0 \leq 2 \gamma \E\left[\left\langle  F(z^*), \bar z^{k+1/3} - z^* \right\rangle\right]
&\leq& r_k(z^*)^2 -  r_{k+1}(z^*)^2 + \xi + \sqrt{\eta}r_k(z^*).
 \end{eqnarray*}
Thus, for all $k\geq 0$,
\begin{eqnarray*}
r_{k+1}(z^*)^2 \leq r_k(z^*)^2   + \xi + \sqrt{\eta}r_k(z^*).
 \end{eqnarray*}
Summing these inequalities from $k=0$ to $K$, we obtain
\begin{eqnarray*}
r_{K+1}(z^*)^2 \leq r_0(z^*)^2 + (K+1)\xi + \sqrt{\eta}\sum_{k=0}^{K}r_k(z^*).
\end{eqnarray*}
Note that this inequality holds for arbitrary $K\geq 0$.
We next use the following technical result.
\begin{lemma}[{Lemma B.2 in \cite{gorbunov2018accelerated}}]\label{Lm:technical_lemma_non_acc}
    Let ${\alpha},a_0,\ldots,a_{N-1}, b, R_1,\ldots, R_{N-1}$ be non-negative numbers and
    \begin{equation*}
        R_{l} \leqslant \sqrt{2}\cdot\sqrt{\left(\sum\limits_{k=0}^{l-1}a_k + b\alpha\sum\limits_{k=1}^{l-1}R_k \right)}\quad l=1,\ldots,N.
    \end{equation*}
    Then, for $l=1,\ldots,N$,
    \begin{equation*}
        \sum\limits_{k=0}^{l-1}a_k + b\alpha\sum\limits_{k=1}^{l-1} R_k \leqslant \left(\sqrt{\sum\limits_{k=0}^{l-1}a_k} + \sqrt{2}b\alpha l\right)^2.
    \end{equation*}
\end{lemma}
Choosing $\alpha=1$, $b=\sqrt{\eta}$, $a_0=r_0(z^*)^2 + \xi$, $a_k=\xi$, $k=1,...,K-1$, $R_k=r_{k}(z^*)$, we obtain
\begin{eqnarray*}
r_0(z^*)^2 + (K+1)\xi + \sqrt{\eta}\sum_{k=0}^{K}r_k(z^*)
&\leq& \left( \sqrt{r_0(z^*)^2 + (K+1)\xi} + (K+1)\sqrt{2\eta} \right)^2 \notag\\
&\leq& 2r_0(z^*)^2 + 2(K+1)\xi + 4 (K+1)^2 \eta.
\end{eqnarray*}

Combining the last inequality with \eqref{eq:monot_proof_3}, we obtain
\begin{eqnarray*}
2 \gamma (K+1) \E\left[\left\langle  F(z), \widehat{z}^{K} - z \right\rangle\right]
&\leq&   r_0(z^*)^2 + 2 \| z -  z^* \|^2 +\sqrt{2\eta}(K+1)\| z -  z^* \|   \\
&&+ \left(2r_0(z^*)^2 + 2(K+1)\xi + 4 (K+1)^2 \eta \right)\\
&\leq& 3r_0(z^*)^2 + 2 \| z -  z^* \|^2+ 2(K+1)\xi \\
&&+ \sqrt{2\eta}(K+1)\| z -  z^* \|  + 6 (K+1)^2 \eta \notag .
 \end{eqnarray*}
Dividing both sides of the inequality by $2 \gamma (K+1)$ and using the definitions $\Delta:=32 \cdot  \frac{ \tau}{ p} \cdot \left( \frac{225 D^2 \tau}{p} + 25 \sigma^2 \right)$, $\xi:=(1+10\gamma^2L^2) \gamma^2 \Delta + \frac{10 \sigma^2\gamma^2}{M}$, $\eta := \gamma^4 L^2 \Delta$, we obtain, for all $z\in \mathcal{C}$
\begin{eqnarray*}
 \E\left[\left\langle  F(z), \widehat{z}^{K} - z \right\rangle\right]
 &\leq& 2\frac{\|z^0-z^*\|^2 +  \| z -  z^* \|^2}{\gamma (K+1)}+ \frac{\xi}{\gamma} + \| z -  z^* \|\sqrt{ \frac{\eta}{ 2 \gamma ^2} }+  3(K+1)  \frac{\eta}{\gamma} \notag \\
&\leq& 2\frac{\|z^0-z^*\|^2 +  \| z -  z^* \|^2}{\gamma (K+1)}+ \frac{10 \sigma^2\gamma}{M} + (1+10\gamma^2L^2) \gamma \Delta\notag \\
&&  + \gamma L \| z -  z^* \| \sqrt{\Delta} +  3(K+1) \gamma^3 L^2 \Delta \notag\\
&\leq& \frac{4\Omega_{\mathcal{C}}^2}{\gamma (K+1)} + \frac{10 \sigma^2\gamma}{M} + \gamma \Delta \notag \\
&& + \gamma L \Omega_{\mathcal{C}} \sqrt{\Delta} + 8(K+1) \gamma^3 L^2 \Delta \notag, 
\label{eq:monot_proof_8}
 \end{eqnarray*}
where in the last inequality we used that $z^0,z,z^* \in \mathcal{C}$, $\max_{z,z' \in \mathcal{C}} \|z-z'\|\leq \Omega_{\mathcal{C}}$, and that $K \geq 1$.

Our choice
\[
\gamma = \min\left\{\frac{1}{3L},\left(\frac{2\Omega_{\mathcal{C}}^2M}{5(K+1)\sigma^2}\right)^{\frac{1}{2}}, \left(\frac{\Omega_{\mathcal{C}}^2}{6(K+1)^2L^2\Delta}\right)^{\frac{1}{4}}\right\},
\]
in the Theorem assumptions implies
\[
\frac{4\Omega_{\mathcal{C}}^2}{\gamma (K+1)} = \mathcal{O} \Bigg(\frac{L\Omega_{\mathcal{C}}^2}{K} + \frac{\sigma\Omega_{\mathcal{C}}}{\sqrt{MK}}+\frac{\sqrt{L\Omega_{\mathcal{C}}^3\sqrt{\Delta}}}{\sqrt{K}}\Bigg),
\]
and  we obtain 
\begin{align*}
 \sup_{z \in \mathcal{C}} \E\left[\left\langle  F(z), \widehat{z}^{K} - z \right\rangle\right]
&=& \mathcal{O} \Bigg( \frac{L\Omega_{\mathcal{C}}^2}{K} + \frac{\sigma\Omega_{\mathcal{C}}}{\sqrt{MK}}+\frac{\sqrt{L\Omega_{\mathcal{C}}^3\sqrt{\Delta}}}{\sqrt{K}}+\sqrt{\frac{(\Delta+L^2\Omega_{\mathcal{C}}^2)\Omega_{\mathcal{C}} \sqrt{\Delta}}{KL}}\Bigg).
 \end{align*}



\subsubsection{Bounded Iterates}
Let us now consider the situation under the additional assumption that for all $k$ the iterations of the algorithm satisfy $\|\bar z^k\| \leq \Omega$. In this case, summing \eqref{eq:monot_proof_1} from $k=0$ to $K$, we obtain, for any $z$,
\begin{eqnarray}
2 \gamma (K+1) \E\left[\left\langle  F(z), \widehat{z}^{K} - z \right\rangle\right]
&\leq& \|z^0-z\|^2 + (K+1)\xi + \sqrt{\eta}\sum_{k=0}^K\sqrt{\E\| \bar z^{k} -  z \|^2}  \notag \\
\hspace{-10em}&\leq&  \|z^0-z\|^2 + (K+1)\xi + 2(K+1) \sqrt{\eta}(\Omega +\| z \|).   \notag
\label{eq:monot_proof_10}
 \end{eqnarray}
%
%
Dividing both sides of this inequality by $2 \gamma (K+1)$ and using the definitions $\Delta:=32 \cdot  \frac{ \tau}{ p} \cdot \left( \frac{225 D^2 \tau}{p} + 25 \sigma^2 \right)$, $\xi:=(1+10\gamma^2L^2) \gamma^2 \Delta + \frac{10 \sigma^2\gamma^2}{M}$, $\eta := \gamma^4 L^2 \Delta$, we obtain, for all $z\in \mathcal{C}$
\begin{eqnarray*}
 \E\left[\left\langle  F(z), \widehat{z}^{K} - z \right\rangle\right]
&\leq& \frac{\|z^0-z\|^2 }{2\gamma (K+1)} + \frac{\xi}{\gamma} +(\Omega +\| z \|) \sqrt{ \frac{\eta}{  \gamma ^2} } \notag \\ 
&\leq& \frac{\|z^0-z\|^2 }{2\gamma (K+1)}+ \frac{10 \sigma^2\gamma}{M} + (1+10\gamma^2L^2) \gamma \Delta \\
&&+ (\Omega +\| z \|)\gamma L \sqrt{\Delta}  \notag \\
&\leq& \frac{\Omega_{\mathcal{C}}^2 }{2\gamma (K+1)}+ \frac{10 \sigma^2\gamma}{M} + 10\gamma^3L^2 \Delta \\
&&+ \gamma((\Omega +\Omega_{\mathcal{C}})L \sqrt{\Delta} + \Delta),
\label{eq:monot_proof_11}
 \end{eqnarray*}
where in the last inequality we used that $z^0,z,z^* \in \mathcal{C}$, $\max_{z,z' \in \mathcal{C}} \|z-z'\|\leq \Omega_{\mathcal{C}}$, and that $K \geq 1$.

Our choice 
\[
\gamma = \min\left\{\frac{1}{3L},\left(\frac{\Omega_{\mathcal{C}}^2M}{20(K+1)\sigma^2}\right)^{\frac{1}{2}}, \left(\frac{\Omega_{\mathcal{C}}^2}{60(K+1)^2L^2\Delta}\right)^{\frac{1}{4}},\left(\frac{\Omega_{\mathcal{C}}^2}{(K+1)((\Omega +\Omega_{\mathcal{C}})L \sqrt{\Delta} + \Delta)}\right)^{\frac{1}{2}}\right\},
\]
in the Theorem assumptions implies
\begin{align*}
 \sup_{z \in \mathcal{C}} \E\left[\left\langle  F(z), \widehat{z}^{K} - z \right\rangle\right]
&=& \mathcal{O} \Bigg( \frac{L\Omega_{\mathcal{C}}^2}{K} + \frac{\sigma\Omega_{\mathcal{C}}}{\sqrt{MK}}+\frac{\sqrt{L\Omega_{\mathcal{C}}^3\sqrt{\Delta}}}{K^{3/4}}+\sqrt{\frac{((\Omega +\Omega_{\mathcal{C}})L \sqrt{\Delta} + \Delta)\Omega_{\mathcal{C}}^2 }{K}}\Bigg).
 \end{align*}


\EndProof

\subsection{Proof of Theorem \ref{app:th1}, Non-Monotone Case}
The proof relies on the general per-iterate estimate in Section \ref{S:proof_SM}, which we refine for the purposes of this section.
We start with the same estimate \eqref{temp1}:
\begin{align}
 \E\left[\|\bar z^{k+1} - z^* \|^2\right]  &\leq \E\left[\|\bar z^k - z^* \|^2\right] - \E\left[\|\bar z^{k+1/3} - \bar z^k \|^2\right] \notag\\
 &\hspace{0.4cm}- 2 \gamma \E\left[\langle \bar g^{k+1/3}, \bar z^{k+1/3} - z^* \rangle\right] + \gamma^2 \E\left[\|\bar g^{k+1/3} - \bar g^{k}\|^2\right]. \label{NM_temp1}
 \end{align}
We use the same Lemma \ref{l3} to bound the last term in the r.h.s., and the following counterpart of Lemma \ref{l2} to deal with the last but one term.
\begin{lemma} \label{l7}
Under Assumptions \ref{a1}, \ref{a2}(NM), \ref{a3} it holds that
\begin{align}
 \label{temp15}
 - 2 \gamma \E\left[\langle \bar g^{k+1/3}, \bar z^{k+1/3} - z^* \rangle\right] &\leq  2 \gamma L \sqrt{\E\left[ \| \bar z^{k} - z^*\|^2\right]} \sqrt{\E\left[\text{\rm Err}(k+1/3)\right]} \notag\\
 &\hspace{0.4cm}+ \gamma L \E\left[ \| \bar z^{k+1/3} - \bar z^k\|^2\right] + \gamma L \E\left[ \text{\rm Err}(k+1/3)\right].
  \end{align}
 \end{lemma}
 
 \textbf{Proof:} First of all, we use the independence of all random vectors $\xi^{i} = (\xi^{i}_1, \ldots , \xi^{i}_m)$ and take only the conditional expectation $\E_{\xi^{k+1/3}}$ w.r.t. the vector $\xi^{k+1/3}$, conditioned on all the other randomness. This gives us the following chain of inequalities: 
 \begin{align*}
  - 2 \gamma \E &\left[\langle \bar g^{k+1/3}, \bar z^{k+1/3} - z^* \rangle\right]
  \stackrel{\eqref{seq1_2}}{=} - 2 \gamma \E\left[\left\langle  \frac{1}{M} \sum\limits_{m=1}^M \E_{\xi^{k+1/3}}[F_m(z_m^{k+1/3}, \xi_m^{k+1/3})], \bar z^{k+1/3} - z^* \right\rangle \right] \\
  &\overset{\eqref{as3}}{=} -2 \gamma \E\left[\left\langle  \frac{1}{M} \sum\limits_{m=1}^M F_m(z_m^{k+1/3}), \bar z^{k+1/3} - z^* \right\rangle\right] \\
  &= - 2 \gamma \E\left[\left\langle  \frac{1}{M} \sum\limits_{m=1}^M F_m(\bar z^{k+1/3}), \bar z^{k+1/3} - z^* \right\rangle\right] \\
  &\hspace{0.4cm}+ 2 \gamma \E\left[\left\langle  \frac{1}{M} \sum\limits_{m=1}^M [F_m(\bar z^{k+1/3}) - F_m(z_m^{k+1/3})], \bar z^{k+1/3} - z^* \right\rangle\right] \\
  &= - 2 \gamma \E\left[\left\langle  F(\bar z^{k+1/3}), \bar z^{k+1/3} - z^* \right\rangle\right] \\
  &\hspace{0.4cm}+ 2 \gamma \E\left[\left\langle  \frac{1}{M} \sum\limits_{m=1}^M [F_m(\bar z^{k+1/3}) - F_m(z_m^{k+1/3})], \bar z^{k+1/3} - z^* \right\rangle\right] \\
  &\overset{\eqref{as6}}{\leq}2 \gamma \E\left[\left\langle  \frac{1}{M} \sum\limits_{m=1}^M [F_m(\bar z^{k+1/3}) - F_m(z_m^{k+1/3})], \bar z^{k+1/3} - z^* \right\rangle\right] \\
  &\leq 2 \gamma \E\left[ \| \bar z^{k+1/3} - z^*\| \cdot \left\|\frac{1}{M} \sum\limits_{m=1}^M    F_m(\bar z^{k+1/3}) - F_m(z_m^{k+1/3}) \right\| \right]\\
  &\leq 2 \gamma \E\left[ \| \bar z^{k+1/3} - z^*\| \cdot \frac{1}{M} \sum\limits_{m=1}^M  \left\|   F_m(\bar z^{k+1/3}) - F_m(z_m^{k+1/3}) \right\| \right] \\
  &\overset{\eqref{as1}}{\leq} 2 \gamma L \E\left[ \| \bar z^{k+1/3} - z^*\| \cdot \frac{1}{M} \sum\limits_{m=1}^M  \left\|  z_m^{k+1/3} - \bar z^{k+1/3} \right\| \right] \\
  &\overset{}{\leq} 2 \gamma L \E\left[ \| \bar z^{k} - z^*\| \cdot \frac{1}{M} \sum\limits_{m=1}^M  \left\|  z_m^{k+1/3} - \bar z^{k+1/3} \right\| \right] \\
  &\hspace{0.4cm}+ 2 \gamma L \E\left[ \| \bar z^{k+1/3} - \bar z^k\| \cdot \frac{1}{M} \sum\limits_{m=1}^M  \left\|  z_m^{k+1/3} - \bar z^{k+1/3} \right\| \right]\\
  &\overset{\eqref{eq:cauchy_schwarz_random},\eqref{eq:cs}}{\leq} 2 \gamma L \sqrt{\E\left[ \| \bar z^{k} - z^*\|^2\right]} \cdot \sqrt{\E\left[ \left(\frac{1}{M} \sum\limits_{m=1}^M  \left\|  z_m^{k+1/3} - \bar z^{k+1/3} \right\| \right)^2 \right]} \\
  &\hspace{0.4cm}+\gamma L \E\left[ \| \bar z^{k+1/3} - \bar z^k\|^2\right] + \gamma L \E\left[ \left(\frac{1}{M} \sum\limits_{m=1}^M  \| \bar z^{k+1/3} - z_m^{k+1/3}\|  \right)^2\right].
  \end{align*}
  
It is easy to see that, by convexity of the squared norm,
$$ \E\left[ \left(\frac{1}{M} \sum\limits_{m=1}^M  \| \bar z^{k+1/3} - z_m^{k+1/3}\|  \right)^2\right]\leq \E\left[ \frac{1}{M} \sum\limits_{m=1}^M   \| \bar z^{k+1/3} - z_m^{k+1/3}\|^2 \right] \stackrel{\eqref{seq3}}{=}\E\text{Err}(k+1/3).$$
This completes the proof.
\EndProof

We next move to the refinement of the general per-iterate estimate \eqref{NM_temp1}. Namely, we substitute \eqref{temp15} and \eqref{temp3} into \eqref{NM_temp1}, and obtain the following counterpart of \eqref{temp16}:
\begin{align*}
 \E&\left[\|\bar z^{k+1} - z^* \|^2\right]  \leq \E\left[\|\bar z^k - z^* \|^2\right] - \E\left[\|\bar z^{k+1/3} - \bar z^k \|^2\right] \notag\\
 &\hspace{0.4cm}+ 2 \gamma L \sqrt{\E\left[ \| \bar z^{k} - z^*\|^2\right]} \sqrt{\E\left[\text{Err}(k+1/3)\right]} \\
 &\hspace{0.4cm}+ \gamma L \E\left[ \| \bar z^{k+1/3} - \bar z^k\|^2\right] + \gamma L \E\left[ \text{Err}(k+1/3)\right]\\
 &\hspace{0.4cm}+ \gamma^2 \left( 5L^2 \E\left[ \|\bar z^{k+1/3} - \bar z^{k}\|^2\right] + \frac{10 \sigma^2}{M}  + 5L^2\E\left[\text{Err}(k+1/3) \right] +5L^2\E\left[\text{Err}(k) \right] \right).
  \end{align*}
Since, by the Theorem assumptions,  $\gamma \leq \frac{1}{5L}$, we have
\begin{align}
   \frac{1}{2}\E&\left[\|\bar z^{k+1/3} - \bar z^k \|^2\right] \leq \E\left[\|\bar z^k - z^* \|^2\right] - \E\left[\|\bar z^{k+1} - z^* \|^2\right] \notag\\
 &\hspace{0.4cm}+ 2 \gamma L \sqrt{\E\left[ \| \bar z^{k} - z^*\|^2\right]} \sqrt{\E\left[\text{Err}(k+1/3)\right]} \notag\\
 &\hspace{0.4cm}+ (5\gamma^2L^2 + \gamma L)\E\left[\text{Err}(k+1/3) \right] +5\gamma^2L^2\E\left[\text{Err}(k) \right]  + \frac{10 \gamma^2 \sigma^2}{M}. \label{NM_temp2}
  \end{align}
We next elaborate the term
\begin{align*}
 \E&\left[\|\bar z^{k+1/3} - \bar z^{k}\|^2 \right] \\
 &= \gamma^2\E\left[\left\|\frac{1}{M} \sum\limits_{m=1}^M (F_m(z^k_m, \xi^{k}_m) - F_m(z^k_m) + F_m(z^k_m) - F_m(\bar z^k) + F_m(\bar z^k)) \right\|^2 \right] \\
 &\overset{\eqref{eq:cs1}}{\geq} \frac{\gamma^2}{2} \E\left\| F(\bar z^k) \right\|^2 - \gamma^2\E\left[\left\|\frac{1}{M} \sum\limits_{m=1}^M (F_m(z^k_m, \xi^{k}_m) - F_m(z^k_m) + F_m(z^k_m) - F_m(\bar z^k)) \right\|^2 \right] \\
 &\overset{\eqref{eq:cs1}}{\geq} \frac{\gamma^2}{2} \E\left\| F(\bar z^k) \right\|^2 - 2\gamma^2\E\left[\left\|\frac{1}{M} \sum\limits_{m=1}^M F_m(z^k_m, \xi^{k}_m) - F_m(z^k_m) \right\|^2 \right] - 2\gamma^2\E\left[\left\|\frac{1}{M} \sum\limits_{m=1}^M  F_m(z^k_m) - F_m(\bar z^k) \right\|^2 \right] \\
 &\overset{\eqref{as1}}{\geq} \frac{\gamma^2}{2} \E\left\| F(\bar z^k) \right\|^2 - \frac{2 \gamma^2 \sigma^2}{M} - \frac{2 \gamma^2 L^2}{M} \sum\limits_{m=1}^M \E\left[\left\|  z^k_m-\bar z^k \right\|^2 \right] \\
 &= \frac{\gamma^2}{2} \E\left\| F(\bar z^k) \right\|^2 - \frac{2 \gamma^2 \sigma^2}{M} - 2 \gamma^2 L^2 \E\left[\text{Err}(k)\right].
\end{align*}
Substituting this into \eqref{NM_temp2} gives
\begin{align*}
   \frac{\gamma^2}{4}\E\left[\|F(\bar z^k)\|^2\right] &\leq \E\left[\|\bar z^k - z^* \|^2\right] - \E\left[\|\bar z^{k+1} - z^* \|^2\right] \notag\\
 &\hspace{0.4cm}+ 2 \gamma L \sqrt{\E\left[ \| \bar z^{k} - z^*\|^2\right]} \sqrt{\E\left[\text{Err}(k+1/3)\right]}\\
 &\hspace{0.4cm}+  (\gamma L + 5\gamma^2L^2)\E\left[\text{Err}(k+1/3) \right] +6\gamma^2 L^2\E\left[\text{Err}(k) \right] + \frac{11 \gamma^2 \sigma^2}{M}.
  \end{align*}
    This is a refined general bound for the per-iteration progress in the non-monotone setting.
Applying Lemma \ref{l6} for the consensus error terms, we get
\begin{align*}
   \frac{\gamma^2}{4}\E\left[\|F(\bar z^k)\|^2\right] &\leq \E\left[\|\bar z^k - z^* \|^2\right] - \E\left[\|\bar z^{k+1} - z^* \|^2\right] \notag\\
 &\hspace{0.4cm}+ 2 \gamma L \sqrt{\E\left[ \| \bar z^{k} - z^*\|^2\right]} \sqrt{\frac{8\gamma^2 \tau}{ p} \cdot \left( \frac{225 D^2 \tau}{p} + 25 \sigma^2 \right)}\\
 &\hspace{0.4cm}+ \gamma^2 \left( \frac{11 \sigma^2}{M}  + \frac{8 (\gamma L + 11\gamma^2 L^2) \tau}{ p} \cdot \left( \frac{225 D^2 \tau}{p} + 25 \sigma^2 \right) \right).
  \end{align*}
Finally, summation over all $k$ from $0$ to $K$ and averaging gives the following bound.
\begin{align}
\label{temp202}
\E\left[ \frac{1}{K+1}\sum\limits_{k=0}^K\|F(\bar z^k)\|^2\right] &\leq \frac{4\| z^0 - z^* \|^2}{\gamma^2 (K+1)} - \frac{4\E\left[\| z^{K+1} - z^* \|^2\right]}{\gamma^2 (K+1)}  +  \frac{44\sigma^2}{M}  \notag\\
&\hspace{0.4cm}+ \sqrt{\frac{32 L^2 \tau}{ p} \cdot \left( \frac{225D^2 \tau}{p} + 25\sigma^2 \right)} \cdot \frac{1}{K+1}\sum\limits_{k=0}^K \sqrt{\E\left[ \| \bar z^{k} - z^*\|^2\right]}\notag\\
&\hspace{0.4cm}+ \frac{8(\gamma L + 11\gamma^2 L^2) \tau}{p} \cdot \left( \frac{225 D^2 \tau}{p} + 25\sigma^2  \right).
\end{align}
Our next goal is to analyze this bound in two situations: bounded iterates and unbounded iterates.

\subsubsection{Unbounded Iterates}
\label{S:app_nonmonot_unbound}
First, we consider the general case when the iterates $\bar z^{k}$ are not assumed to be bounded. We carefully analyze this sequence and prove that it can not go too far from the solution to the variational inequality. This allows us to obtain the final convergence rate bound.
To that end, first, we write the following corollary of \eqref{temp202}:
\begin{align*}
\E\left[\| \bar z^{K+1} - z^* \|^2\right] &\leq \| z^0 - z^* \|^2  +  \frac{11\gamma^2 (K+1)\sigma^2}{M}  \notag\\
&\hspace{0.4cm}+ \sqrt{\frac{2 \gamma^4 L^2 \tau}{ p} \cdot \left( \frac{225D^2 \tau}{p} + 25\sigma^2 \right)} \cdot \sum\limits_{k=0}^K \sqrt{\E\left[ \| \bar z^{k} - z^*\|^2\right]}\notag\\
&\hspace{0.4cm}+ \frac{\gamma^2(K+1)(\gamma L + 11\gamma^2 L^2) \tau}{2p} \cdot \left( \frac{225 D^2 \tau}{p} + 25\sigma^2  \right).
\end{align*}
Next, we apply Lemma \ref{Lm:technical_lemma_non_acc} with $R_k = \sqrt{\E\left[ \| \bar z^{k} - z^*\|^2\right]}$, $b = \sqrt{\frac{2 \gamma^4 L^2 \tau}{ p} \cdot \left( \frac{225D^2 \tau}{p} + 25\sigma^2 \right)}$, $a_k = \frac{\gamma^2(\gamma L + 11\gamma^2 L^2) \tau}{2p} \cdot \left( \frac{225 D^2 \tau}{p} + 25\sigma^2  \right) + \frac{11\gamma^2\sigma^2}{M}$, $a_0 = \| z^0 - z^* \|^2 + \frac{\gamma^2(\gamma L + 11\gamma^2 L^2) \tau}{2p} \cdot \left( \frac{225 D^2 \tau}{p} + 25\sigma^2  \right) + \frac{11\gamma^2\sigma^2}{M}$,
and get 
\begin{equation*}
        \sum\limits_{k=0}^{K}a_k + b\sum\limits_{k=1}^{K} R_k \leq \left(\sqrt{\sum\limits_{k=0}^{K}a_k} + \sqrt{2}b (K+1)\right)^2 \leq 2\sum\limits_{k=0}^{K}a_k + 4 b^2 (K+1)^2,
\end{equation*}
which gives
\begin{equation*}
        \sum\limits_{k=1}^{K} R_k \leq \frac{1}{b}\sum\limits_{k=0}^{K}a_k + 4 b (K+1)^2.
\end{equation*}
Substituting this in \eqref{temp202} with the same notation, we have
\begin{align*}
    \frac{\gamma^2 (K+1)}{4}\E\left[ \frac{1}{K+1}\sum\limits_{k=0}^K\|F(\bar z^k)\|^2\right] &\leq \sum\limits_{k=0}^{K}a_k  + b\left( \frac{1}{b}\sum\limits_{k=0}^{K}a_k + 4 b (K+1)^2 \right),
\end{align*}
and, hence, 
\begin{align*}
    \E\left[ \frac{1}{K+1}\sum\limits_{k=0}^K\|F(\bar z^k)\|^2\right] &\leq \frac{8}{\gamma^2 (K+1)}\sum\limits_{k=0}^{K}a_k  + \frac{16 b^2 (K+1)}{\gamma^2}.
\end{align*}
Finally, we get
\begin{align*}
    \E\left[ \frac{1}{K+1}\sum\limits_{k=0}^K\|F(\bar z^k)\|^2\right] &= \mathcal{O} \Bigg(\frac{\| z^0 - z^* \|^2}{\gamma^2 (K+1)} + \frac{(\gamma L +\gamma^2 L^2) \tau}{p} \cdot \left( \frac{D^2 \tau}{p} + \sigma^2  \right)  \\
    & \hspace{1cm} +  \frac{\sigma^2}{M}  + \frac{(K+1)\gamma^2 L^2 \tau}{ p} \cdot \left( \frac{D^2 \tau}{p} + \sigma^2 \right)\Bigg).
\end{align*}
As before, we denote $\Delta:=32 \cdot  \frac{ \tau}{ p} \cdot \left( \frac{225 D^2 \tau}{p} + 25 \sigma^2 \right)$.
Our choice $\gamma = \min\left\{\frac{1}{5L}, \left(\frac{\|z^0 - z^* \|^2}{(K+1)^2L^2 \Delta}\right)^{1/4}\right\}$ in the Theorem assumptions further implies
\begin{align*}
    \E\left[ \frac{1}{K+1}\sum\limits_{k=0}^K\|F(\bar z^k)\|^2\right] &= \mathcal{O} \Bigg(\frac{L^2\| z^0 - z^* \|^2}{K} + L\| z^0 - z^* \|\sqrt{\Delta}   \\
    & \hspace{1cm} +  \frac{\sigma^2}{M}  + \frac{\sqrt{L\| z^0 - z^* \|\Delta^{3/4}}}{\sqrt{K}} \Bigg).
\end{align*}

\subsubsection{Bounded iterates}
Under the additional assumption that $\|z^*\| \leq \Omega$ and $\|\bar z^k\| \leq \Omega$, using \eqref{temp202}, we obtain
\begin{align*}
\E\left[ \frac{1}{K+1}\sum\limits_{k=0}^K\|F(\bar z^k)\|^2\right] &= \mathcal{O} \Bigg(\frac{\| z^0 - z^* \|^2}{\gamma^2 (K+1)} + \frac{(\gamma L + \gamma^2 L^2) \tau}{p} \cdot \left( \frac{D^2 \tau}{p} + \sigma^2  \right)   \\
&\hspace{1cm}+ \frac{\sigma^2}{M} + \sqrt{\frac{L^2 \Omega^2 \tau}{ p} \cdot \left( \frac{D^2 \tau}{p} + \sigma^2 \right)} \Bigg).
\end{align*}
Our choice $\gamma = \min\left\{\frac{1}{5L}, \left(\frac{\Omega^2}{(K+1)L \Delta}\right)^{1/3}\right\}$ in the Theorem assumptions further implies
\begin{align*}
\E\left[ \frac{1}{K+1}\sum\limits_{k=0}^K\|F(\bar z^k)\|^2\right] &= \mathcal{O} \Bigg(\frac{L^2 \Omega^2}{K} + \frac{\sigma^2}{M} + \frac{(L \Omega \Delta)^{2/3}}{K^{1/3}}   + L \Omega \sqrt{\Delta} \Bigg).
\end{align*}

\EndProof


\section{Anytime Convergence via a Restart Technique} \label{sec:restart}
In this section, we propose a simple procedure that gives our method more flexibility by avoiding the fixed budget $K$ for the number of iterations that needs to be set before the start of the method.
We start with a generic interpretation of Theorem \ref{app:th1} that combines all the cases considered in the theorem. In all the cases, there is some optimality measure $\rho(K)$, e.g., $\rho(K)=\E\left[\|\bar z^{K+1} - z^* \|^2\right]$ in the strongly-monotone case. Further, in all the cases there is some function $\Xi(K)$ which bounds $\rho(K)$ from above after $K$ iterations.
Theorem \ref{app:th1} states that if we fix the budget of $K$ iterations and set the stepsize $\gamma(K)$, then after $K$ iterations it is guaranteed that $\rho(K) \leq \Xi(K)$. Let us refer to the iterations of Algorithm \ref{alg4} as basic iterations. We organize the restart procedure as follows. We construct a sequence of the budgets $K_t=2^t$ for $t\geq 0$. For each restart $t$ we set the stepsize $\gamma(K_t)$, run the algorithm for $K_t$ basic iterations and use the obtained point as a warm-start for the next restart. We can also use the same starting point for all the restarts.

Let us now assume that the algorithm has made $N$ basic iterations. This means that it made at least $T = \lfloor \log_2(N+1)\rfloor-1$ restarts. Since at the end of the last restart it made $K_T$ basic iterations with the stepsize $\gamma(K_T)$, we obtain, by Theorem \ref{app:th1}, that we guarantee that
\[
\rho(K_T) \leq \Xi(K_T) = \Xi\left(2^{T} \right) = \Xi\left(2^{\lfloor \log_2(N+1)\rfloor-1} \right) = \Xi\left(\mathcal{O}(N) \right).
\]
Since in all the cases in Theorem \ref{app:th1}, we have that the dependence of $\Xi(K)$ on $K$ is either exponential or polynomial, we obtain that $\rho(K_T) =\Xi\left(\mathcal{O}(N) \right)= \mathcal{O}\left(\Xi\left(N \right)\right)$. Thus, we have obtained an anytime-convergent algorithm with the convergence rates, up to constant factors, similar to that of Algorithm \ref{alg4}. This algorithm does not require to fix the number of basic steps $K$ in advance.
\end{document}